\RequirePackage{fix-cm}
\documentclass[smallextended,envcountsect]{svjour3}       
\smartqed  
\usepackage{graphicx}
\usepackage{amssymb}
\usepackage{amsfonts}
\usepackage{amsmath}
\usepackage{graphics}
\usepackage{csquotes}
\usepackage[caption=false]{subfig}
\usepackage[pdfencoding=auto, psdextra]{hyperref}
\hypersetup{colorlinks=true,linkcolor=blue, anchorcolor=red, citecolor=blue, filecolor = red, urlcolor = blue,
            pdfauthor=author}
\usepackage{bookmark}
\usepackage[T1]{fontenc}
\usepackage{array, booktabs, makecell}
\usepackage{siunitx}
\usepackage{xurl}
\usepackage[numbers]{natbib}
\usepackage{enumerate}
\usepackage[misc]{ifsym}
\usepackage[title]{appendix}
\spnewtheorem{rem}{Remark}[section]{\bfseries}{\itshape}
\spnewtheorem{examp}{Example}[section]{\bfseries}{\itshape}
\spnewtheorem{cond}{Condition}{\bfseries}{\itshape}

\journalname{Stoch Environ Res Risk Assess}

\begin{document}

\title{Analysis of Spherical Monofractal and Multifractal Random Fields
}


\author{Nikolai~Leonenko\and Ravindi~Nanayakkara\and{Andriy~Olenko}}


\institute{Nikolai~N.~Leonenko \at
           School of Mathematics, Cardiff University, \\
           Senghennydd Road, Cardiff CF24 4AG, UK \\
           \email{LeonenkoN@Cardiff.ac.uk} \\
           ORCiD: 0000-0003-1932-4091
           \and
           Ravindi~Nanayakkara \at
           Department of Mathematics and Statistics, La Trobe University, \\ Melbourne, VIC 3086, Australia \\
           \email{D.Nanayakkara@latrobe.edu.au} \\
           ORCiD: 0000-0001-7141-6917
           \and
           \Letter \ Andriy~Olenko \at
           Department of Mathematics and Statistics, La Trobe University, \\ Melbourne, VIC 3086, Australia \\
           \email{A.Olenko@latrobe.edu.au} \\
           ORCiD: 0000-0002-0917-7000}

\date{\vspace{-0.2cm}Received: 15 May 2020 / Accepted: date}

\maketitle
\begin{abstract}
The R\'enyi function plays an important role in the analysis of multifractal random fields. For random fields on the sphere, there are three models in the literature where the R\'enyi function is known explicitly. The theoretical part of the article presents multifractal random fields on the sphere and develops specific models where the R\'enyi function can be computed explicitly. For all considered models explicit expressions of their multifractal spectrum are obtained. Properties of the models and dependencies of their characteristics on parameters are investigated. Then these results are applied to the Cosmic Microwave Background Radiation data collected from the Planck mission. The main statistical model used to describe these data in the literature is isotropic Gaussian fields. We present numerical multifractality studies and methodology based on simulating random fields, computing the R\'enyi function and the multifractal spectrum for different scenarios and actual CMB data.
The obtained results can also find numerous potential applications for other geoscience, environmental and directional data.
\end{abstract}

\section{Introduction}
\label{S1:1}

Recent years have witnessed an enormous amount of attention, in the environmental, earth science, biological and astrophysical  literature, on investigating
spherical random fields. Excellent overviews of some novel geostatistics directions and applications can be found in \cite{christakos2017, jeong2017, marinucci2011random, porcu2018} and references therein.  From a statistical point of view,  random fields on Euclidean spaces is a rather well studied area. However, the majority of available results is not directly translatable to manifolds (where a sphere is an obvious first important candidate for investigations) and requires new stochastic models and tools, see, for example, \cite{porcu2019_2,porcu2019, lang2015,malyarenko2012invariant,marinucci2011random}. This research investigates multifractal properties of spherical random fields and provides practical methodology and examples of applications to actual data.

The concept of multifractality initially emerged in the context of physics. B.~Mandelbrot showed the significance of scaling relations in turbulence modelling. Subsequently this concept developed to mathematical models and examining their fine scale characteristics. A multifractal pattern is a type of a fractal pattern that scales with multiple scaling rules in contrast to monofractals that have only scaling rule. A fractal dimension explores the change in characteristics with respect to the change in the scale used. In general, a multifractal scheme is a fractal scheme where its dynamics cannot be explained by a single fractal dimension. More details and references can be found in \cite{harte2001multifractals}.

Multifractal structures are typical in nature. Multifractal models have been extensively used in the fields of geophysics, genomics, image modelling, finance, hydrodynamic turbulence, meteorology, internet traffic, etc., see references in~\cite{ angulo2008, anh2008multifractality}. Multifractal behaviour has been discovered in stochastic processes as well, see~\cite{angulo2015, grahovac2014detecting}. Multifractal products of stochastic processes have been investigated by \cite{mannersalo2002multifractal} with applications of time series in economics and teletraffic. New teletraffic models have been explored and random multifractal measure constructions by considering the stationarity of the processes' increments were proposed. This methodology has been first introduced in the groundwork \cite{Kahan1987} on multiplicative chaos and $T$-martingales of positive type. \cite{jaffard1999multifractal} has shown the multifractal nature of specific L{\'e}vy processes and demonstrated that the multifractal spectra of such processes depicts a linear pattern rather than a concave pattern which has been noticed in the actual teletraffic data. \cite{molchan1996scaling} has studied the multifractal properties such as scaling exponents of the structure function and the R\'enyi function of random cascade measures under various conditions. Multifractal analysis has been an important technique in the examination of singular measures and the multifractal spectrum of the random measures based on self-similar processes~\cite{falconer1994multifractal}. The main methods to construct random multifractal structures are based on stochastic processes, branching processes and binomial cascades~\cite{doukhan2002theory}. The R\'enyi function plays an important role in the analysis of multifractal random fields. There are several scenarios where the R\'enyi function was computed for the one-dimensional case and time-series. However, there are very few multidimensional models where it is given in an explicit form. \cite{leonenko2013} computed the R\'enyi function for three classes of multifractal random fields on the sphere. It showed some major schemes with regard to the R\'enyi function which reveal the multifractality of random fields that are~homogeneous~and~isotropic.

The Cosmic Microwave Background(CMB) is the radiation from the universe since 380,000 years from the Big Bang. This elongated time period is very short compared to the age of the universe which is of 14 billion years. The CMB is an electromagnetic radiation residue from it's earliest stage. The CMB depicts variations which corresponds to different regions and represents the roots for all future formation including the solar system, stars and galaxies. At the beginning, the universe was very hot and dense and formation of atoms was impossible. The atoms were split as electrons and protons. That time the universe constituted of a plasma or ionised gas. Then the universe started to expand and cool down. Thus, it had been possible for the atoms to reconcile. This phenomenon is known as \enquote{Epoch of combination} and since that time photons have been able to move freely escaping from the opaque of the early universe. The first light which eliminated from this process is termed as the~cosmic~microwave~background,~see~\cite{planck2019}.

In 2009, the European Space Agency launched the mission Planck to study the CMB. The frequency range captured by the Planck is much wider and its resolution is higher than that of the previous space mission WMAP. The CMB's slight variations were measured with a high precision, see~\cite{esa2019}. One of the aims of the mission Planck was to verify the standard model of cosmology using this achieved greater resolution and to find out fluctuations from the specified standard model of cosmology. According to the standard model of the CMB, the Universe is homogeneous  and isotropic. This means that almost every part of the universe has very similar properties and that they do not differ based on the direction of the space. However, various research argue that it's not the case \cite{hill2018foreground, kogut1996tests, marinucci2004testing, Minkov_2019, novikov20, starck2004detection}. The motivation of this paper is to develop several multifractal models and the corresponding statistical methodology and use them and other existing models to study whether the Cosmic Microwave Background Radiation data has a multifractal behaviour.

The aim of this paper is to present and study three known multifractal models for random fields defined on the sphere and suggest several simpler models for which the R\'enyi function can be explicitly computed. 

The first novelty of the obtained results is that for all these models, we derive singularity spectrum and study dependence of their R\'enyi functions on the scaling parameter. We provide several plots that illustrate typical multifractal behaviour of the models. Note, that even for the three known models their singularity spectrum was not computed and analysed before. Secondly, in Section~\ref{S1:5} and the example in Section~\ref{S1:7} we demonstrate the direct probability approach that can be employed to check whether assumptions on models' parameters guarantee the form of the R\'enyi function. This approach is less general than the one that is based on martingales for $q \in [1,2]$, see~\cite{leonenko2013, mannersalo2002multifractal}. The advantage is that the proposed methodology is simple and can also be used for $q>2.$ Third, Section~\ref{S1:7} suggests four simple new models, explicitly computes their singularity spectrum and R\'enyi functions and investigates their properties.  Finally, we discuss the methodology of computing the R\'enyi functions and provide various numerical studies of the actual CMB data. 

The proposed models and methodology can find various applications to other spherical data. The obtained results and discussion in the paper provide detailed guidance how the multifractal modelling and analysis can be done for general  spherical data. It could be very useful for various earth,  environmental and image analysis problems. In particular, the recent paper \cite{olenko2} discusses methodology and provides R code for transforming various spherical and directional data to the HEALPix format of the CMB data. Then, the results of  this research can be directly applied. 

The plan of the paper is as follows. Section~\ref{S1:2} provides main notations and definitions related to the theory of random fields. Section~\ref{S1:3} introduces spherical random fields. Section~\ref{S1:4} gives results related to the theory of multifractality and the R\'enyi function. Section~\ref{S1:5} describes the direct probability approach to get the conditions on the limit random measure~$\mu$. Section~\ref{S1:6} provides results for three known models of the R\'enyi function for spherical random fields of the exponential type. Section~\ref{S1:7} proposes new models based on power transformations of Gaussian fields. Section~\ref{S1:8}~presents numerical studies including computing and fitting the empirical R\'enyi functions for CMB data from different sky windows and models. The conclusions and some new problems are given in~Section~\ref{S1:9}. Finally, proofs of all key results can be found in Appendix~\ref{appen}.

\section{Main notations and definitions}
\label{S1:2}

This section presents background materials in the random fields theory and multifractal analysis methodology. Most of the material included in this and next two sections are based on \cite{lang2015, leonenko1999limit, malyarenko2012invariant, mannersalo2002multifractal,  marinucci2011random}.

   Let $S\subset \mathbb{R}^n,$ $n\in{\mathbb N,}$ be a multidimensional set, $\Vert\cdot\Vert$ denote the Euclidean distance in ${\mathbb{R}}^n,$  $s_{n-1}(1) = \{ u \in \mathbb{R}^n: \Vert u \Vert =1 \},$ and $ SO(n)$ be the group of rotations in~${\mathbb{R}}^n.$ The notation~$|\cdot|$ will be used for the~Lebesgue measure~on~$\mathbb{R}^n$.~$\{\cdot\} \overset{d}{=} \{\cdot\}$ will stand for the equality of finite dimensional~distributions.

   The Kronecker delta is a function defined as:
   $$\delta _{i}^{j}={\begin{cases}0,&{\text{if }}i\neq j,\\1,&{\text{if }}i=j.\end{cases}}$$

   For $\nu > -{\frac{1}{2}}$, we use the Bessel function of the first kind of order $\nu$
\[
    J_{\nu}(z) = \sum_{m=0}^{\infty} (-1)^m \left(\frac{z}{2}\right)^{2m+\nu}[m! \Gamma(m+\nu+1)]^{-1}, \quad z > 0,
\]
where $\Gamma(\cdot)$ is the gamma function.

\begin{definition}
A random field is a function $\xi(\omega, x) : \Omega \times S \rightarrow \mathbb{R}^m$ such that $\xi(\omega, x)$ is a random vector for each $x \in S$. For simplicity it will also be denoted by $\xi(x),$ $x \in S$.
\end{definition}

   When $n=1$, $\xi(x)$ is a random process. When $S \subseteq \mathbb{R}^n, \:n>1$, then $\xi(x)$ is termed as a random field. It is called a vector random field for $m > 1$. In this paper, we concentrate on scalar random fields $\xi(x), \: x \in S,$ $n > 1$ and~$m=1$.

   If $\{\xi(x_1),...,\xi(x_N),x_1,...,x_N \in S\}$ is a set of random variables belonging to a Gaussian system for each $N \geq 1$, then $\xi(x), x \in S$, is called Gaussian.

   We assume that all random variables $\xi(x)$ are defined on the same probability space $(\Omega, \mathcal{A}, P)$.

\begin{definition}
A second order random field is a random function $\xi : S \rightarrow L_{2}(\Omega, \mathcal{A}, P)$, $S \subset \mathbb{R}^n$.
\end{definition}

   In other words, the random variables $\xi(x)$, $x \in S$, satisfy $E|\xi(x)|^2 < +\infty$. Thus, a second order random field over $S$ is  a family $\{ \xi(x), x \in S \}$ of square integrable random variables.

\begin{definition}
A second order random field $\xi(x), x \in \mathbb{R}^n$,  is homogeneous (in the wide sense) if its mathematical expectation $m(x) = E[\xi(x)]$  and covariance function $B(x,y) = cov(\xi(x), \xi(y))$ are invariant with respect to the Abelian group $G = (\mathbb{R}^n, +)$ of shifts in $\mathbb{R}^n$, that is $$m(x) = m(x+\tau), \: B(x,y) = B(x+\tau, y+ \tau),$$ for any $x,y,\tau \in \mathbb{R}^n$.
\end{definition}

   That is, for homogeneous random fields $E[\xi(x)] = const$, and the covariance function $B(x,y) = B(x-y)$ depends only on the difference $x-y$.

   The covariance function $B(x-y)$ of a homogeneous random field is a non-negative definite kernel on $\mathbb{R}^n \times \mathbb{R}^n$, that is, for any $r \geq 1$, $x^{(j)} \in \mathbb{R}^n$, $z_j \in~\mathbb{C}$, $j=1,...,r$,
$$\sum_{i,j=1}^{r} B(x^{(i)}-x^{(j)}) z_i \Bar{z_j} \geq 0.$$

   If the covariance function $B(x)$ is continuous at $x=0$, then the field is mean-square continuous for each $x \in \mathbb{R}^n$ and vice versa. 

\begin{definition}
The second order random field $\xi(x)$ is isotropic (in the wide sense) on $\mathbb{R}^n$ if its mathematical expectation and covariance function  are invariant with respect to the group of rotations on the sphere, i.e.
$$m(x) = m(gx), \; B(x, y) =  B(gx, gy),$$
for every  $x,y,\tau \in \mathbb{R}^n$ and $g \in SO(n).$  
\end{definition}

In the following we will be considering real-valued second order random fields. It will be also assumed that $E[\xi(x)]=0$ without loss of generality.
 
If a real-valued second order random field $\xi(x), x \in \mathbb{R}^n$ is homogeneous and isotropic, then its mathematical expectation and the covariance function depend only on the Euclidean distance $\rho_{xy} = \Vert x-y \Vert$ between $x$ and $y$. It means that its mathematical expectation $m(x)$ and covariance function $B(x,y)$ are invariant with respect to shifts, rotations and reflections in~$\mathbb{R}^n.$

\begin{definition}
A stochastic process $\{ X(t), t \geq 0 \}$ is self-similar if for any non-random constant $a > 0$, there exists non-random constant $b > 0$ such that
$\{ X(at) \} \overset{d}{=} \{ b X(t) \}.$
\end{definition}

   For self-similar, continuous at 0 and non-trivial $X(t)$, the constant $b$ must be equal $a^H, \: a>0,$ where $H \geq 0$. Thus,
$\{ X(at) \} \overset{d}{=} \{ a^H X(t) \}.$
   The constant $H$ is known as the Hurst parameter. The process $\{ X(t), t \geq 0 \}$ is called $H$-$ss$ (self-similar) or $H$-$sssi$ (self-similar stationary increments) if its increments are stationary.

   The concept of multifractal processes was motivated by establishing the following scaling rule of self-similar processes.

\begin{definition}
A stochastic process $X(t)$ is multifractal if it holds $\{ X(ct) \} \overset{d}{=} \{ M(c) X(t) \},$
where $M(c)$ is a random variable independent of $X(t)$ for every $c > 0$ and the distribution of $M(c)$ does not depend on $t$.
\end{definition}

   The process is self-similar if $M(c)$ is non-random for every $c > 0$ and $M(c) = c^H$. The scaling factor $M(c)$ satisfies
$\{ M(ab) \} \overset{d}{=} \{ M_1(a) M_2(b) \}$ for every selection of constants $a$ and $b$ and random $M_1$ and $M_2$ that are independent copies of $M$. This establishes the characteristic of the deterministic factor $H$-$ss$ processes $(ab)^H = a^H b^H.$

   Another definition of multifractality is

\begin{definition}
A stochastic process $X(t)$ is multifractal if there exist non-random functions $c(q)$ and $\tau(q)$ such that for all $t, s \in \mathcal{T}, q \in \mathcal{Q},$
\[
    E|X(t) - X(s)|^q = c(q) |t-s|^{\tau(q)},
\]
where $\mathcal{T}$ and $\mathcal{Q}$ are intervals on the real line with positive length and $0 \in \mathcal{T}$.
\end{definition}

   The function $\tau(q)$ is known as the scaling function. The interval $\mathcal{Q}$ may include negative values. Instead of the increments of the process, the definition can also be established on the moments of the process. i.e. $E|X(t)|^{q} = c(q) t^{\tau(q)}$. Above definitions coincide if the increments are stationary. If $\{ X(t) \}$ is $H$-$sssi$, then it holds that $\tau(q) = Hq$.

\section{Spherical random fields}
\label{S1:3}

This section introduces some basic notations of the theory of random fields on a sphere. The sphere is a simplest case of a manifold in $\mathbb{R}^n$. For simplicity, we consider only the case $n=3$.

Let us denote the 3-dimensional unit ball as $B^3 = \{ x \in \mathbb{R}^3: \Vert x \Vert \leq 1 \}$. The spherical surface in $\mathbb{R}^3$ with a given radius $r > 0$ is
$s_2(r) = \{x \in \mathbb{R}^3 : \Vert x \Vert = r\}$,
with the corresponding Lebesgue measure on the sphere
$\sigma_r(du)
    = \sigma_r(d\theta \cdot d\varphi) = r^2 \sin{\theta}d\theta d\varphi$, $(\theta,\varphi) \in s_2(1).$ For two points on $s_2(r)$ we use  $\Theta$ to denote the length of the angle formed between two rays originating at the origin and pointing at these two points. $\Theta$ is called the angular distance between these two points.

A spherical random field $T = \{ T(r,\theta,\varphi) : 0 \leq \theta \leq \pi, 0 \leq \varphi \leq 2\pi, r > 0\}$ is a random function, which is defined on the sphere $s_2(r)$. We deal with a spherical real-valued mean-square continuous random field $T$ with a constant mean  and finite second order moments. 

\begin{definition}
A real-valued second order random field $T(x), x \in s_2(r),$ with $E[T(x)]=0$ is isotropic if $E[T(x_1)T(x_2)] = B(\cos{\Theta}),$ $x_1, x_2 \in s_2(r)$, depends only on the angular distance $\Theta$ between $x_1$ and $x_2$.
\end{definition}

For the considered mean-square continuous isotropic random fields, the covariance function $B(\cos{\Theta})$ is a continuous function on $[0,\pi)$.

An isotropic spherical random field on $s_2(r)$ can be expanded in a Laplace series in the mean-square sense.
\begin{equation}
    T(r,\theta,\varphi) =  \sum_{l=0}^{\infty}\sum_{m=-l}^{l}Y_l^m (\theta,\varphi)a_l^m(r), \label{eq:1a}
\end{equation}
where $\{Y_l^m (\theta,\varphi)\}$ represents the spherical harmonics defined as
$$Y_l^m (\theta,\varphi) = c_l^m\exp{(im \varphi)}P_l^m(\cos{\theta}), \quad l=0,1,..., \; m=0, \pm 1, ..., \pm l,$$
with
$$c_l^m = (-1)^m \left[\frac{2l+1}{4\pi}\frac{(l-m)!}{(l+m)!}\right]^{1/2},$$ the associated Legendre polynomials $P_l^m(\cdot)$ having degree $l$ and order $m$ 
\[
P_{l}^{m}(x)=(-1)^{m}\left(1-x^{2}\right)^{m / 2} \frac{d^{m}}{d x^{m}} P_{l}(x), \]
and the $l$-th Legendre polynomials, see \cite{leonenko2013},
\[P_{l}(x)=\frac{1}{2^{l} l !} \frac{d^{l}}{d x^{l}}\left(x^{2}-1\right)^{l}.\] 

The spherical harmonics have the following properties
$$
\int_{0}^{\pi} \int_{0}^{2 \pi} \overline{Y_l^m(\theta, \varphi)} Y_{l^{\prime}}^{ m^{\prime}}(\theta, \varphi) \sin \theta d \varphi d \theta=\delta_{l}^{l^{\prime}} \delta_{m}^{m^{\prime}},$$

$$\overline{Y_l^m(\theta, \varphi)}=(-1)^{m} Y_l^{(-m)}(\theta, \varphi).$$

   The notation $\Tilde{T}(x) = T(r,\theta,\varphi), x \in \mathbb{R}^3,$ will be used to highlight the random field's dependence on Euclidean coordinates.

   The random coefficients of the Laplace series can be computed as the mean-square stochastic integrals via the inversion arguments as
\begin{equation}
    a_l^m(r)=\int_0^{\pi} \int_0^{2\pi} T(r,\theta,\varphi)\overline{Y_l^m(\theta, \varphi)}r^2 \sin{\theta}d\theta d\varphi = \int_{s_2(1)}\Tilde{T}(ru)\overline{Y_l^m(u)}\sigma_1(du), \label{eq:2a}
\end{equation}
   where $u= {\frac{x}{\Vert x \Vert}} \in s_2(1), \quad r= \Vert x \Vert.$

   The covariance functions $E(T(r,\theta,\varphi)T(r,\theta^{'},\varphi^{'}))$ of the isotropic random fields depend only on the angular distance $\Theta = \Theta_{PQ}$ between the points $P = (r,\theta,\varphi)$ and $Q = (r,\theta^{'},\varphi^{'}).$ For spherical isotropic random fields it possesses
\begin{equation}
    E a_l^{m}(r)\overline{a_{l^{'}}^{{m}^{'}}}(r) = \delta_l^{l^{'}}\delta_m^{m^{'}}C_l(r), \quad E |a_l^m(r)|^2 = C_l(r), \quad m = 0, \pm 1,..., \pm l.\label{eq:3a}
\end{equation}

   The angular power spectrum of the isotropic random field $T(r,\theta,\varphi)$ is defined as the functional series $\{C_0(r), C_1(r),..., C_l(r),...\}.$  

   From \eqref{eq:1a} - \eqref{eq:3a} and the addition theorem for spherical harmonic functions we obtain
\[
    Cov(T(r,\theta,\varphi),T(r,\theta^{'},\varphi^{'})) = \frac{1}{4\pi}\sum_{l=0}^{\infty} (2l+1)  C_l(r) P_l(\cos{\Theta}),
\]
where for every $r > 0$ it holds $\sum_{l=0}^{\infty} (2l+1)  C_l(r) < \infty.$ 

   If $T(r,\theta,\varphi)$ is an isotropic Gaussian field defined on the sphere $s_2(r),$ then the coefficients $a_l^m(r)$ are independent Gaussian random variables that are complex-valued with $E a_l^m(r) = 0.$

 For the homogeneous and isotropic random field $\Tilde{T}(x),$ $x \in \mathbb{R}^3,$ it holds
\[ E a_l^m(r)\overline{a_{l^{'}}^{m^{'}}}(s) =\delta_l^{l^{'}}\delta_m^{m^{'}}C_l(r,s),\ r>0,s>0,
\]
   where
\[    C_l(r,s) = 2\pi^2 \int_0^{\infty} \frac{J_{l+\frac{1}{2}}(\mu r){J_{l+\frac{1}{2}}}(\mu s)}{(\mu r)^{1/2}(\mu s)^{1/2}}G(d\mu), \quad l = 0, 1,2,...,
\]
   and $G(\cdot)$ is a finite measure defined on the Borel sets of $[0, \infty)$ satisfying $$\sigma^2 = Var\{\Tilde{T}(x)\} = \int_0^{\infty} G(d\mu) < \infty, \quad x \in {\mathbb{R}^3}.$$

\section{R\'enyi function and multifractal spectrum}
\label{S1:4}

   This section introduces basic notations and concepts of the multifractal theory and R\'enyi functions. 
   
   The  R\'enyi function which is also known as the index of diversity is used in multifractal analysis to assess the randomness of many natural phenomena. It can be used to detect the multifractal behaviour of a given random process. The R\'enyi function computes how the measure/mass/intensity on a surface varies with the resolution or the block size of an image. That is, it calculates the change in detail of a pattern according to the change in scale. The R\'enyi function characterises the distortion in the mean of a pattern's probability distribution of pixel values. R\'enyi functions of non-fractal and monofractal processes exhibit a flatter curve than ones of  multifractal processes.  R\'enyi functions of multifractal processes typically have quadratic shapes that suggest the presence of different fractal dimensions.

   Consider a random field $\Lambda(x, \omega), \; x \in \mathbb{R}^3, \omega \in \Omega,$ that is measurable, homogeneous and isotropic (HIRF) on the 3-dimensional Euclidean space $\mathbb{R}^3$. It will be called the mother field. For simplicity it will be denoted as $\Lambda(x) = \Lambda(x,\omega).$
 
\begin{cond}\label{cond1}
Let a random field $\Lambda(x),$ $x \in \mathbb{R}^3$, satisfy
\[ E[\Lambda(x)] = 1, \quad \Lambda(x) > 0,\quad Cov(\Lambda(x), \Lambda(y))=R_{\Lambda}(\Vert x-y \Vert)=\sigma_{\Lambda}^2 \rho_{\Lambda}(\Vert x-y \Vert),\]
where $\rho(0)=1 \; \text{and} \; \sigma_{\Lambda}^2 < \infty.$
\end{cond}

Let $\Lambda^{(i)}(x), x \in \mathbb{R}^3, i=0,1,2,...,$ be a sequence of independent copies of the random field $\Lambda(\cdot)$. We consider the re-scaling of $\Lambda(\cdot)$ defined as $\Lambda^{(i)}(b^{i}x),$
where $b > 1$ is a constant called a scaling factor and $b^{i}x$ is the product of a vector $x$ by a scalar $b^{i}$.

   A finite-product field on $B^3$ is defined by
\[
    \Lambda_k(x) = \prod_{i=0}^{k}\Lambda^{(i)}(b^ix), \quad k=1,2,....
\]

   Then one can introduce the random measure $\mu_k(\cdot)$ on the Borel $\sigma$-algebra $\mathcal{B}$ of a unit ball $B^3$ by
$$\mu_k(A) = \int_{y \in A}\Lambda_k(y) dy, \: A \in \mathcal{B}, \: k=0,1,2,....$$

   We denote by $\mu_k \xrightarrow{d} \mu, \quad k \rightarrow \infty,$ the weak almost surely convergence of the measures $\mu_k$ to some measure $\mu.$ It means that for all continuous functions $g(y), y \in B^3,$~it~holds with probability 1 that
$$ \int_{B^3}g(y)\mu_k(dy) \rightarrow \int_{B^3}g(y)\mu(dy), \: k \rightarrow \infty.$$

\begin{rem}\label{muconv}
The weak almost surely convergence of random measures implies that for  a finite or countable family of sets $A_j$ from $\mathcal{B},$  with probability 1,
\[\mu_k(A_j) \to \mu(A_j), \quad k\to +\infty,\]  
for all $j,$ see {\rm \cite{Kahan1987, mannersalo2002multifractal}}. Moreover, it was shown in {\rm \cite{denisov2016limit, leonenko2013, mannersalo2002multifractal}} that for mother random fields with $\rho_{\Lambda}(r)$ possessing an exponentially decaying bound, the random variables $\mu_k(B^3)$ converge to $\mu(B^3)$ in $L_2$ (and hence in $L_q$ for $q\in (0, 2]$) when $k\to +\infty.$ 
In the following, for all models considered in this paper, it will be assumed that $|\rho_\Lambda(r)| \leq Ce^{-\gamma r}$ for some positive constants $C,$ and $\gamma$.
\end{rem}

\begin{definition}
The R\'enyi function of a random measure $\mu$ is a non-random function defined by 
$$T(q) = \liminf \limits_{m \rightarrow \infty} \frac{\log_2 E \sum_l \mu (B_l^{(m)})^{q}}{\log_2 |B_l^{(m)}|},$$
   where $\{ B_l^{(m)}, l = 0,1,...,{2^m -1}, \: m=1,2,...,\}$ denotes the mesh formed by the $m^{th}$ level dyadic decomposition of the unit ball $B^3$.
\end{definition}
   The~key result~about~the~form~of~the~R\'enyi~function~is~the~following theorem. 

\begin{theorem}\textup{\cite{leonenko2013}}\label{theo:4.1} Suppose that Condition~{\rm\ref{cond1}} holds.
\begin{enumerate}[(i)]
    \item Assume that the correlation function $\rho_{\Lambda}(\Vert x-y \Vert) = \rho(r)$ of the field $\Lambda(\cdot)$ satisfies the following condition
    \begin{equation}
        |\rho_{\Lambda}(r)| \leq Ce^{-\gamma r}, \: r > 0, \label{eq:4a}
    \end{equation}
    for some positive constants $C$ and $\gamma$. Then, for the scaling factor $b > \sqrt[3]{1+ \sigma_{\Lambda}^2},$
    the measures $\mu_k \xrightarrow{d} \mu,$ $k \rightarrow \infty,$ on $B^3$.
    
    \item If for some range $q \in Q = [q_{-}, q_{+}],$ both $E^q \Lambda(0) < \infty$ and $E\mu^q(B^3) < \infty,$ then the R\'enyi function $T(q)$ of $\mu$ is given by
    \[
        T(q) = q-1-\frac{1}{3}\log_{b}E\Lambda^q(x), \quad q \in Q.
    \]
\end{enumerate}
\end{theorem}

   Similarly, for spherical random fields on $s_2(1)$, one can introduce an analogous~approach.

\begin{cond}\label{cond2}
Let the random field $\Tilde{\Lambda}(x), \: x \in s_2(1),$ satisfy
$$E\Tilde{\Lambda}(x) = 1, \quad Var\Tilde{\Lambda}(x) = \sigma_{\Tilde{\Lambda}}^2 < \infty, \quad \Tilde{\Lambda}(x) > 0,$$
$$Cov(\Tilde{\Lambda}(\theta, \varphi), \Tilde{\Lambda}(\theta^{'}, \varphi^{'})) 
    = \frac{1}{4\pi}\sum_{l=0}^{\infty} (2l+1)  C_l P_l(\cos{\theta}), \quad \sum_{l=0}^{\infty} (2l+1) C_l < \infty.$$
 \end{cond}   
 
   Let $ \Tilde{\Lambda}^{(i)}(x), x \in s_2(1)$, $i=0,1,2,...$, be a sequence of independent copies of the field $\Tilde{\Lambda}(\cdot)$. Let us use the following spherical coordinate notations for points on $s_2(1):$  $x =(1, \theta, \varphi) \in s_2(1).$ Consider $\Tilde{\Lambda}^{(i)}(b^{i}\times x),$ where $b>1$ is a scaling factor, $b^{i}\times x := (1,b^i \underset{\pi}{\times} \theta, b^i \underset{2\pi}{\times} \varphi) \in s_2(1),$ and the modulus algebra is used to compute the products $b^i \underset{\pi}{\times} \theta $ and $ b^i \underset{2\pi}{\times} \varphi.$

   Define the finite product fields on $s_2(1)$ by
$$ \Tilde{\Lambda}_k(x) = \prod_{i=0}^{k}\Tilde{\Lambda}^{(i)}(b^i\times x), \quad k=1,2,....$$

   Let us introduce the random measure $\mu_k(\cdot)$ on the Borel $\sigma$-algebra $\mathcal{B}$ of $s_2(1)$ as
\begin{equation}
    \mu_k(A) = \int_{A}\Tilde{\Lambda}_k(y) dy, \: k=0,1,2,..., \: A \in \mathcal{B}. \label{eq:5a}
\end{equation}
    
   We denote by $\mu_k \xrightarrow{d} \mu, \: k \rightarrow \infty,$ the weak convergence of the measures $\mu_k$ to some non-degenerate measure $\mu$. It means that for all continuous functions $g(y), y \in s_2(1),$
$$ \int_{s_2(1)}g(y)\mu_k(dy) \rightarrow \int_{s_2(1)}g(y)\mu(dy), \quad k \rightarrow \infty.$$

   The R\'enyi function of the random measure $\mu$ defined on $s_2(1)$ is defined as

\begin{equation}
    T(q) = \liminf \limits_{m \rightarrow \infty} \frac{\log_2 E \sum_l \mu(S_l^{(m)})^{q}}{\log_2 |S_l^{(m)}|},\label{eq:6a}
\end{equation}
   where $\{ S_l^{(m)}, l = 0,1,...,{2^m-1}\}$ is the mesh constructed by $m^{th}$ level dyadic decomposition of the spherical surface of $s_2(1)$.

\begin{theorem}\textup{\cite{leonenko2013}}\label{theo:4.2}
Suppose that Condition~{\rm\ref{cond2}} holds and the isotropic random field $\Tilde{\Lambda}(\cdot)$ is the restriction to the sphere $s_2(1)$ of the HIRF $\Lambda(x), x \in \mathbb{R}^3,$ with the correlation function $\rho_{\Lambda}(\Vert x-y \Vert) = \rho(r)$. Under similar assumptions to Theorem~{\rm \ref{theo:4.1}}, the R\'enyi function $T(q)$ of the limit measure $\mu$ on $s_2(1)$ is given~by
\[T(q) = q-1-\frac{1}{2}\log_{b}E\Lambda^q(t), \quad q \in Q.\]   
\end{theorem}

\begin{rem} If $x$ and $y$ are two locations  on the unit sphere  $s_2(1)$ and  $\Theta$ is the angle between them, then the Euclidean distance between these two points is $2\sin({\Theta}/{2})$, which gives a direct correspondence between the covariance function $\rho_{\Lambda}(\Vert x - y\Vert)$ in the Euclidean space and the covariance function  $\rho(\cos{\Theta}) = \rho_{\Lambda}(2\sin{{\Theta}/{2}})$ on the sphere. Thus, the restriction of the HIRF $\Lambda(x)$  to $s_2(1)$ is an isotropic spherical random field. 
\end{rem}

   The multifractal or singularity spectrum is defined via the Legendre transform~as
\begin{equation}
    f(h) = \inf_{q}(hq-T(q)).\label{eq:7a}
\end{equation}
   and is used to describe local fractal dimensions of random fields.

\section{Conditions on measure  \texorpdfstring{$\mu$}{u}}
\label{S1:5}

The random measure $\mu$ in the previous section was defined as a weak limit of the measures $\mu_k$. Therefore, it would be difficult to check moment conditions on $\mu$ as its probability distribution is not explicitly known. This section provides sufficient conditions on the scaling factor $b$ and the variance $\sigma^2_{\Lambda}$ that guarantee $E\mu^q(B^3)< \infty$. The general method to obtain such conditions for the range $q \in [1, 2]$ uses martingale $L^2$ convergence, see, for example, \cite{mannersalo2002multifractal}. The proof of the main result of this section in Appendix~\ref{appen} demonstrates the direct probability approach, which is more elementary.

\normalsize{\begin{theorem}\label{theo:5.1}
Let the mother field $\Lambda(x) > 0, \: x \in \mathbb{R}^{3},$ satisfy the conditions
$$E\Lambda(x) = 1, \quad Var\Lambda(x) = \sigma_{\Lambda}^2 < +\infty,\quad Cov(\Lambda(x), \Lambda(y)) = \sigma_{\Lambda}^2 \rho_{\Lambda}(\Vert x-y \Vert),$$
$$|\rho_\Lambda(\tau)| \leq Ce^{-\gamma \tau}, \: \tau>0,$$
and the scaling factor $b > \max(\sqrt[3]{1+\sigma_{\Lambda}^2}, \: e^{\frac{\sigma^2_{\Lambda}C}{3}}).$

   Then the measures $\mu_k \xrightarrow{d} \mu, \: k \rightarrow \infty,$ and $E\mu^q(B^3) < + \infty,$ \: for $q \in [1,2].$ 
\end{theorem}}

\begin{rem}\label{rem:5.1}
The direct probability approach can be used to obtain conditions on the mother field that guarantee $E \mu_k^q(B^3) < + \infty,$ \: for $q$ in the range $[1, Q],$ where $Q > 2.$

   For example, using the Lyapunov's inequality for $q \in [1, 4],$ see ~{\rm \cite[p.162]{Loeve},}
$$E \mu_k^q(B^3) \leq (E \mu_k^4(B^3))^{q/4},$$
the conditions on $b$ and $\sigma_{\lambda}^2$ that guarantee $E \mu_k^4(B^3) < + \infty$ are also sufficient for $E \mu^q(B^3) < +\infty,$ $q \in [1, 4]$. Then, it follows from $$E \mu_k^4(B^3) = \int_{B^3}\int_{B^3}\int_{B^3}\int_{B^3} \prod_{i=0}^{k}E\left(\prod_{j=1}^4 \Lambda^{(i)}(y_i b^i)\right)\prod_{j=1}^4 dy_i,$$ that one can impose some additional assumptions on the fourth order moments $E(\prod_{j=1}^4 \Lambda(y_j b^i))$ or cumulants of the mother field $\Lambda(\cdot)$.
We will provide an example of such conditions in Section~{\rm\ref{S1:7}} for Model~{\rm 4.}
\end{rem}

\section{R\'enyi functions of exponential models}
\label{S1:6}

For the random fields on the sphere, there are three models where the R\'enyi function is known explicitly, see \cite{leonenko2013}. These models were obtained for exponential type spherical random fields. This sections introduces these models, derives their singularity spectrum and studies dependence of their R\'enyi functions on the scaling parameter.\\[0mm]

\noindent \textbf{Model 1} 
Let the random field $\Lambda(x)$ be given as $$\Lambda(x) = exp\left \{ Y(x) - \frac{1}{2}\sigma_Y^2 \right\},$$ where $Y(x), \: x \in \mathbb{R}^3,$ is a zero-mean Gaussian, measurable, separable random field with the covariance function $\sigma_Y^2\rho_Y(r), \: \rho_Y(0)=1.$\\[0mm]

   The following result provides the conditions and the explicit form of the R\'enyi function for Model 1.

\begin{theorem}\textup{\cite{leonenko2013}}\label{theo:6.1}
Let for Model {\rm 1} the correlation function satisfy
$$0 < |\rho_Y(r)| \leq Ce^{-\gamma r},\: r > 0,$$
for some positive $C$ and $\gamma$ and $b > exp\{\frac{\sigma_Y^2}{3}\}.$ 

   If $Y(x),$ $x \in s_2(1),$ is a spherical isotropic random field that is a restriction of $Y(x), \: x \in \mathbb{R}^3$, on the sphere $s_2(1),$ then the random measures \eqref{eq:5a} generated by the spherical fields $\Tilde{\Lambda}(x) = exp\left \{ Y(x) - \frac{1}{2}\sigma_Y^2 \right\},$ $x \in s_2(1),$ converge weakly to the random measure $\mu.$ The corresponding R\'enyi function is
\begin{equation}
    T(q) = q\left(1+\frac{\sigma_Y^2}{4\ln b}\right) - q^2\left(\frac{\sigma_Y^2}{4\ln b}\right)-1, \quad q \in [1,2].\label{eq:8a}
\end{equation}
\end{theorem}
\noindent    \textbf{Model 2} Let the random field $\Lambda(x)$ be of the form $$\Lambda(x) = exp\left \{ Z(x) - c_Z \right\}, \quad c_Z = -\ln {\left(1-\frac{1}{\lambda}\right)^{\beta}},$$
where $Z(x), \: x \in \mathbb{R}^3,$ is a gamma-correlated HIRF with the correlation function~$\rho_Z(r).$\\[0mm]

   The field $Z(x)$ has the marginal density
\begin{equation}
    p(u) = \frac{\lambda^{\beta}}{\Gamma (\beta)}u^{\beta-1}e^{-\lambda u}, \quad u, \lambda, \beta \in (0, + \infty), \label{eq:9a}
\end{equation}
 and the bivariate density
\begin{equation}
    p_0(u_1, u_2;\alpha) = \frac{{\left({u_1u_2}/{\alpha}\right)^{\frac{\beta -1}{2}}}}{\Gamma(\beta)(1-\alpha)} \exp\left\{-\frac{u_1+u_2}{1-\alpha}\right\} I_{\beta -1}\left(2\frac{\sqrt{u_1\cdot u_2 \cdot \alpha}}{1-\alpha}\right),\label{eq:10a}
\end{equation}
 where $I_v(z) = \sum_{k=0}^{\infty}\left(\frac{z}{2}\right)^{2k+v} {(k!\Gamma(k+v+1))}^{-1}$ is the modified Bessel function of the first kind, $\alpha\in [0,1],\: \lambda,\:\beta,$ and $\gamma$ are constant parameters.

   Then the covariance function of the mother random field is
$$\rho_{\Lambda}(r) = \left({\frac{e^{-2c_z}}{\left({1-{\frac{2}{\lambda}}+{\frac{2}{{\lambda}^2}}(1-\rho_ {Z}(\tau))}\right )^{\beta}}}-1\right) {\left({\frac{e^{-2c_z}}{\left (1-\frac{1}{\lambda}\right)^{\beta}}}-1\right)}^{-1}.$$

   The following result gives the R\'enyi function and the corresponding conditions for Model~2.

\begin{theorem}\textup{\cite{leonenko2013}}\label{theo:6.2}
Suppose that for Model {\rm 2} the parameter $\lambda > 2$ and the correlation function satisfies
\[0 < |\rho_Z(r)| \leq Ce^{-\gamma r},\: r > 0,\]
for some positive constants $C$ and $\gamma.$ Then, for the parameters $(\beta, \lambda)$ from~the~set 
\[L_{\beta, \lambda} = \left\{(\beta, \lambda):b > \left(1+ \frac{\frac{1}{\lambda^2}}{1-\frac{2}{\lambda}}\right)^{\frac{\beta}{2}}, \lambda > 2, \beta >0 \right\},\]
   the measures $\mu_k \xrightarrow{d} \mu, \: k \rightarrow \infty$. The R\'enyi function of $\mu$ is given by
\begin{equation}
    T(q) = q\left(1-\frac{\beta}{2}\log_{b} {\left(1-\frac{1}{\lambda}\right)}\right) + \left(\frac{\beta}{2}\right)\log_{b}\left({1-\frac{q}{\lambda}}\right)-1. \label{eq:11a}
\end{equation}
   where $q \in Q = \{ 0 < q < \lambda\} \cap [1,2] \cap L_{\beta, \lambda}.$ 
\end{theorem}

\noindent    \textbf{Model 3} Let the mother random field be $$\Lambda(x) = exp\left \{ U(x) - c_U \right\}, \: x \in \mathbb{R}^3,$$
where $U(x) = -{Z^{-1}(x)},$ and $Z(x), \: x \in \mathbb{R}^3,$ is a gamma-correlated HIRF with the densities given by \eqref{eq:9a} and~\eqref{eq:10a} and the correlation function $\rho_Z(r).$\\[0mm]

\begin{theorem}\textup{\cite{leonenko2013}}
Suppose that for Model {\rm 3} the correlation function satisfies
$$0 < |\rho_Z(r)| \leq Ce^{-\gamma r},\: r > 0,$$
for some positive constants $C$ and $\gamma.$
Then, for any $(\beta, \lambda) \in L_{\beta, \lambda}$ and $b > \left( \frac{\Gamma(\beta) 2^{\frac{\beta}{2}-1}K_{\beta}(2\sqrt{2\lambda})}{\lambda^{\beta/2}[K_{\beta}(2\sqrt{\lambda})]^2}\right)^{\frac{1}{2}}$ the measures $\mu_k \xrightarrow{d} \mu$ when $\: k \rightarrow \infty.$ The R\'enyi function of measure $\mu$ is
\begin{equation}
   T(q) = q\left(1+\frac{c_U}{2\ln b}\right)- \frac{1}{2}\log_{b}\left(q^{\beta/2}K_{\beta}(2\sqrt{q\lambda})\right)- 
   \left (1 + \frac{1}{2}{\log_{b} \left(\frac{2\lambda^{\beta/2}}{\Gamma(\beta)}\right)}\right),\label{eq:12a} 
\end{equation}
   where $q \in Q = [1,2] \, \cap \, L_{\beta, \lambda},\: K_{\lambda}(x)$ is the modified Bessel function of the third kind and 
$$c_U = \ln {\left(\frac{2\lambda^{\beta/2}K_{\beta}(2\sqrt{\lambda})}{\Gamma(\beta)} \right)}.$$ \label{theo:6.3}
\end{theorem}

   Let $\alpha(q)$ denote the $q^{th}$ order singularity exponent and be defined by
\begin{equation}
    \alpha(q) = \frac{d}{dq}T(q). \label{eq:13a}
\end{equation}
Then the multifractal spectrum defined by \eqref{eq:7a} can be expressed as a function of $\alpha$ by 
\begin{equation}
    f(\alpha(q))= q \cdot \alpha(q) - T(q). \label{eq:14a}
\end{equation}

   For Model 1 it is easy to see from \eqref{eq:8a} that
$$\alpha(q) = 1+ \frac{\sigma_Y^2}{4\ln(b)}-\frac{\sigma_Y^2}{2\ln(b)}q,$$
\begin{equation}
    f(\alpha(q))=1-\frac{\sigma_Y^2}{4\ln(b)}q^2, \quad q \in [1,2]. \label{eq:15a}
\end{equation}

   By \eqref{eq:11a} we obtain for Model~2
$$\alpha(q)= 1-\frac{\beta}{2}\log_b\left(1-\frac{1}{\lambda}\right)+\frac{\beta}{2 \ln(b)(q-\lambda)},$$

\begin{equation}
    f(\alpha(q))= 1+\frac{\beta}{2}\left(\frac{q}{\ln(b)(q-\lambda)}-\log_b \left( 1-\frac{q}{\lambda}\right) \right). \label{eq:16a}
\end{equation}

   For Model~3 it follows from \eqref{eq:12a} and ${K_{\beta}}^{'}(q)=-\frac{1}{2}\left(K_{\beta-1}(q)+K_{\beta+1}(q)\right)$, see 9.6.26 in \cite{abramowitz1948handbook}, that 
\[ \alpha(q) = 1 + \frac{c_U}{2\ln(b)}-\frac{\beta}{4\ln(b)q} + \frac{\sqrt{\lambda}(K_{\beta-1}(2\sqrt{q\lambda})+K_{\beta+1}(2\sqrt{q\lambda}))}{2\ln(b)K_{\beta}(2\sqrt{q\lambda})\sqrt{q}},\]
\[    f(\alpha(q))= 1+{\frac{\beta}{2}\log_b\left( \frac{2{\lambda}^{\beta/2}}{\Gamma({\beta})}\right) {-\frac{\beta}{4\ln(b)}}}+{\frac{1}{2}\log_b(q^{\beta/2}K_{\beta}(2\sqrt{q\lambda}))} \]
\begin{equation}+ {\frac{\sqrt{q\lambda}(K_{\beta-1}(2\sqrt{q\lambda})+K_{\beta+1}(2\sqrt{q\lambda}))}{2\ln(b)K_{\beta}(2\sqrt{q\lambda})}}.\label{eq:17a}
\end{equation}

   Summarising the above results we obtain

\begin{theorem}\label{theo:6.4}
 Let the corresponding conditions of Theorems~{\rm\ref{theo:6.1}, \ref{theo:6.2}} and {\rm \ref{theo:6.3}} are satisfied for Models~{\rm 1, 2} and \rm{3.} Then the multifractal spectra of these models are given by \eqref{eq:15a}, \eqref{eq:16a} and \eqref{eq:17a} respectively.
\end{theorem}
\begin{figure}[htb!]
    \centering
    \subfloat[R\'enyi function of Model 1]{\label{fig11a}
    \includegraphics[width=0.32\textwidth]{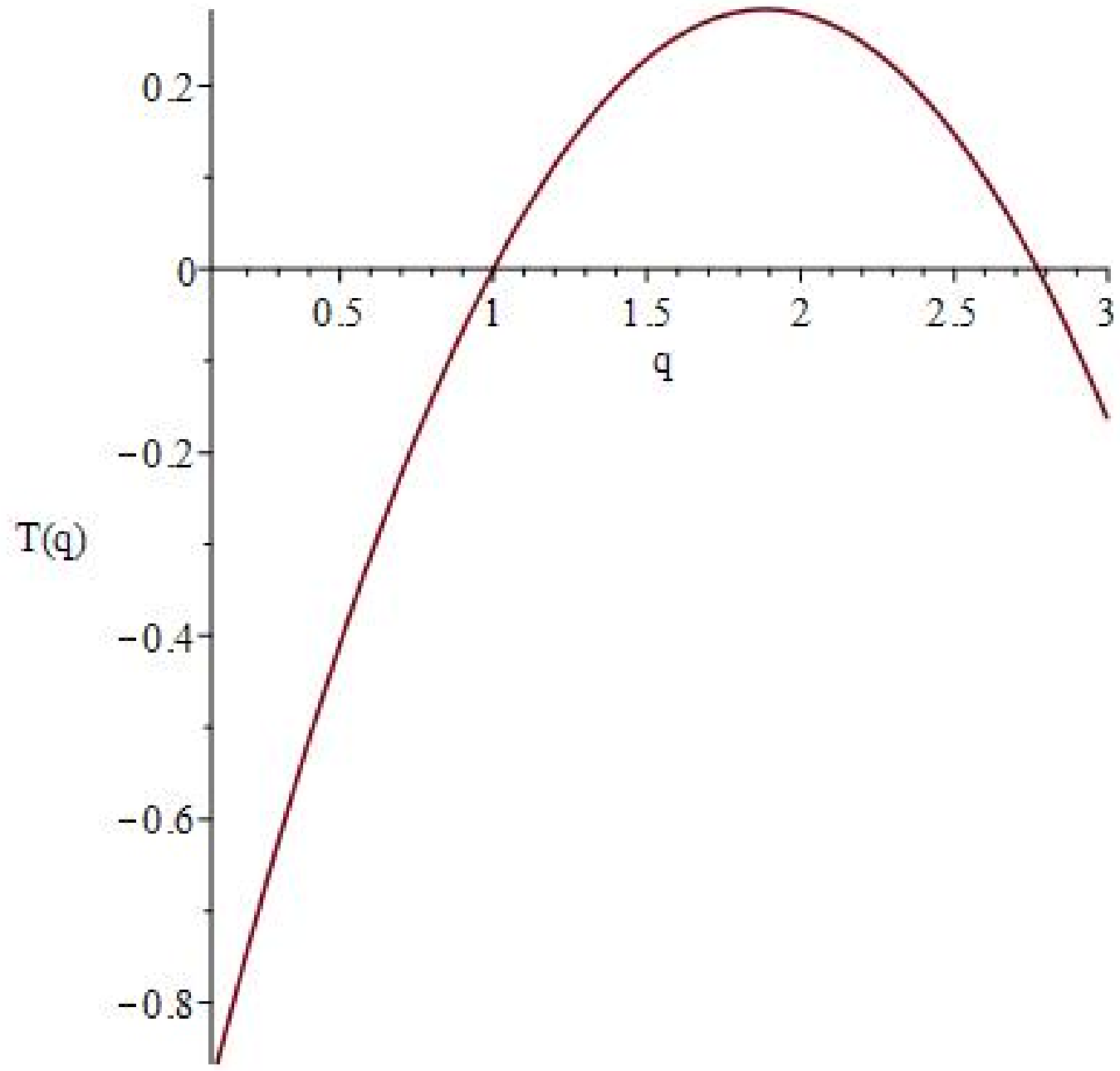}}
    \hfill
    \centering
    \subfloat[R\'enyi function of Model 2]{\label{fig11b}
    \includegraphics[width=0.32\textwidth]{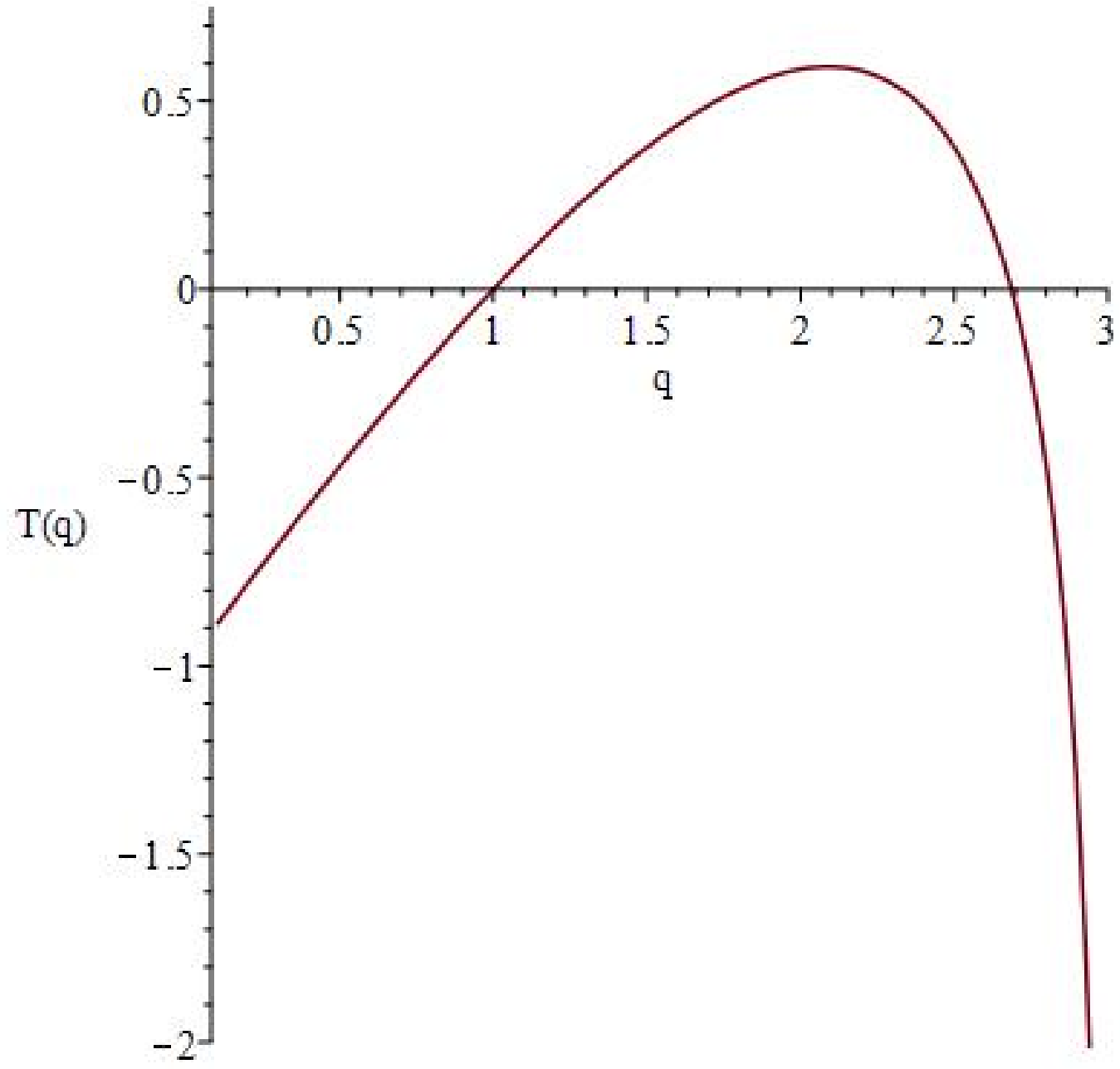}}
    \hfill
    \centering
    \subfloat[R\'enyi function of Model 3]{\label{fig11c}
    \includegraphics[width=0.32\textwidth]{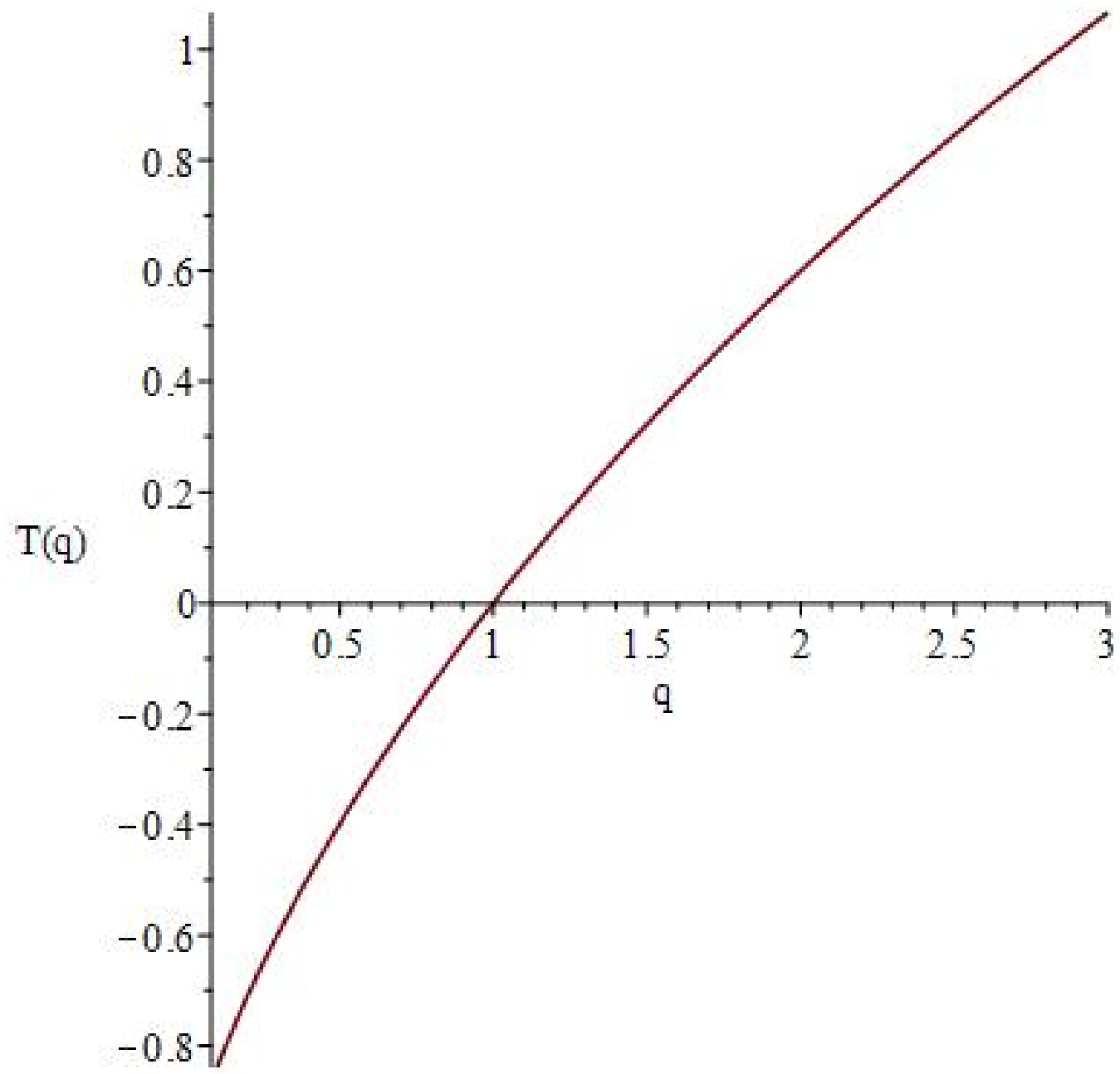}}
    \medskip \\
    \centering
    \subfloat[Spectrum of Model 1]{\label{fig11d}
    \includegraphics[width=0.32\textwidth]{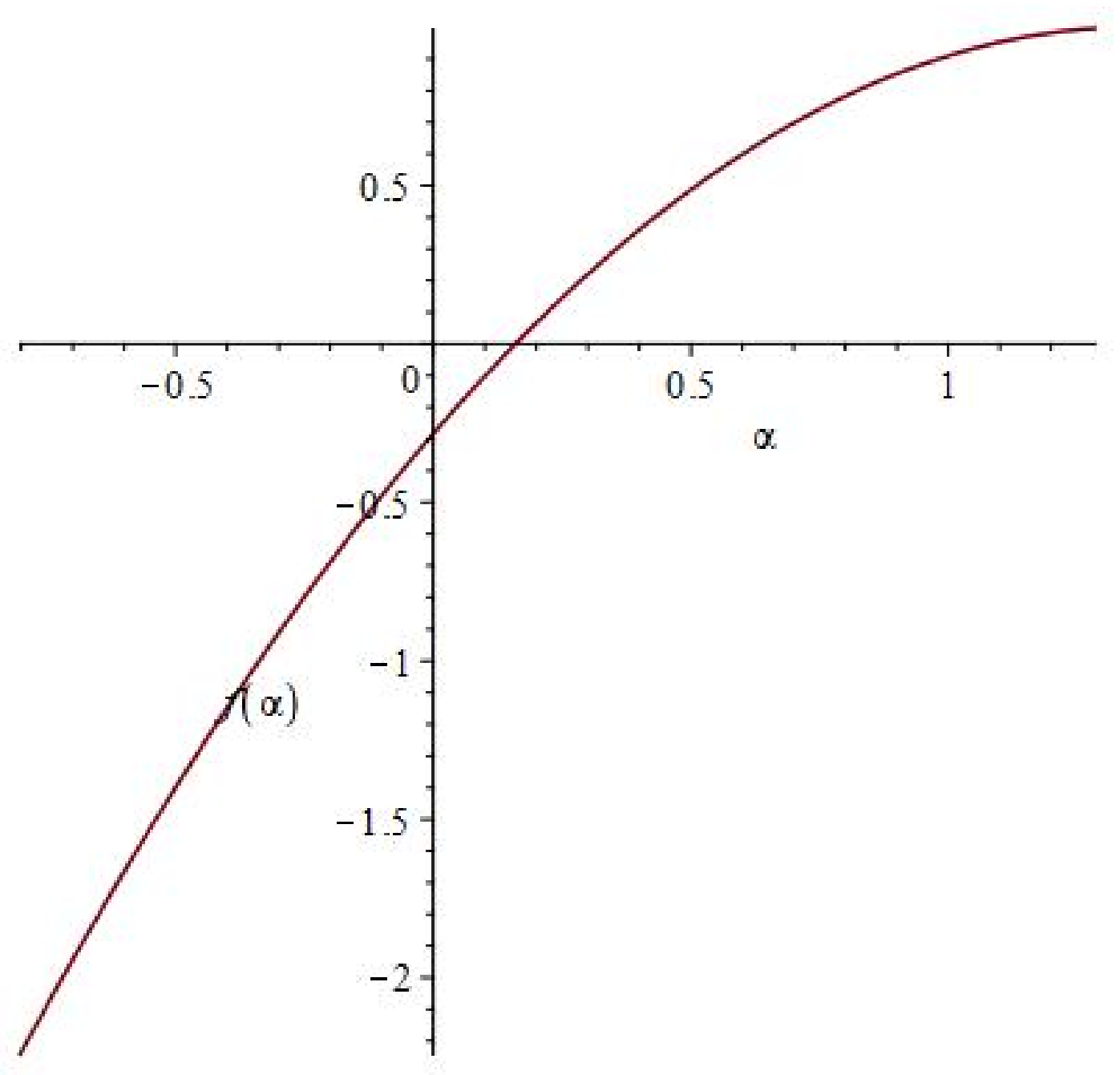}}
    \hfill
    \centering
    \subfloat[Spectrum of Model 2]{\label{fig11e}
    \includegraphics[width=0.32\textwidth]{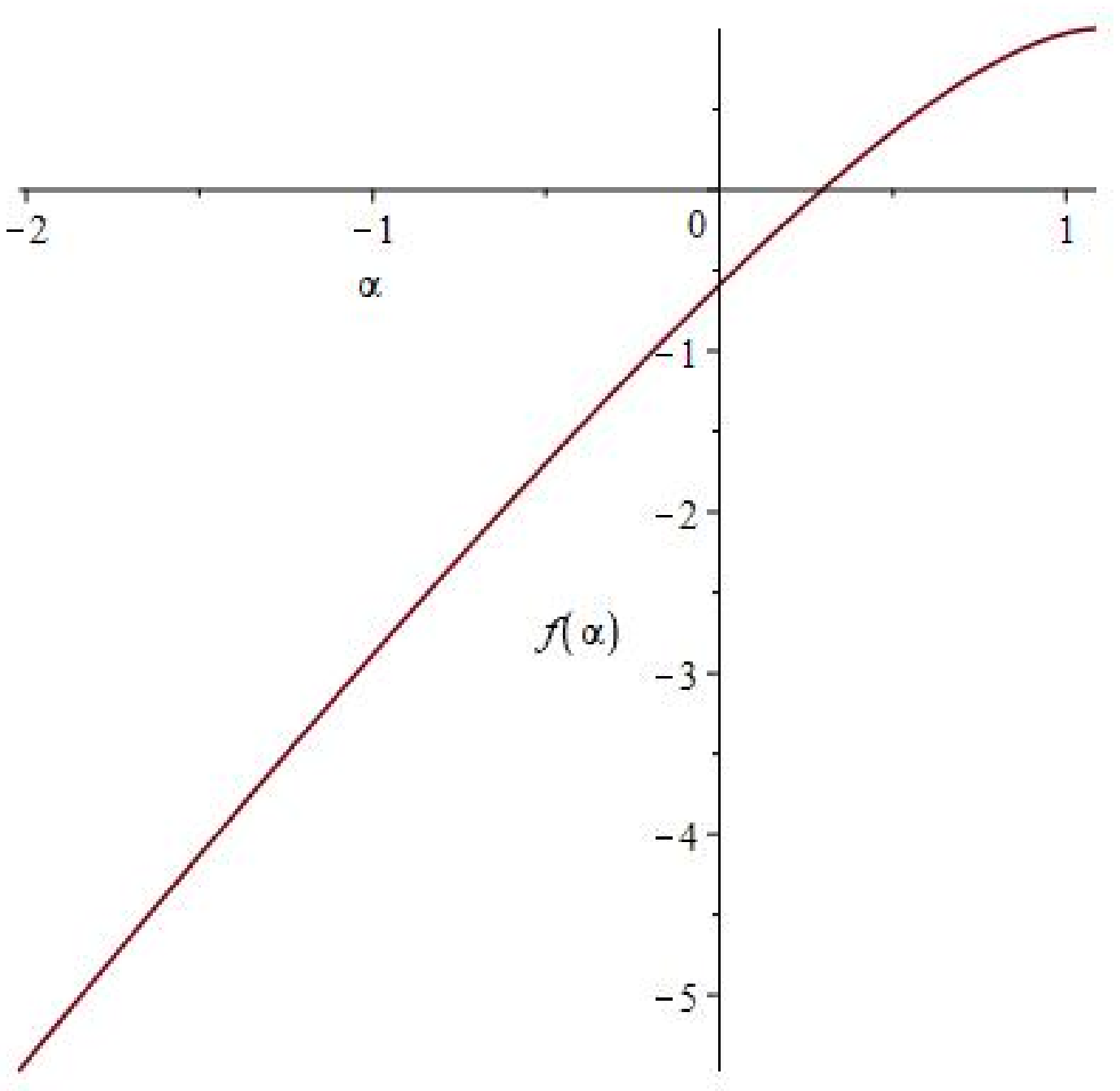}}
    \hfill
    \centering
    \subfloat[Spectrum of Model 3]{\label{fig11f}
    \includegraphics[width=0.32\textwidth]{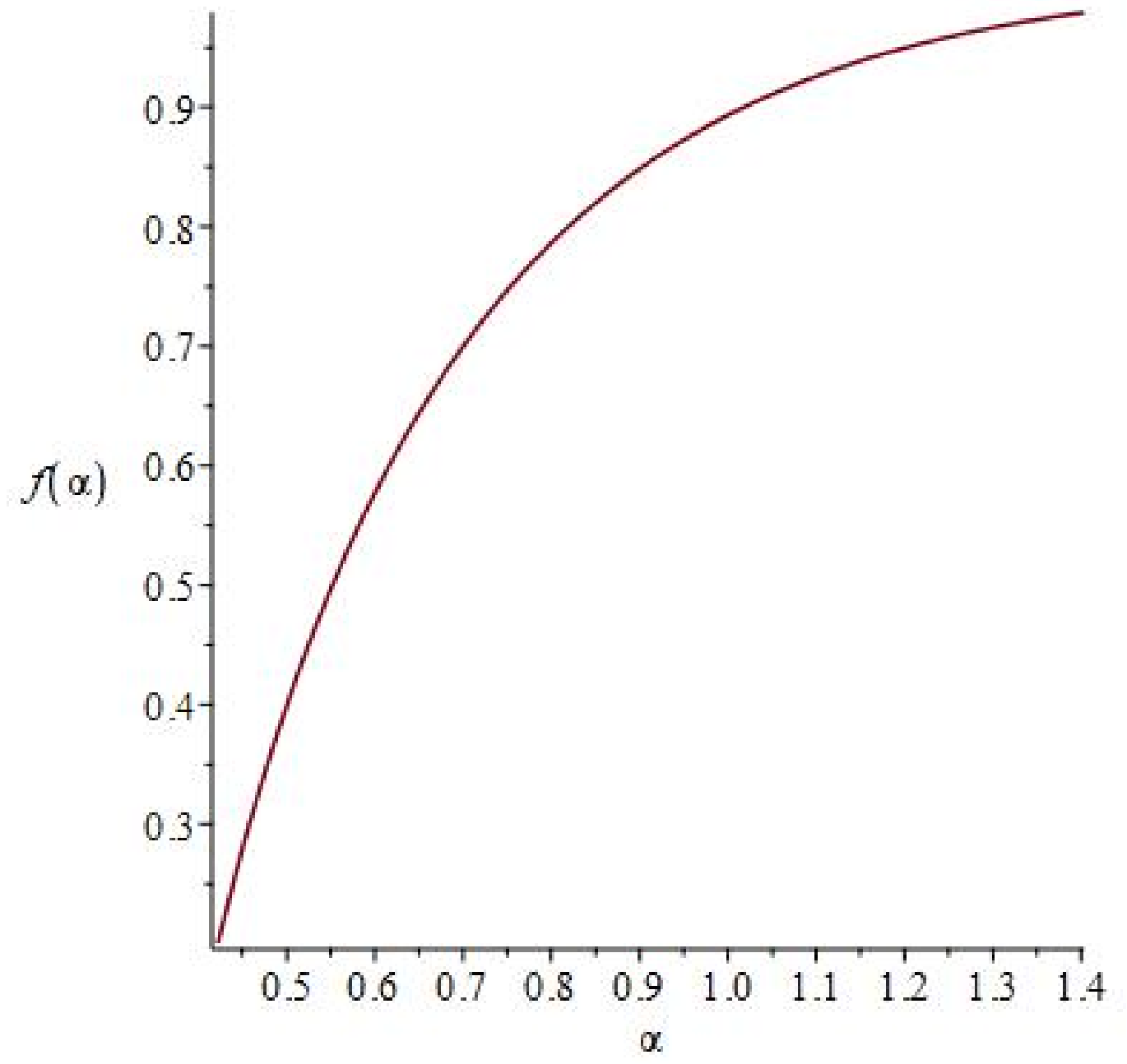}}
    \caption{Examples of R\'enyi functions and multifractal spectra for Models 1, 2 and 3}
    \label{fig11}
\end{figure}

The plots shown in Fig.~\ref{fig11} illustrate behaviours of the R\'enyi functions and multifractal spectra for Models 1, 2 and 3. For these numerical examples, we used the following values of the parameters: $b=2, \: \sigma_Y=1, \: \lambda=3, \; \text{and} \; \beta=2.$ Notice that these values of $b, \: \lambda \: \text{and} \: \beta$ satisfy the conditions in $L_{\beta, \lambda}$. We also selected $(0, 3)$ as the range of $q$ values. It is slightly wider than the range $[1, 2]$ in the theorems and allows better visualisation of $T(q)$ and $f(\alpha)$, see Section~\ref{S1:81} on the way to check its validity.

Fig.~\ref{fig11} shows that the R\'enyi functions of the Models~1 and~2 have parabolic shapes while the  R\'enyi function of the Model~3 is closer to a linear shape on the interval $(0, 3)$. Also, comparing the plots for Models~1, 2 and 3, we can see that the R\'enyi functions of Model~1 and 2 exhibit a concave down increasing and decreasing behaviour within $q \in (0, 3),$ whereas for Model~3 it increases. The multifractal spectra of Models~1, 2 and 3 show a concave down increasing behaviour within $q \in (0, 3).$

\section{Models based on power transformations of Gaussian fields}
\label{S1:7}

   In the previous section, we considered three models based on an exponential transformation of Gaussian or gamma-correlated HIRF. This section proposes few much simpler scenarios where conditions of the theorems from Section~\ref{S1:4} are satisfied.

   First, note that the condition $\Lambda(x)>0$ in \cite{leonenko2013} can be relaxed to $\Lambda(x)>0$ almost~sure.\\[0mm]

  \noindent  \textbf{Model~4} Let $\Lambda(x) = Y^2(x)$, where $Y(x), \: x \in \mathbb{R}^3,$ is a zero-mean unit variance Gaussian HIRF with a covariance function $\rho_Y(\tau), \: \tau \geq 0.$\\[0mm]

For this model we obtain the following result, see the proof in Appendix~\ref{appen}.

\begin{theorem}
Suppose that for Model {\rm 4,} the correlation function of $Y(x)$ satisfies $|\rho_{Y}(r)| \leq Ce^{-\gamma r}, r>0, \gamma>0, \: \text{and} \; b> \max(\sqrt[3]{1+\sigma^2_{\lambda}}, e^{{\sigma^2_{\Lambda}C/3}}).$

Then the measures $\mu_k \xrightarrow{d} \mu, \: k \rightarrow \infty,$ and the corresponding R\'enyi function is equal to

\begin{equation}
    T(q) = q-1-\frac{1}{2}\log_b\left( \frac{2^q \Gamma(q+\frac{1}{2})}{\sqrt{\pi}}\right), \: q \in [1,2].\label{eq:21a}
\end{equation}\label{theo:7.1}
\end{theorem}
\begin{examp}\label{examp1}
The approach developed in Section~{\rm \ref{S1:5}} can be used to obtain the moment conditions for $q \in [1, 4]$ in the case of Model~{\rm 4}. As it is shown in Appendix~\ref{appen} it is enough to require that $b > e^{{\frac{\sigma{(\max(\sigma_{\Lambda}C, 1))}^4}{3}}}.$
\end{examp}
Now we show that the assumption $\Lambda(x)>0$ almost surely is indeed not restrictive and it is easy to construct a modification of Model~4 with $\Lambda(x) > 0.$\\[0mm]

   \noindent \textbf{Model 4'} Let $$\Lambda(x) = Y^{2}(x)\cdot (1-\varepsilon) + \varepsilon, \: \varepsilon \in (0,1),$$ where $Y(x), \: x \in \mathbb{R}^n,$ is a zero-mean unit variance Gaussian HIRF with a covariance function $\rho_Y(\tau), \: \tau\geq 0.$\\[0mm]

   It is easy to see that
$$E \Lambda(x) = (1-\varepsilon) EY^2(x)+\varepsilon =1, \quad \sigma_{\Lambda}^2 = Var \Lambda(x) = 2(1-\varepsilon)^2 < +\infty,$$
$$Cov(\Lambda(x),\Lambda(y)) = 2(1-\varepsilon)^2 \rho^2_Y(\Vert x-y \Vert).$$

   Hence, Model 4' satisfies Conditions~{\rm\ref{cond1}} and~{\rm\ref{cond2}} and $|\rho(r)| \leq Ce^{-\gamma r}, \: r>0,$ $\gamma>0,$ if $|\rho_Y(r)| \leq C'e^{-{\gamma}^{'}r}, \: r>0, \gamma^{'}>0.$

   Therefore, we obtain
\begin{align*}
    T_{\varepsilon}(q) &= q-1-\frac{1}{2}\log_b E((1-\varepsilon) Y^2(x) + \varepsilon)^q \\
    &= q-1-\frac{q}{2}\log_b(1-\varepsilon)-\frac{1}{2}\log_b E\left(Y^2(x)+ \frac{\varepsilon}{1-\varepsilon}\right)^q
\end{align*}
   and 
\[ \lim_{\varepsilon \to 0} T_{\varepsilon}(q) = q-1-\frac{1}{2}\log_b E(Y^{2q}(x)),\]
  which coincides with \eqref{eq:21a}.

   The next model generalizes Model~4 to an arbitrary even power of a Gaussian random field.\\[0mm]

 \noindent   \textbf{Model~5} Let $\Lambda(x)=Y^{2k}(x)$, $k\in{\mathbb N,}$ where $Y(x), \: x \in \mathbb{R}^3,$ is a zero-mean Gaussian HIRF with the variance $\sigma^{2}=\left( \frac{\sqrt{\pi}}{2^k \Gamma(k+ \frac{1}{2})}\right)^{-\frac{1}{k}}$ and a covariance function $\rho_Y(r), r \geq 0.$

\begin{theorem}
Suppose that for Model~{\rm 5} the correlation function of $Y(x)$ satisfies 
$|\rho_Y(r)| \leq Ce^{-\gamma r},$ $r>0,$ $\gamma>0, \: \text{and} \;  b>\max\left(\sqrt[3]{1+\sigma_{\Lambda}^2},\: e^{\frac{\sigma^2_{\Lambda}C}{3}}\right).$

   Then the measures $\mu_k \xrightarrow{d} \mu, \: k \rightarrow \infty,$ and the R\'enyi function is given by 
   \begin{equation}
      T(q) = q-1-\frac{1}{2}\log_{b} E Y^{2kq}(x)
       = q-1-\frac{1}{2}\log_{b}\left(\frac{2^{kq}\Gamma(kq+\frac{1}{2})}{\sqrt{\pi}}\right). \label{eq:24a}
   \end{equation} for $q \in [1,2]$. \label{theo:7.2}
\end{theorem}

   The following model shows how vector-valued random fields can be used to construct mother fields.\\[0mm]
   
 \noindent   \textbf{Model~6} Let $\Lambda(x) = \frac{2}{k}Y(x)$, $k\in{\mathbb N,}$ where $Y(x) \sim {\chi}^{2}(k),$ and the HIRF field $Y(x), \: x \in \mathbb{R}^3,$ has a covariance function $\rho_Y(r), \: r\geq 0.$

\begin{theorem}
Suppose that the correlation function in Model {\rm 6} satisfies the inequality
$|\rho_Y(r)| \leq Ce^{-\gamma r},$ $r>0,$ $\gamma>0,$ and $b>\max\left(\sqrt[3]{1+\sigma_{\Lambda}^2},\: e^{\frac{\sigma^2_{\Lambda}C}{3}}\right).$

   Then the measures $\mu_k \xrightarrow{d} \mu, \: k \rightarrow \infty,$ and for $q \in [1,2]$ the R\'enyi function is  equal to
    \begin{equation} T(q) = q\left(1-\frac{1}{2}\log_b\left(\frac{2}{k}\right)\right)-1-\frac{1}{2}\log_b\left(2^q \frac{\Gamma(q+\frac{k}{2})}{\Gamma(\frac{k}{2})}\right).\label{eq:25a}
   \end{equation} 
\label{theo:7.3}
\end{theorem}

   Note that $\Gamma^{'}(x) = \psi(x)\Gamma(x),$ where $\psi(x)$ is the digamma function defined by 
\[ \psi(x) = \int_{0}^{\infty}\left(\frac{e^{-t}}{t}-\frac{e^{-xt}}{1-{e}^{-t}}\right)dt.\]

   Then, it follows from \eqref{eq:13a}, \eqref{eq:14a} and \eqref{eq:21a} that for Model~4
\[ \alpha(q) = 1 - \frac{1}{2}\log_{b}2-\frac{\psi(q+\frac{1}{2})}{2\ln 2},\]
\begin{equation}
    f(\alpha(q)) = 1+\frac{1}{2}\log_b\left(\frac{\Gamma(q+\frac{1}{2})}{\sqrt{\pi}}\right)-\frac{q\psi(q+\frac{1}{2})}{2\ln 2}.\label{eq:26a}
\end{equation}

   Analogously, for Model~5 one gets from \eqref{eq:24a}
\[ \alpha(q) = 1 - \frac{k}{2}\log_{b}2-\frac{k\psi(kq+\frac{1}{2})}{2\ln 2},\]
\begin{equation}
    f(\alpha(q)) = 1+\frac{1}{2}\log_b\left(\frac{\Gamma(kq+\frac{1}{2})}{\sqrt{\pi}}\right)-\frac{kq\psi(kq+\frac{1}{2})}{2\ln 2}.\label{eq:27a}
\end{equation}

   Finally, it follows from \eqref{eq:25a} that for Model~6
\[ \alpha(q) = 1 - \frac{1}{2}\log_{b}\left(\frac{2}{k}\right)-\frac{1}{2}\log_{b}2-\frac{\psi(q+\frac{k}{2})}{2\ln 2},\]
\begin{equation}
    f(\alpha(q)) = 1+\frac{1}{2}\log_b\left(\frac{\Gamma(q+\frac{k}{2})}{\Gamma(\frac{k}{2})}\right)-\frac{q\psi(q+\frac{k}{2})}{2\ln 2}.\label{eq:28a}
\end{equation}

   Summarising, we obtain

\begin{theorem}\label{theo:7.4}
Let the corresponding conditions of Theorems~{\rm\ref{theo:7.1}}, {\rm\ref{theo:7.2}} and {\rm\ref{theo:7.3}} are satisfied for Models~{\rm 4, 5} and~{\rm 6.} Then the multifractal spectra of these models are given by \eqref{eq:26a}, \eqref{eq:27a} and \eqref{eq:28a} respectively.
\end{theorem}

\begin{figure}[!htb]
    \centering \vspace{-0.3cm}
    \subfloat[R\'enyi function of Model 4]{\label{fig12a}
    \includegraphics[width=0.32\textwidth]{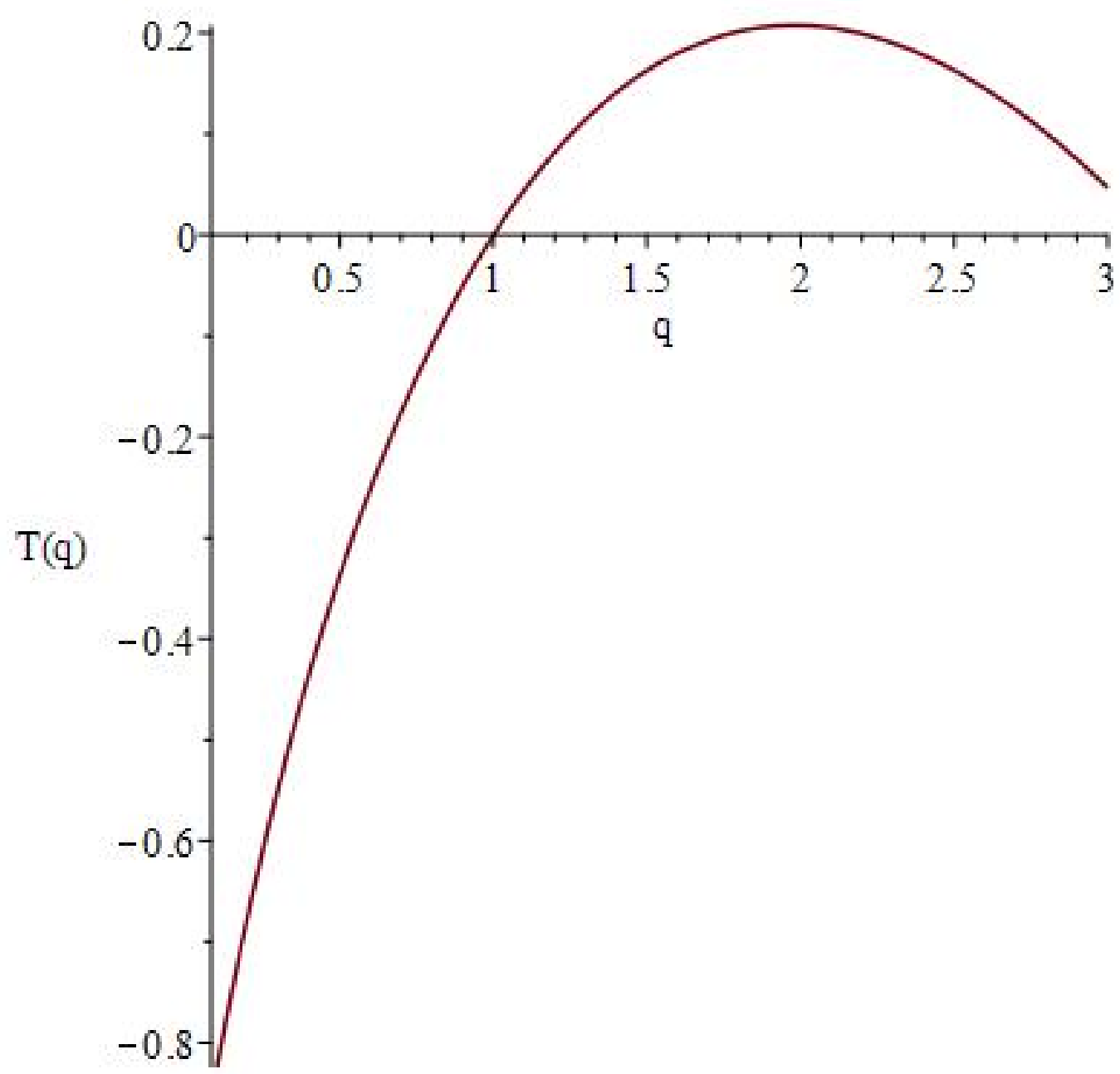}}
    \hfill
    \centering
    \subfloat[R\'enyi function of Model 5]{\label{fig12b}
    \includegraphics[width=0.32\textwidth]{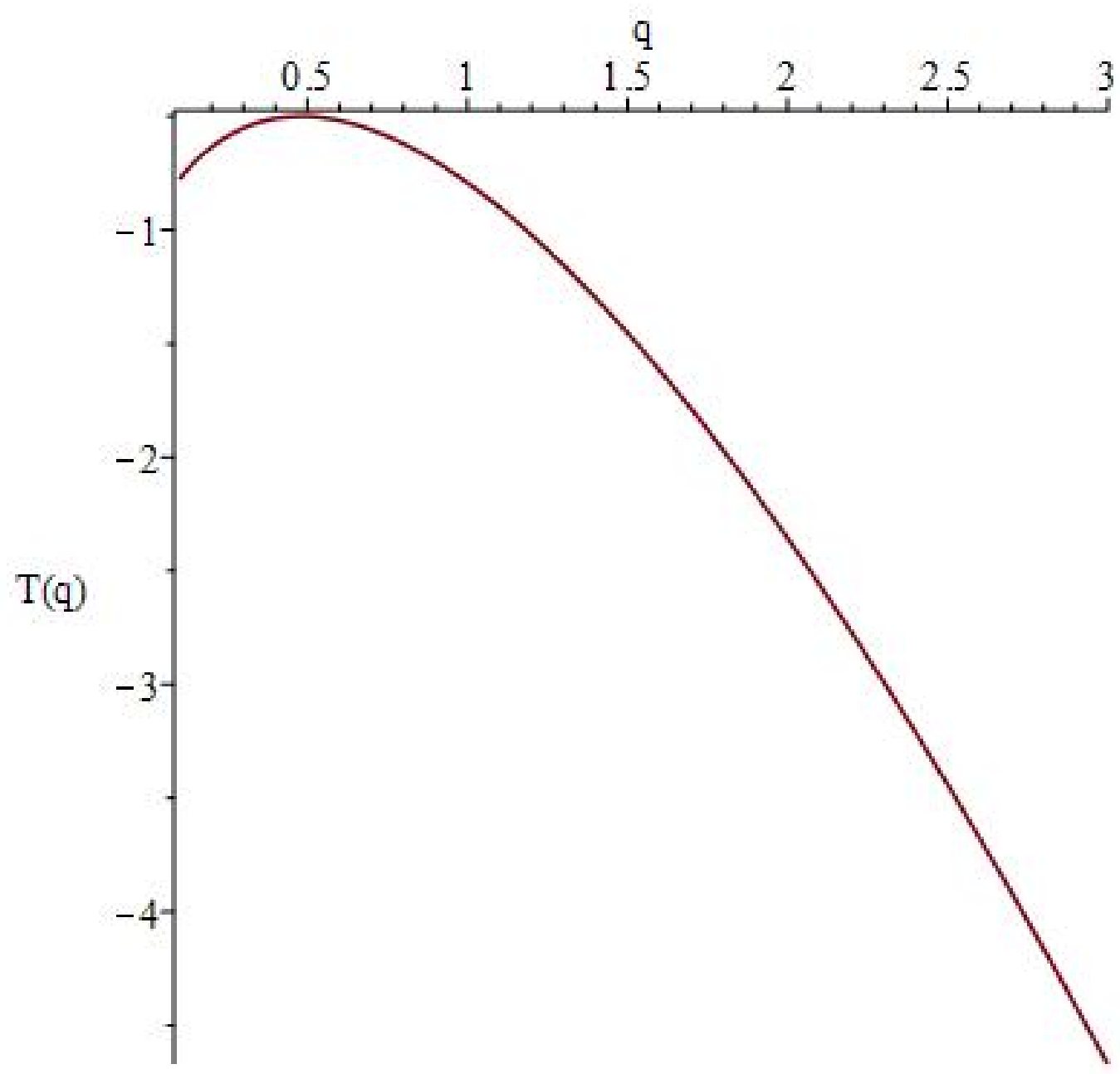}}
    \hfill
    \centering
    \subfloat[R\'enyi function of Model 6]{\label{fig12c}
    \includegraphics[width=0.32\textwidth]{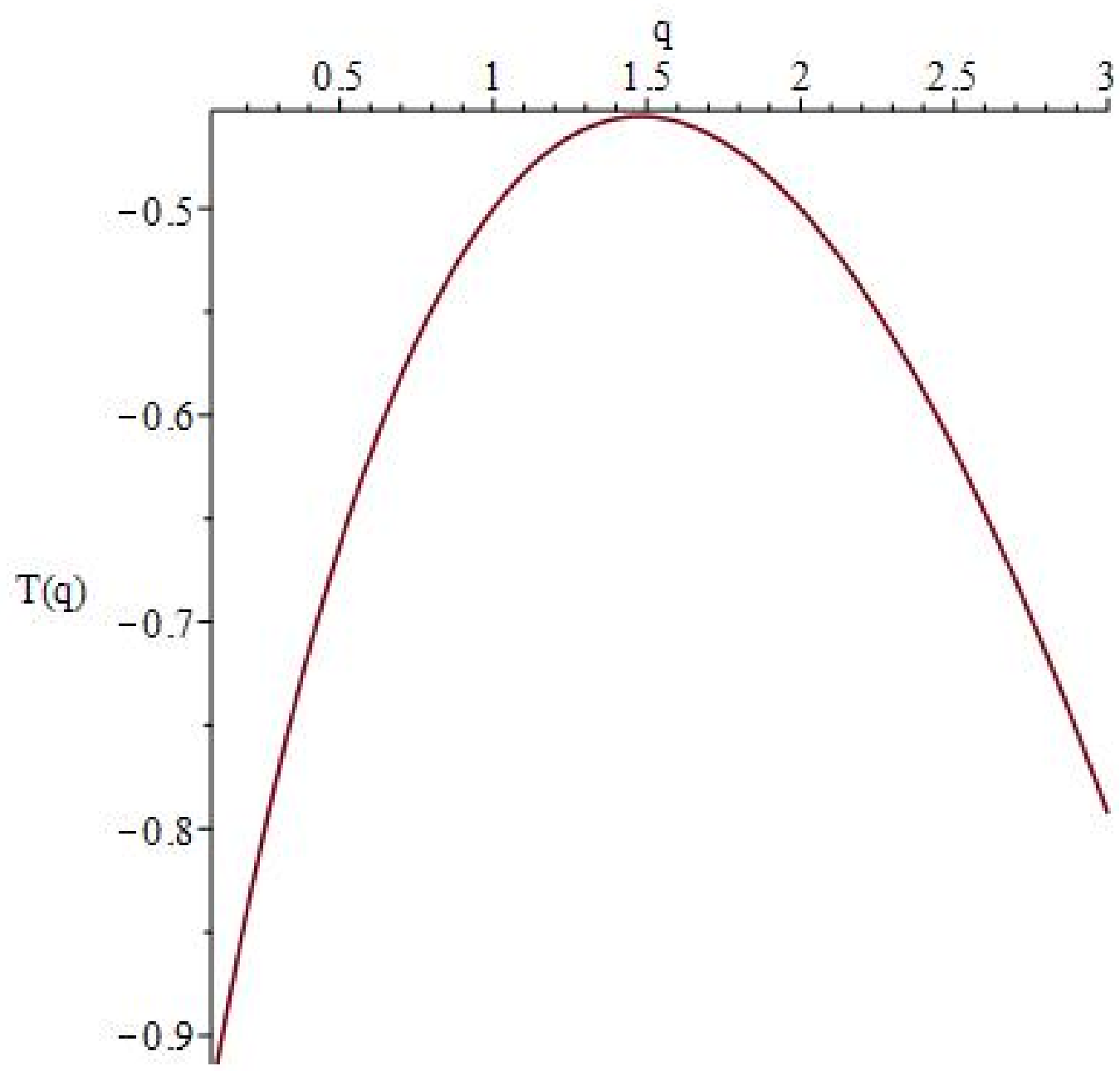}}
    \medskip \\
    \centering \vspace{-0.4cm}
    \subfloat[Spectrum of Model 4]{\label{fig12d}
    \includegraphics[width=0.32\textwidth]{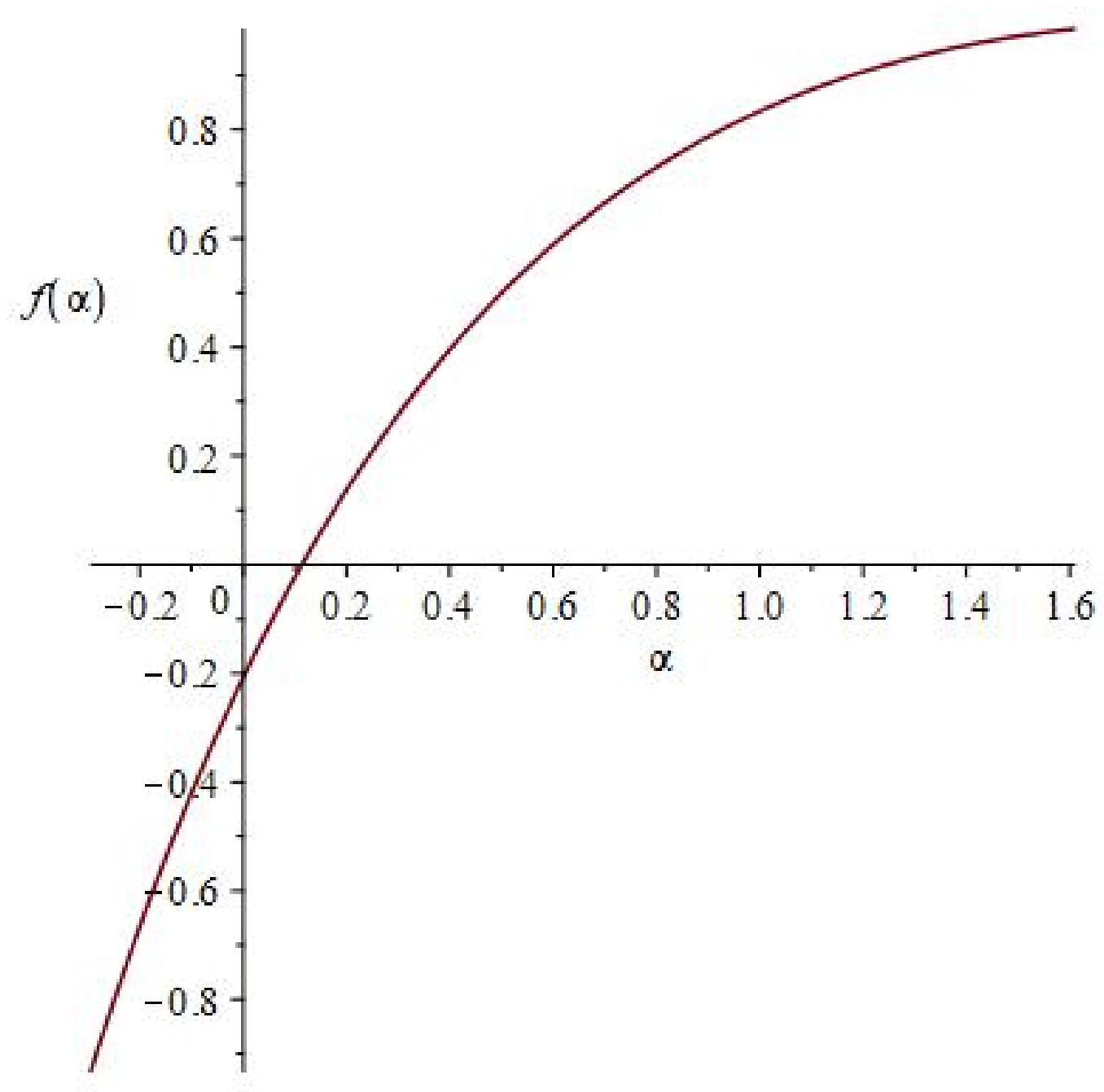}}
    \hfill
    \centering
    \subfloat[Spectrum of Model 5]{\label{fig12e}
    \includegraphics[width=0.32\textwidth]{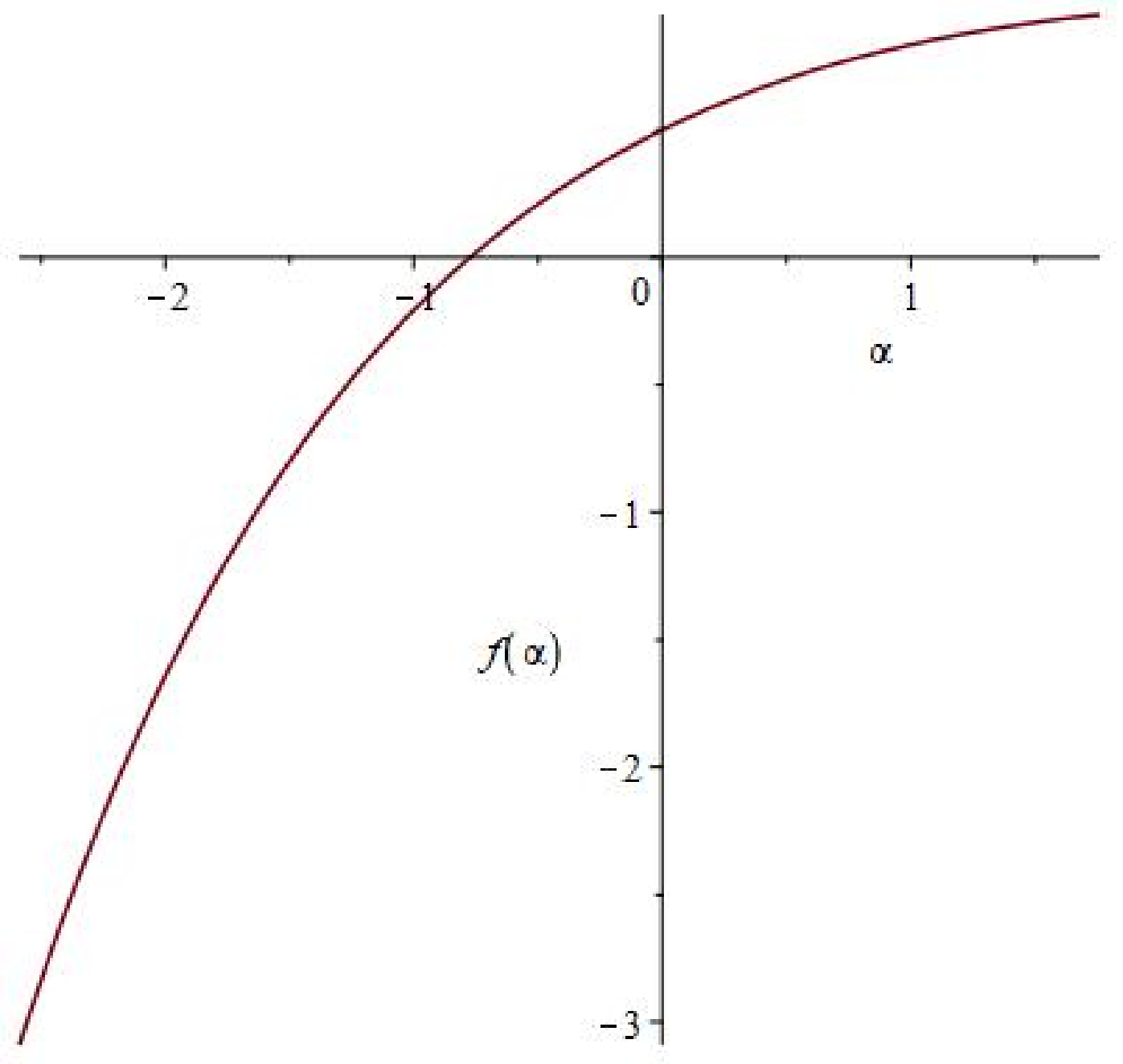}}
    \hfill
    \centering
    \subfloat[Spectrum of Model 6]{\label{fig12f}
    \includegraphics[width=0.32\textwidth]{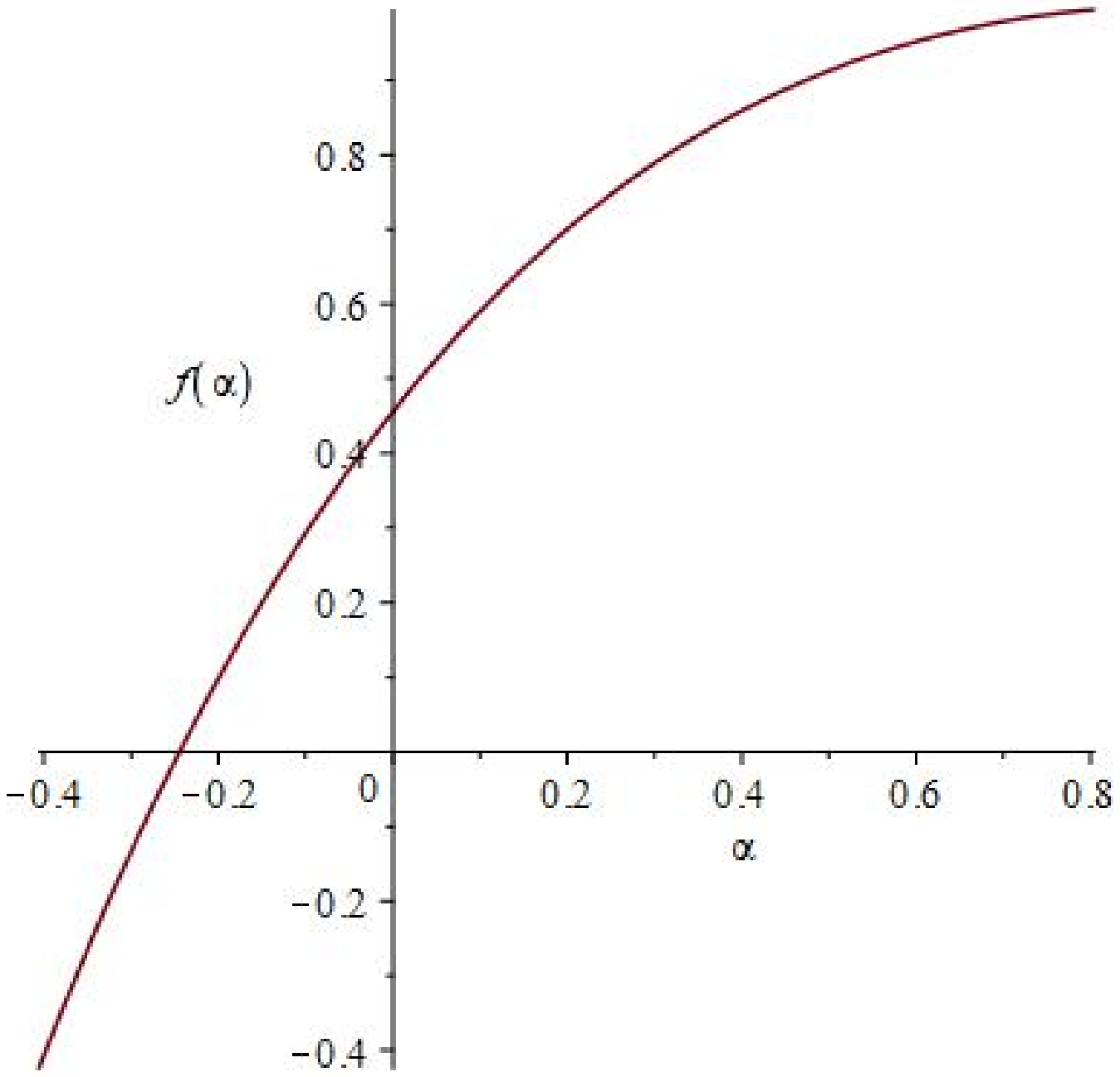}}
    \caption{Examples of R\'enyi functions and multifractal spectra for Models 4, 5 and 6}
    \label{fig12}
\end{figure}

The plots in Fig.~\ref{fig12} demonstrate the behaviours of the R\'enyi functions and multifractal spectra of Models 4, 5 and 6. To plot Fig.~\ref{fig12} we used similar settings and coordinate ranges as in Fig.~\ref{fig11}. The following values of parameters were chosen to produce the plots:
$b=2, \: k=2$ and $\sigma_Y=1$. 

Fig.~\ref{fig12} shows that similar to Models~1 and~2, Models~4, 5 and 6 exhibit a parabolic-type behaviour. The spread of $T(q)$ values is wider for Model~5 when compared to Models~4 and 6 within $q \in (0, 3)$. The R\'enyi functions of Models~4, 5 and 6 manifest a concave down increasing and decreasing behaviour within $q \in (0, 3)$. Similar to Models~1, 2 and 3, the multifractal spectra of Models~4, 5 and 6 show a concave down increasing behaviour within $q \in (0, 3)$. 

\begin{figure}[!htb]
    \centering
    \vspace{-0.2cm}
    \subfloat[R\'enyi functions of Model~1]{\label{fig13a}
    \includegraphics[trim={0cm 0cm 0cm 0cm},clip,width=0.32\textwidth]{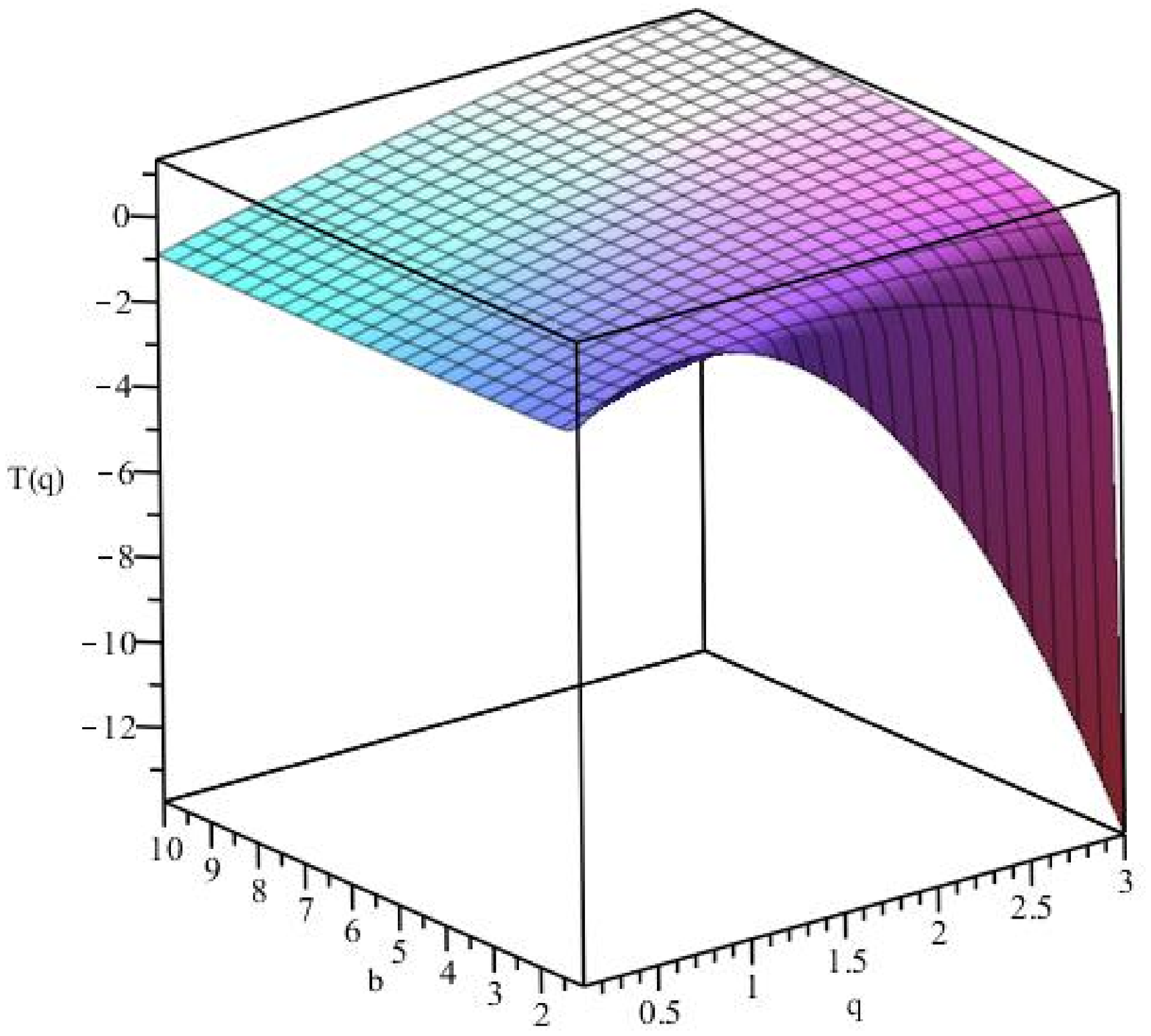}}
    \hfill
    \centering
    \subfloat[R\'enyi functions of Model~2]{\label{fig13b}
    \includegraphics[trim={1.3cm 0cm 1.3cm 0cm},clip,width=0.32\textwidth]{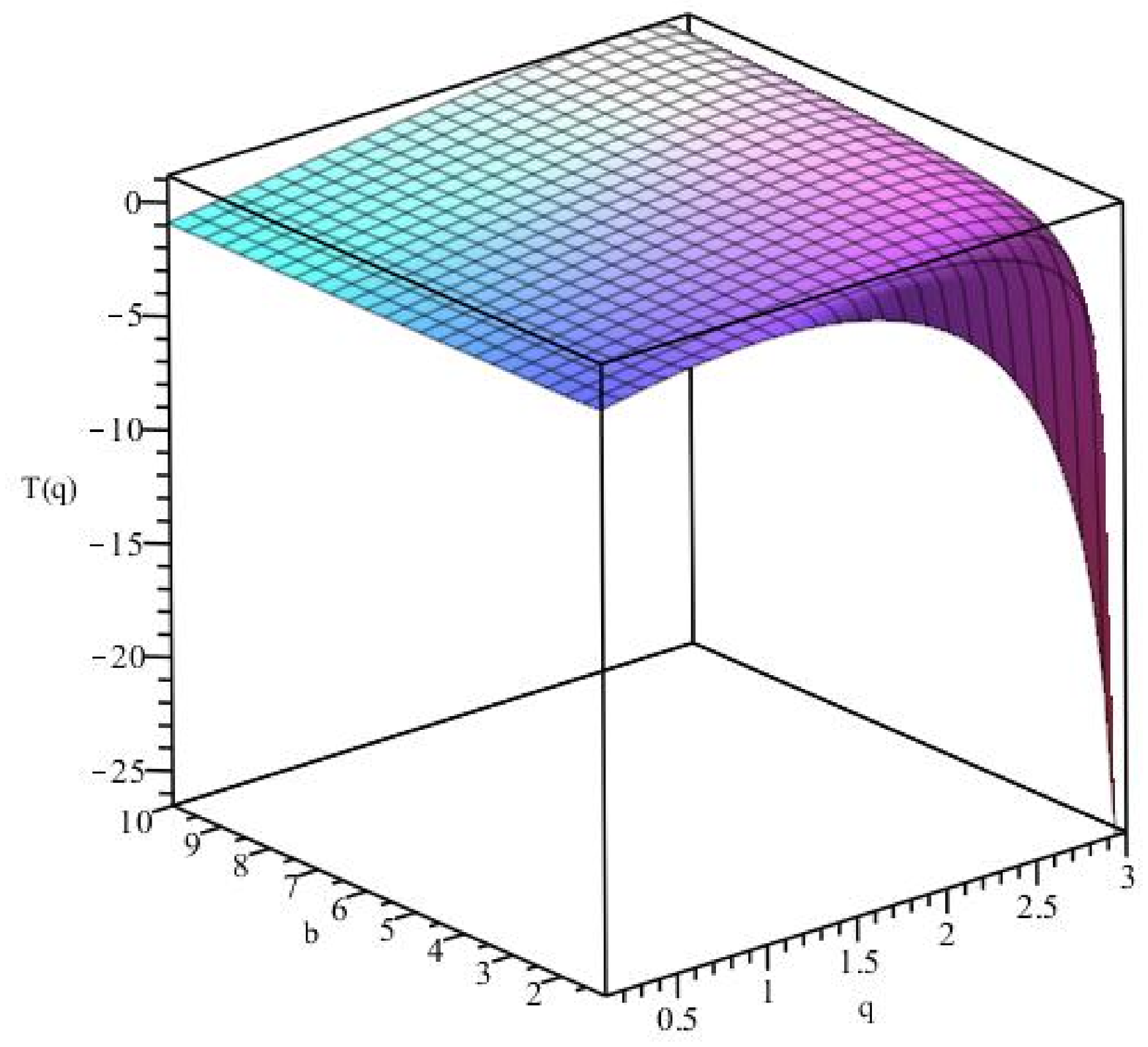}}
    \hfill
    \centering
    \subfloat[R\'enyi functions of Model~3]{\label{fig13c}
    \includegraphics[trim={0cm 0cm 0cm 0cm},clip,width=0.32\textwidth]{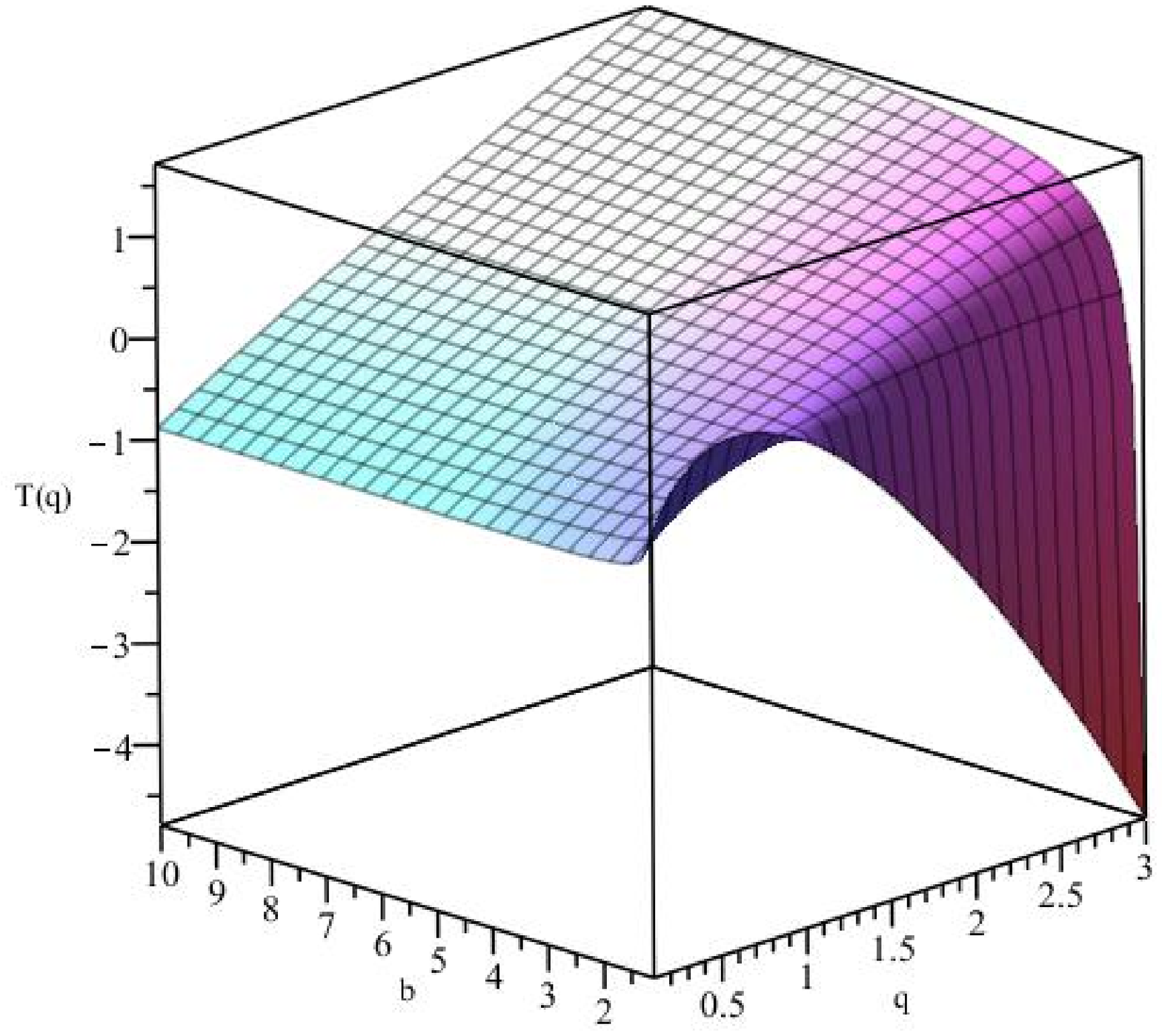}}
    \medskip \\
    \centering\vspace{-0.4cm}
    \subfloat[R\'enyi functions of Model~4]{\label{fig13d}
    \includegraphics[trim={0.5cm 0cm 0.5cm 0cm},clip,width=0.32\textwidth]{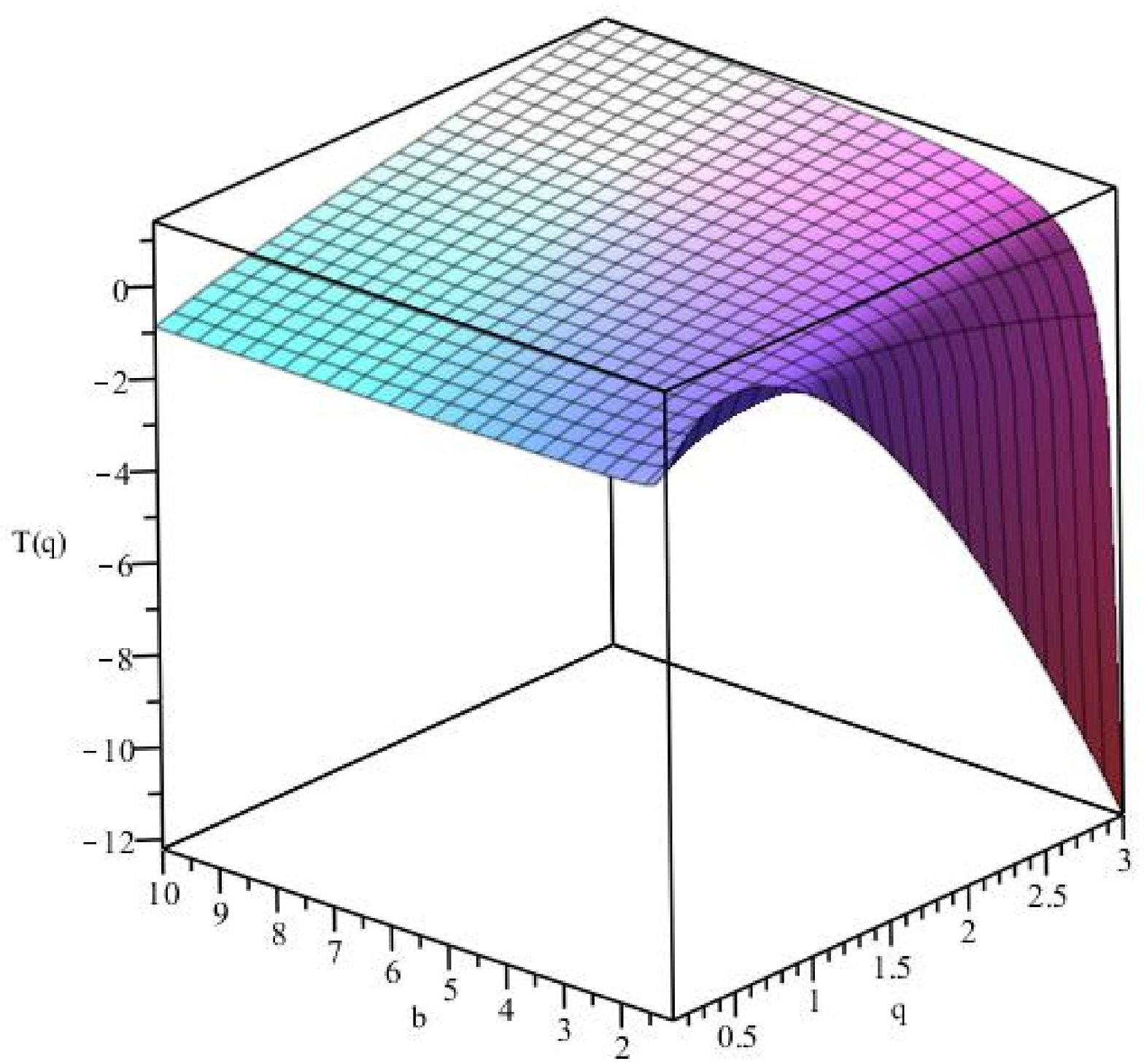}}
    \hfill
    \centering
    \subfloat[R\'enyi functions of Model~5]{\label{fig13e}
    \includegraphics[trim={0cm 0cm 0cm 0cm},clip,width=0.32\textwidth]{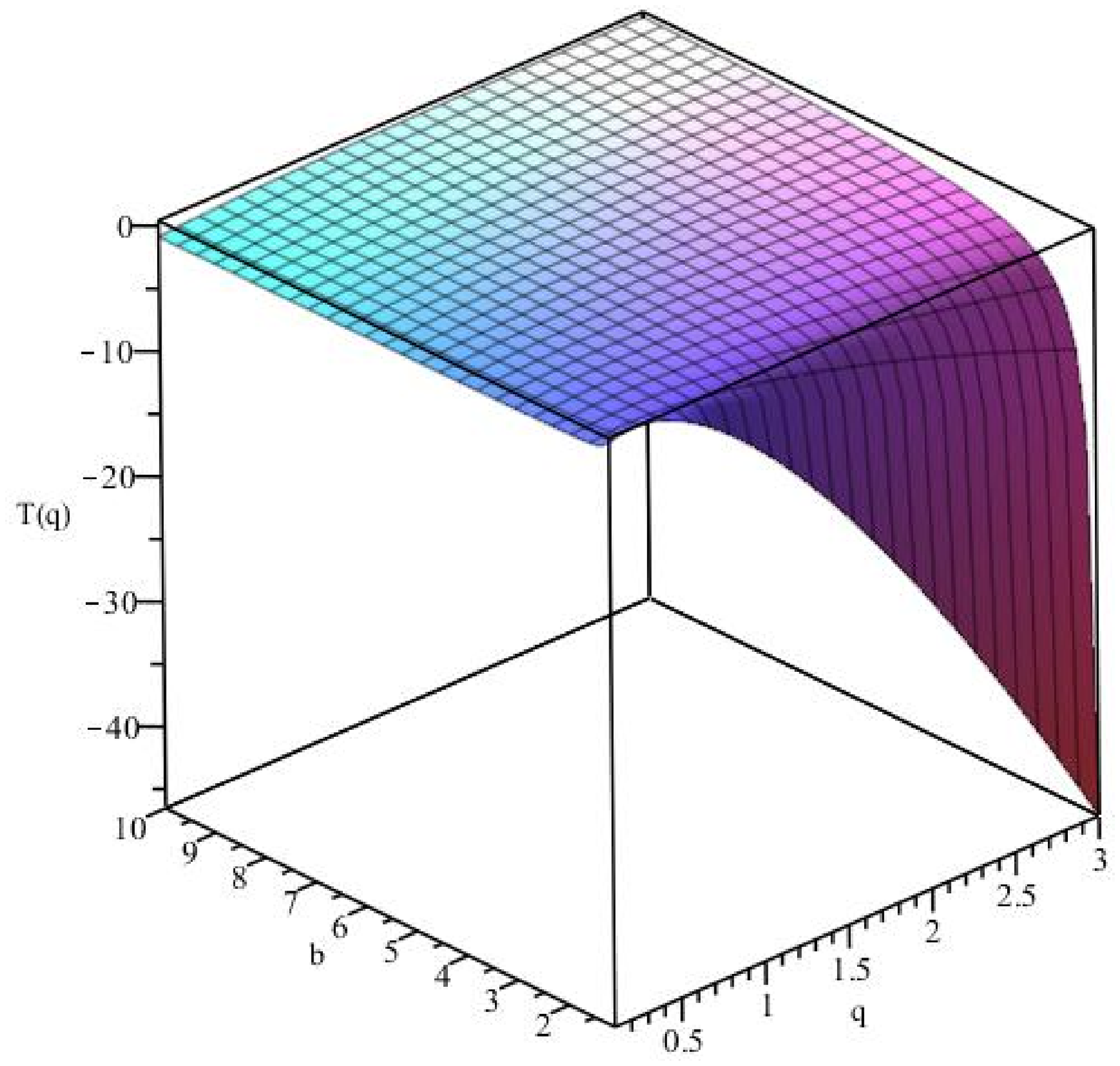}}
    \hfill
    \centering
    \subfloat[R\'enyi functions of Model~6]{\label{fig13f}
    \includegraphics[trim={0cm 0cm 0cm 0cm},clip,width=0.32\textwidth]{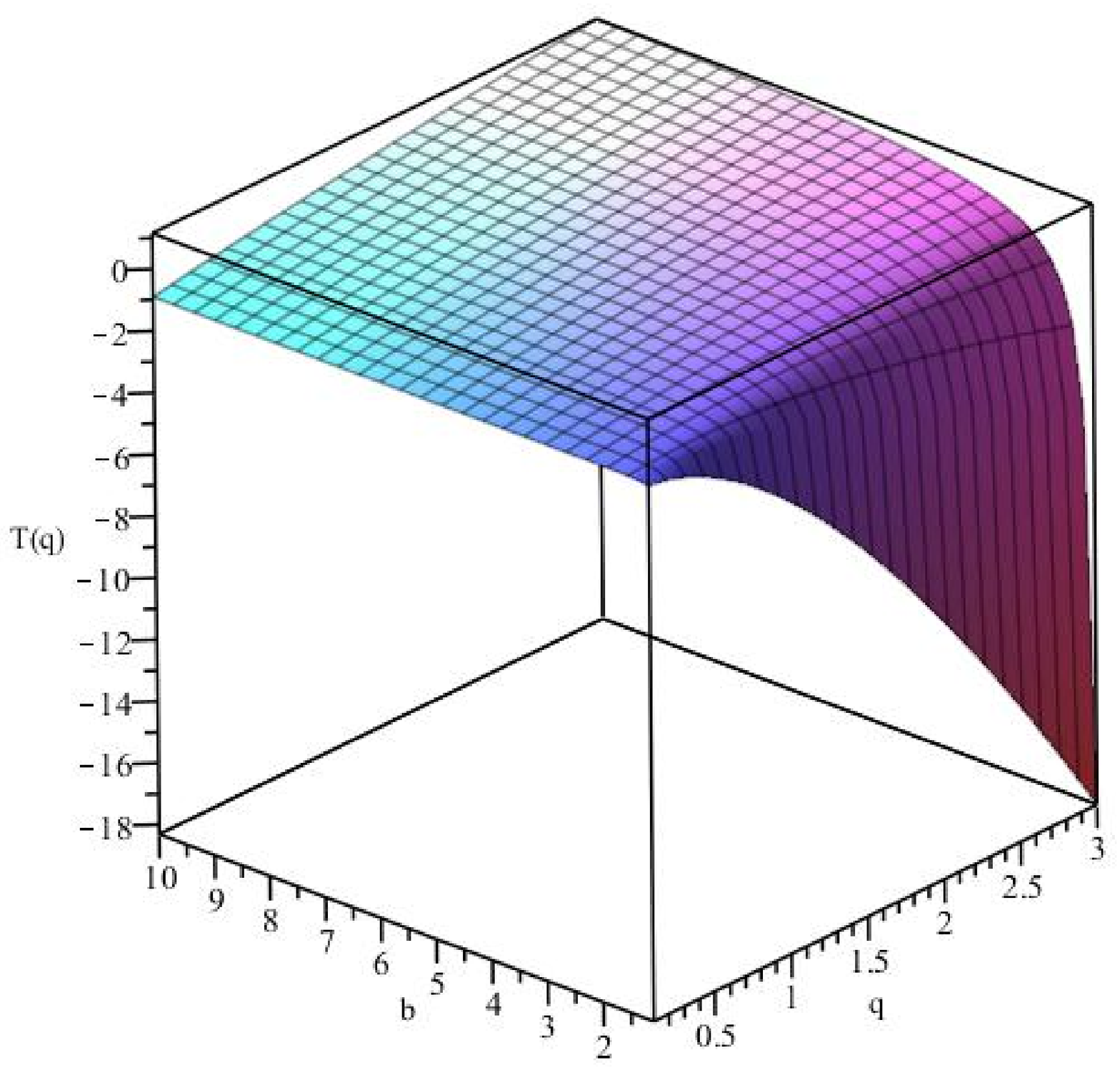}}
    \caption{Dependence of the R\'enyi function on the parameter b}
    \label{fig13}\vspace{-0.3cm}
\end{figure}

Finally, we illustrate the impact of parameter~$b$ on the R\'enyi function using Models 1, 2, 3, 4, 5 and 6. For Model 1, $\sigma_Y=1$ was selected. Then, for Models~2 and 3, the parameters $\lambda=3$ and $\beta=2$ were used. The parameter $k=2$ was chosen for Models~5 and 6. Fig.~\ref{fig13} suggests that the R\'enyi functions for all the models exhibit a similar pattern. It suggests that the R\'enyi functions for Models~1, 3 and 4 are more concave than for other models for small values of~$b$. The dispersion of $T(q)$ values on the interval $(0,3)$ is smallest for Model~3 and largest for Model~5. Concavities of the R\'enyi functions become smaller when $b$ increases and the other functions have almost linear behaviours for $q \in (0,3)$ for large values of $b$.\\[-5mm]
\section{Numerical studies}
\label{S1:8}
\subsection{Simulation methodology}
\label{S1:81}

There are numerous models for which explicit expressions for the R\'enyi function in terms of elementary functions or even series are not available. Also, in the majority of cases obtaining explicit mathematical formulae for R\'enyi functions is a difficult problem and rigorous results were derived only for some ranges of the parameter~$q$. For example, for all models in \cite{leonenko2013} and this paper, $T(q)$ was derived for $q \in [1,2]$ only. These results are also likely to be true for wider ranges of $q$, but showing it requires new proof strategies, see Section~\ref{S1:5} and the example in Section~\ref{S1:7}. For such difficult cases, random field simulations can be used to obtain realizations of random fields from theoretical models and compute empirical R\'enyi functions. These empirical R\'enyi functions can be used as a substitute of exact mathematical functions for verifications whether real data are consistent with considered theoretical~models.

The R package RandomFields provides a wide range of simulation techniques and algorithms for random fields, see~\cite{randomfield} for more detail. The function RFSimulate simulates Gaussian random fields for a given covariance or variogram model and parameters defined by the arguments RMnugget and MTrend.

To compute a sample R\'enyi function the ratio \large{$\frac{\log_2 E \sum_{l} \mu(S_l^{(m)})^q}{\log_2|S_l^{(m)}|}$} \normalsize from \eqref{eq:6a} for large values of $m$ was used. This ratio was replaced by an empirical estimator that employs the HEALPix structure. For the HEALPix resolution $n=$~1024 with 12582912 pixels (in the following denoted by $i$), 196608 sets $S_l^{(m)}$ with 64 pixels per set were used to estimate the ratio. As CMB measurements $M(i)$ can take negative values they were transformed to non-negative ones by subtracting their minimum value: $\tilde{M}(i)$=$M(i)-\min{(M(i))}$. Then, the terms $E(\mu( S_l^{(m)}))^q$ were estimated by ${\hat{\mu_{l}}}^q$=${\left(\sum_{i \in S_l^{(m)}}{{\tilde{M}(i)}}\mathbin{/}{\sum {\tilde{M}(i)}}\right)}^q$. Finally, the empirical R\'enyi function was computed as $$\hat{T}(q) = \frac{\log_2(\sum_{l}{\hat{\mu_{l}}}^q)}{\log_2|S_l^{(m)}|}.$$

   The empirical multifractal spectrum was estimated by
\[\hat{f}(\alpha)= q \cdot \hat{\alpha}(q)-\hat{T}(q),\quad \quad \hat{\alpha}(q) = \frac{\sum_l \left(\left({\hat{\mu_l}^q}\mathbin{/}{\sum \hat{\mu_l}^q}\right) \cdot \ln(\hat{\mu_l})\right)}{\log_2|S_l^{(m)}|}.\]

For the selected large number of pixels, $\hat{T}(q)$ and $\hat{f}(q)$ provide reliable estimations and can be used for wider intervals of $q$ values than $[1, 2]$. In this paper we considered intervals $(0.5, 3)$ and $(-10, 10)$ when it was required.

\begin{figure}[!htb]
    \centering\vspace{-1cm}
    \subfloat[Realization of a multifractal random field]{\label{fig14a}
    \includegraphics[width=0.32\textwidth, height=0.20\textheight]{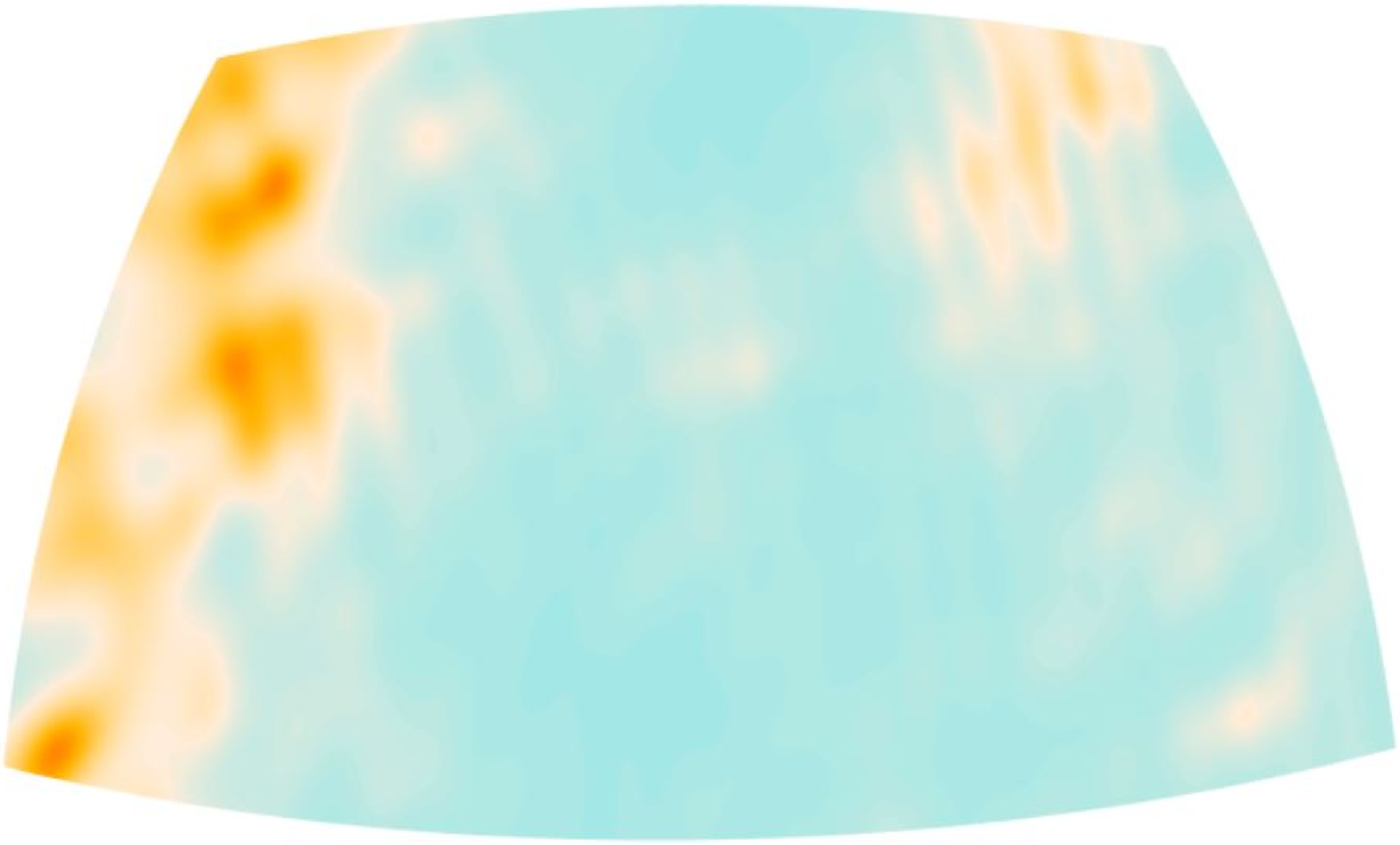}}
    \hfill
    \centering
    \subfloat[Sample R\'enyi function with the fitted Log-Normal model]{\label{fig14b}
    \includegraphics[width=0.32\textwidth, height=0.20\textheight]{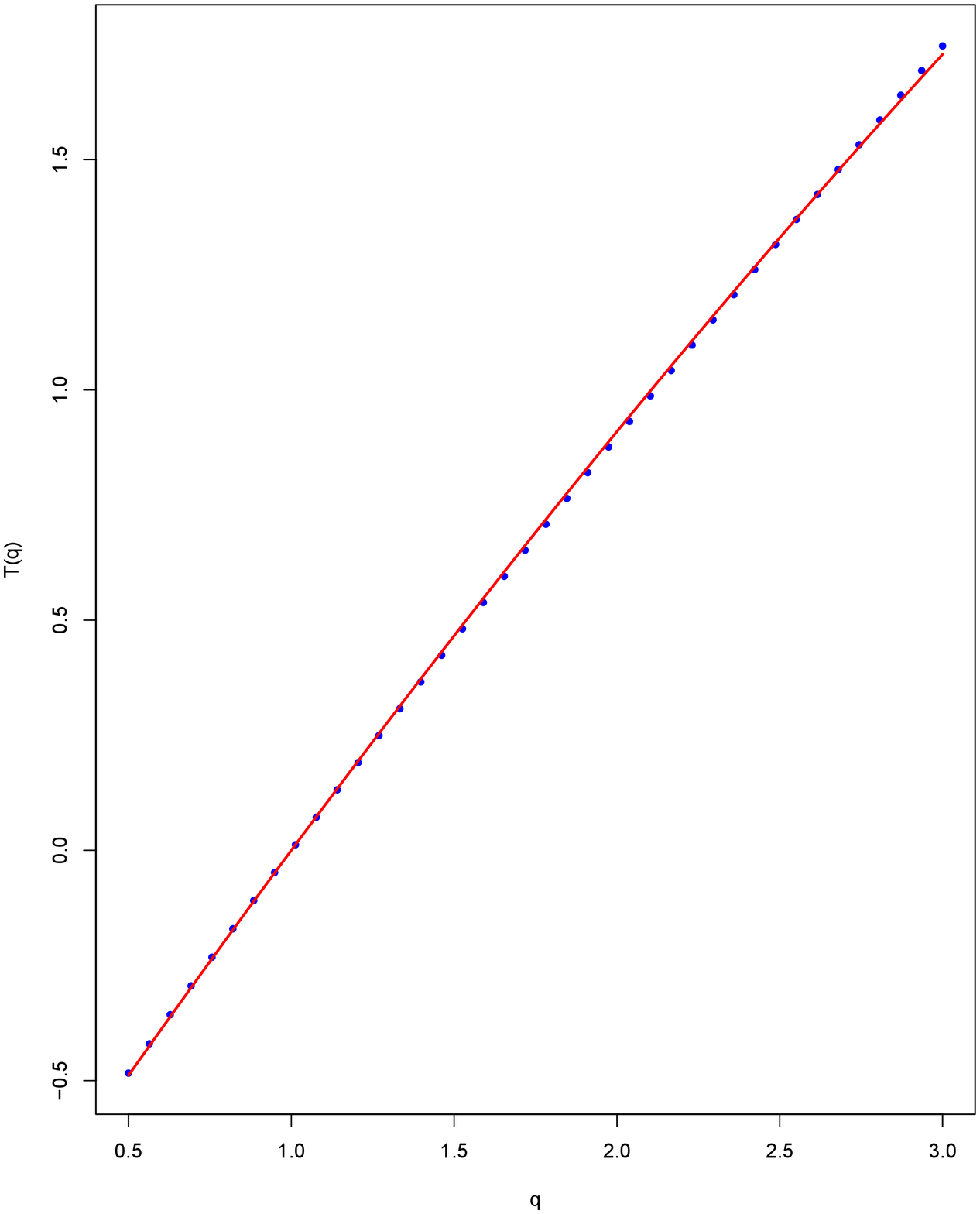}}
    \hfill
    \centering
    \subfloat[$f(\alpha)$ versus $\alpha$]{\label{fig14c}
    \includegraphics[width=0.32\textwidth, height=0.20\textheight]{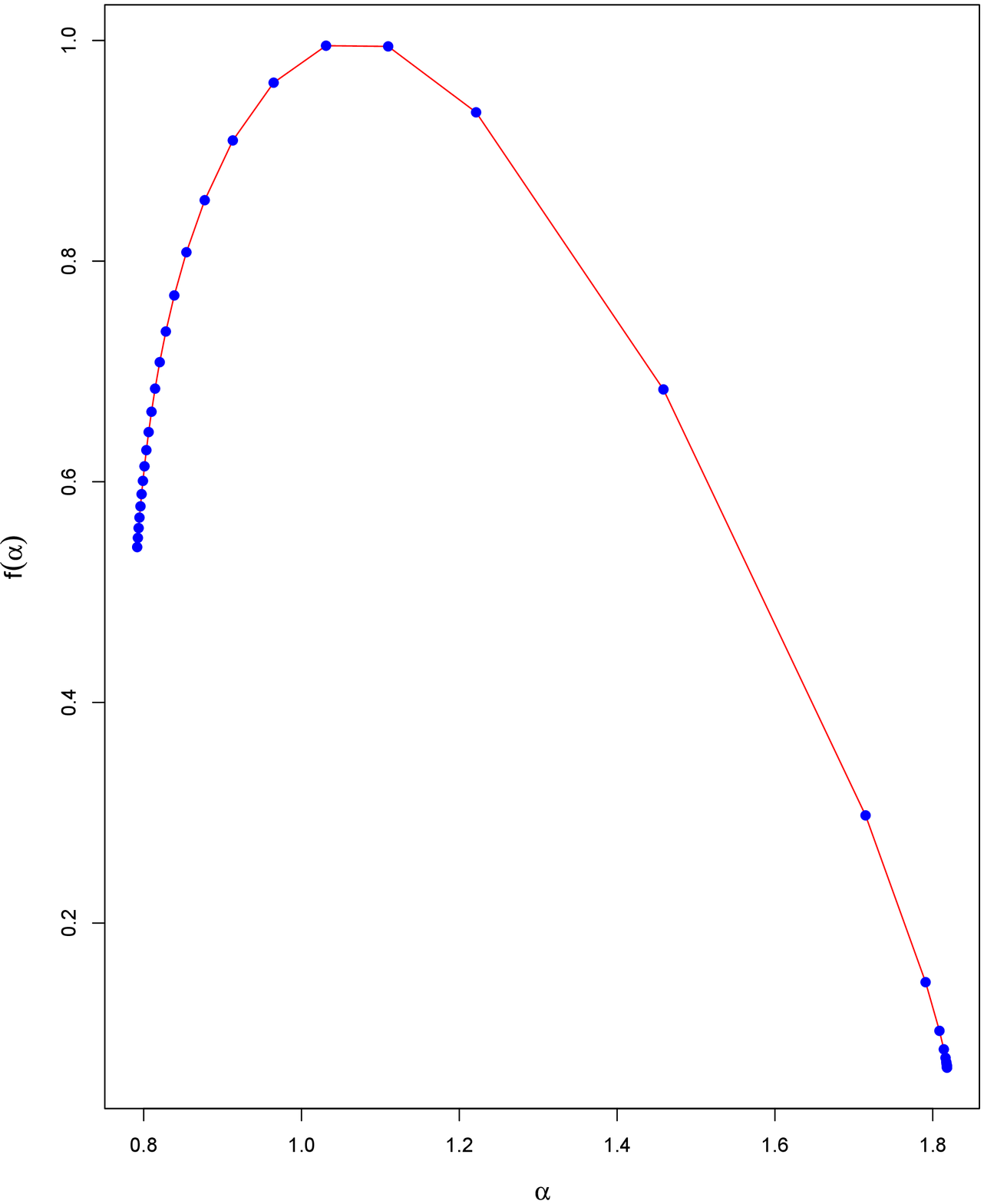}}
    \caption{Analysis of a simulated multifractal random field}\label{fig14}
\end{figure}

Fig.~\ref{fig14a} shows a realization of a multifractal random field in a large spherical window. The field was obtained from a Gaussian mother random field $Y(x)$ with the exponential covariance model and its variance equals 2. This covariance function has an exponential form and obviously satisfies the inequality in~\eqref{eq:4a}. As an approximation of the limit field, a finite product field $\Lambda_{40}(x)$ with $b=3$ was~used.

First 40 realizations $Y_{i}(b^{i}x), \; i=1,...,40,$ were simulated on a sphere by the package~RandomFields. Then the finite-product field ${\Lambda}_{40}(x)$ was computed by transforming the simulated values according to the formula $\exp(\sum_{i=1}^{40}Y_{i}(b^{i}x)-40)$. The dot plot of the empirical  R\'enyi function is shown in Fig.~\ref{fig14b}. The solid straight line is used as a reference to see departures from the fitted Model~1. It is clear that the empirical R\'enyi function and the theoretical one from \eqref{eq:8a} are very close on an interval that is wider than $[1,2]$. Fig.~\ref{fig14c} shows the spread of the multifractal spectrum.

The simulation studies suggest that the theoretical results from previous sections also hold for intervals wider than in the theorems. 

\subsection{Computing the R\'enyi function for CMB data}
\label{S1:82}

\begin{figure}[!htb]
    \centering\vspace{-0.8cm}
    \subfloat[Large window]{\label{fig15a}
    \includegraphics[width=0.32\textwidth, height=0.20\textheight]{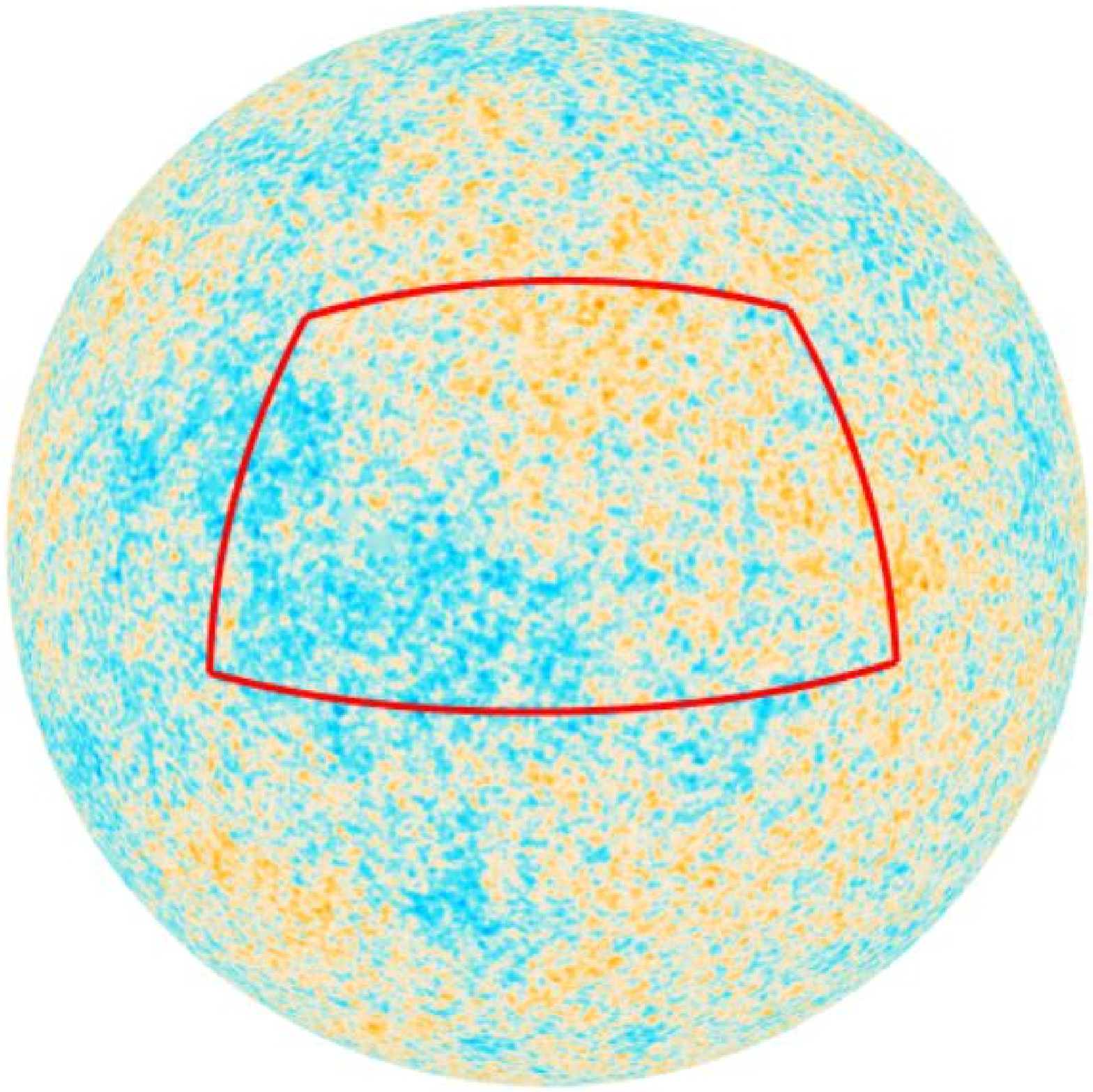}}
    \hfill
    \centering
    \subfloat[Medium window]{\label{fig15b}
    \includegraphics[width=0.32\textwidth, height=0.20\textheight]{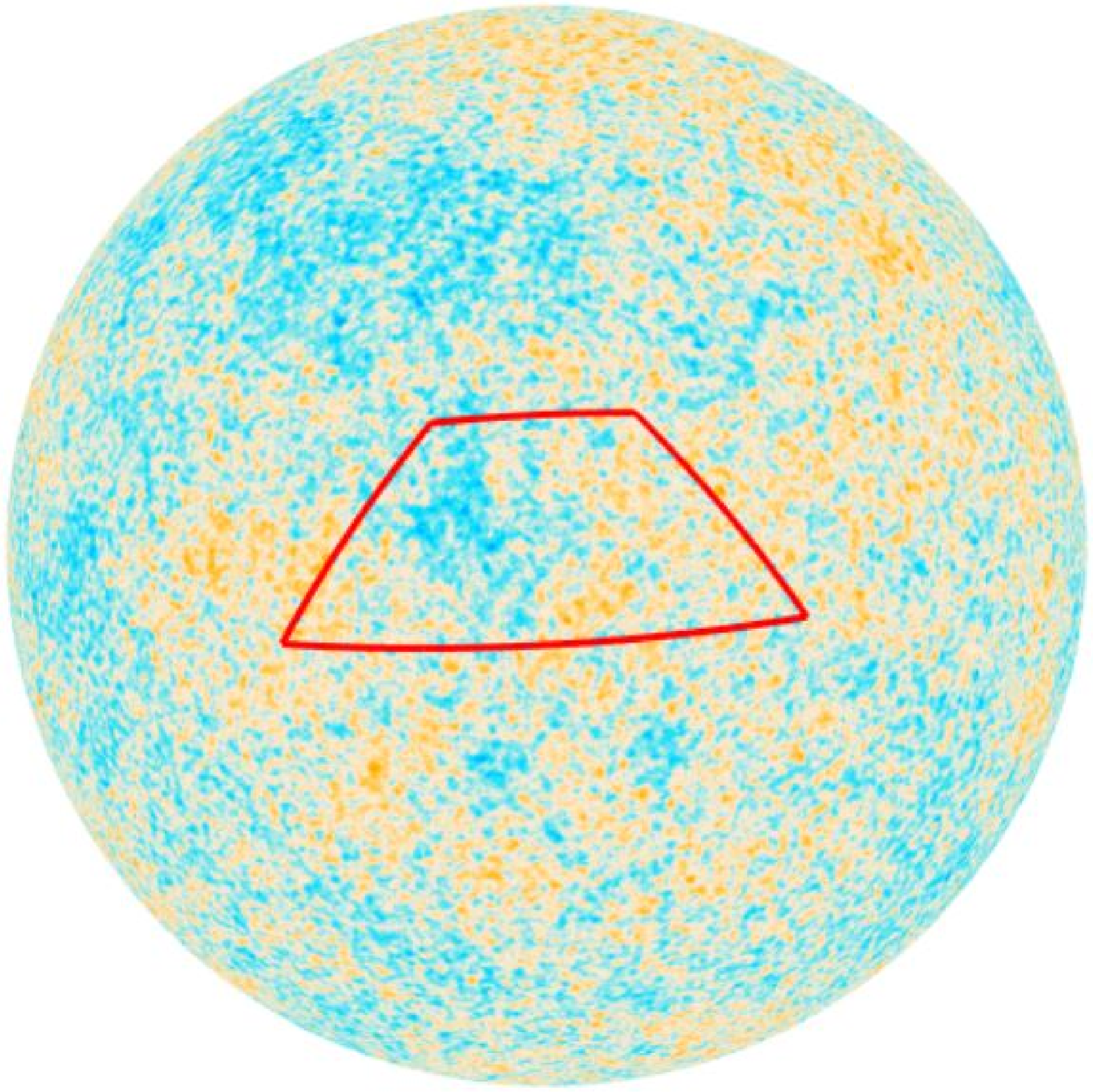}}
    \hfill
    \centering
    \subfloat[Small and very small windows]{\label{fig15c}
     \includegraphics[width=0.32\textwidth, height=0.20\textheight]{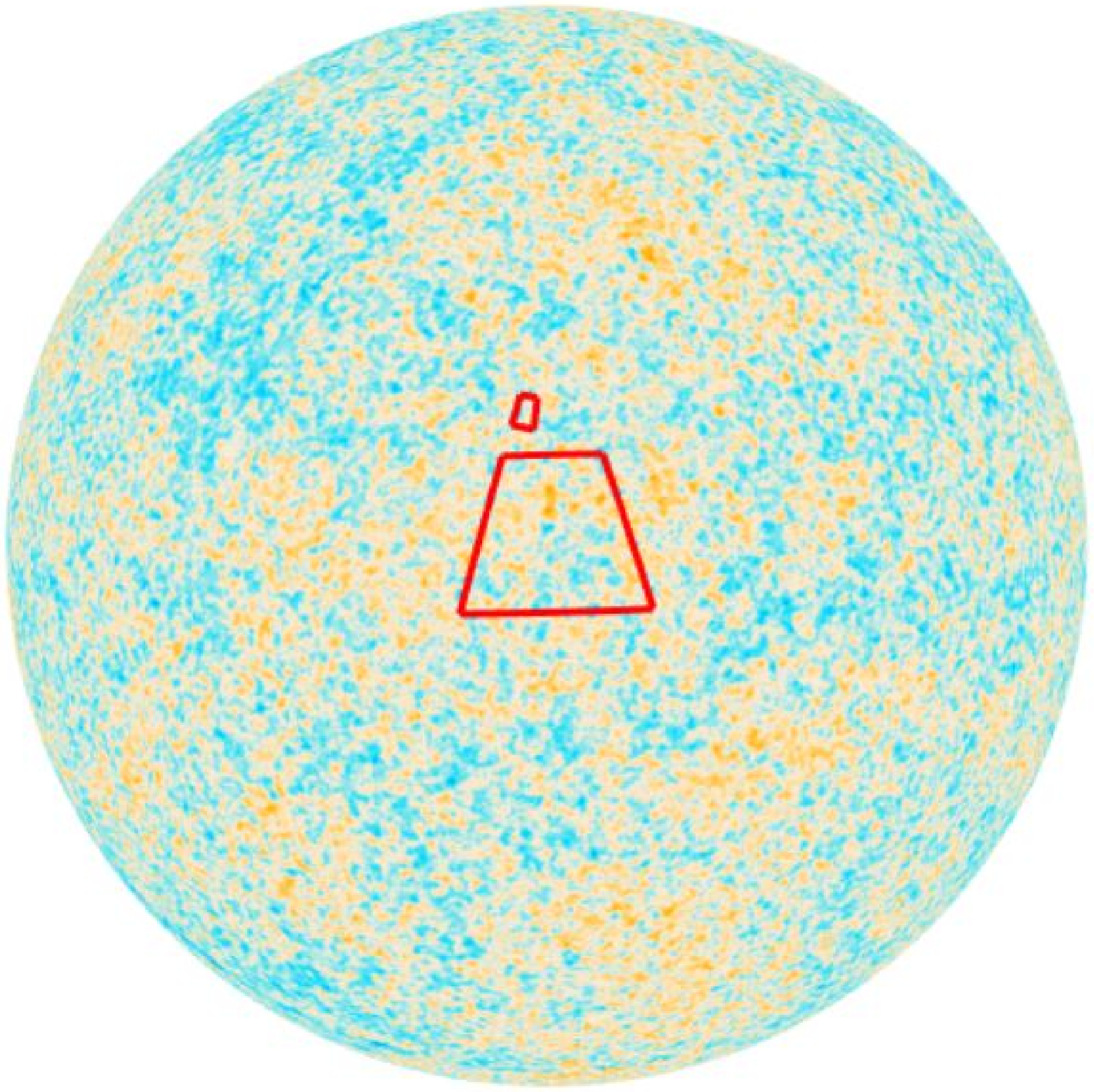}}
    \caption{Different sky windows of CMB data}\label{fig15}
\end{figure}

In this section, empirical R\'enyi functions were calculated for real cosmological data obtained from the NASA/IPAC Infrared Science Archive {\cite{[dataset]}}. Fig.~\ref{fig15} gives examples of sky windows CMB data from which were used to get empirical R\'enyi functions in the following examples. 
\begin{figure}[!htb]
    \centering\vspace{-0.5cm}
    \subfloat[Sample R\'enyi function versus linear function]{\label{fig16a}
    \includegraphics[width=0.32\textwidth, height=0.20\textheight]{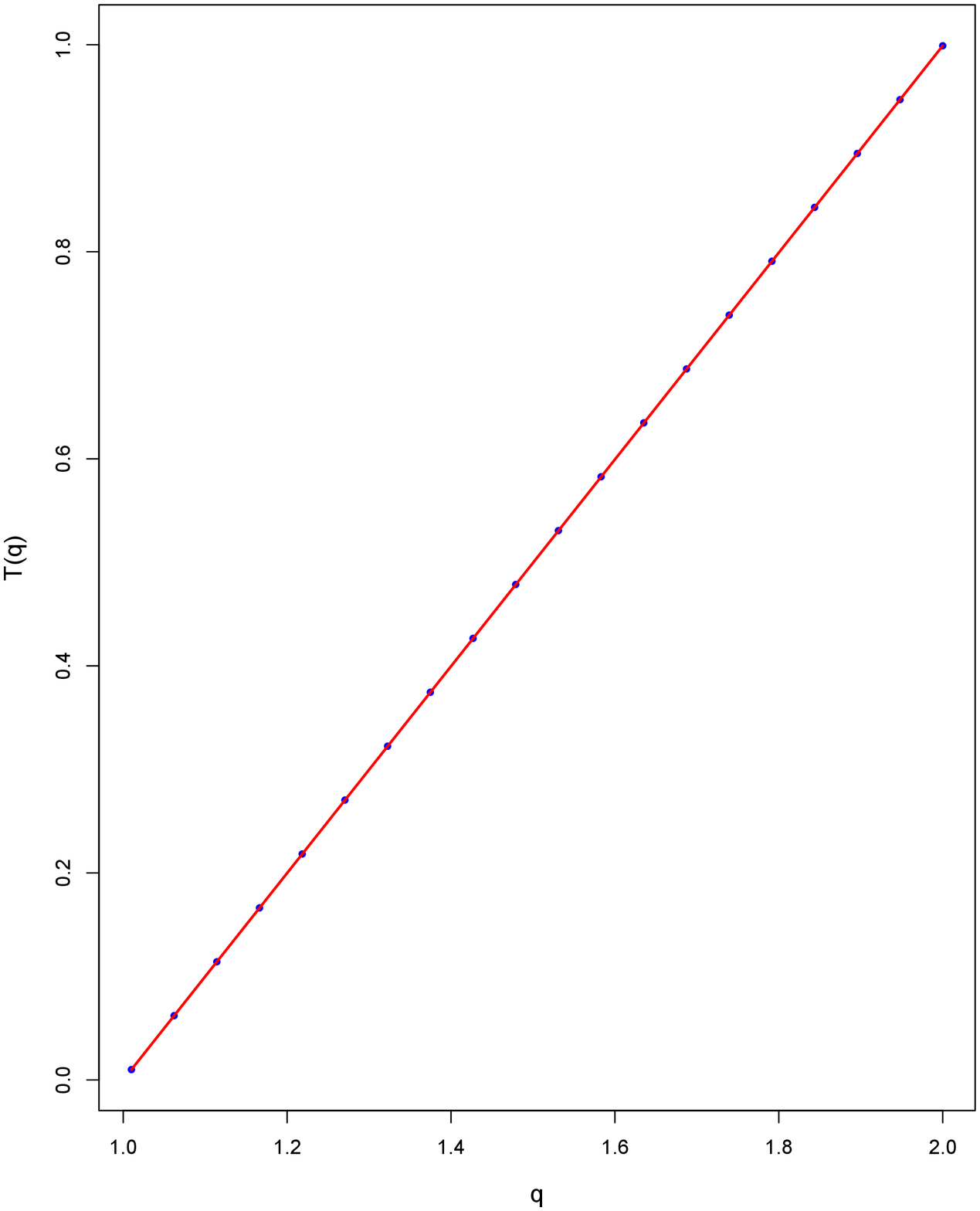}}
    \hfill
    \centering
    \subfloat[Difference of sample R\'enyi function and linear function]{\label{fig16b}
    \includegraphics[width=0.32\textwidth, height=0.20\textheight]{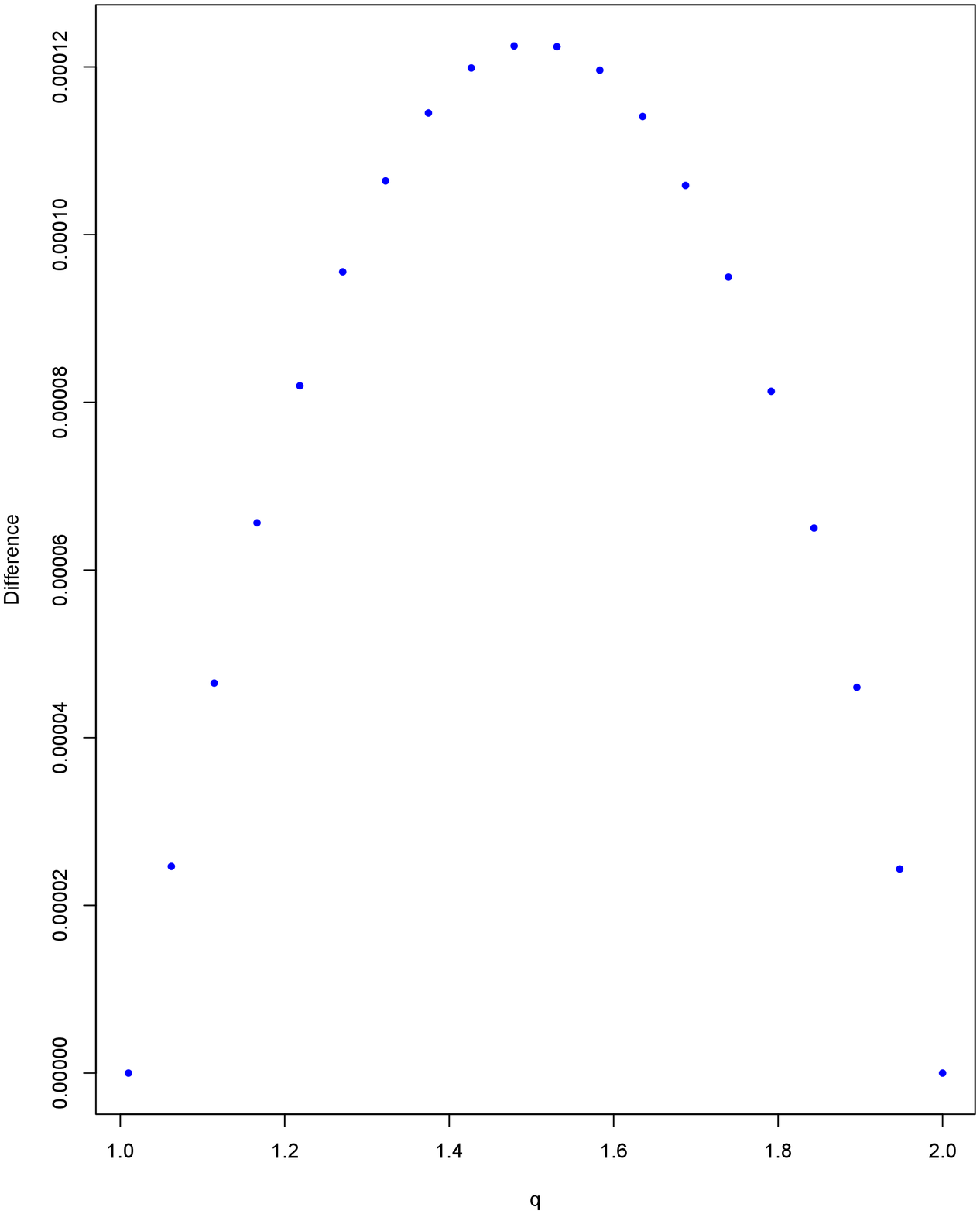}}
    \hfill
    \centering
    \subfloat[$\alpha(q)$ versus $q$]{\label{fig16c}
    \includegraphics[width=0.32\textwidth, height=0.20\textheight]{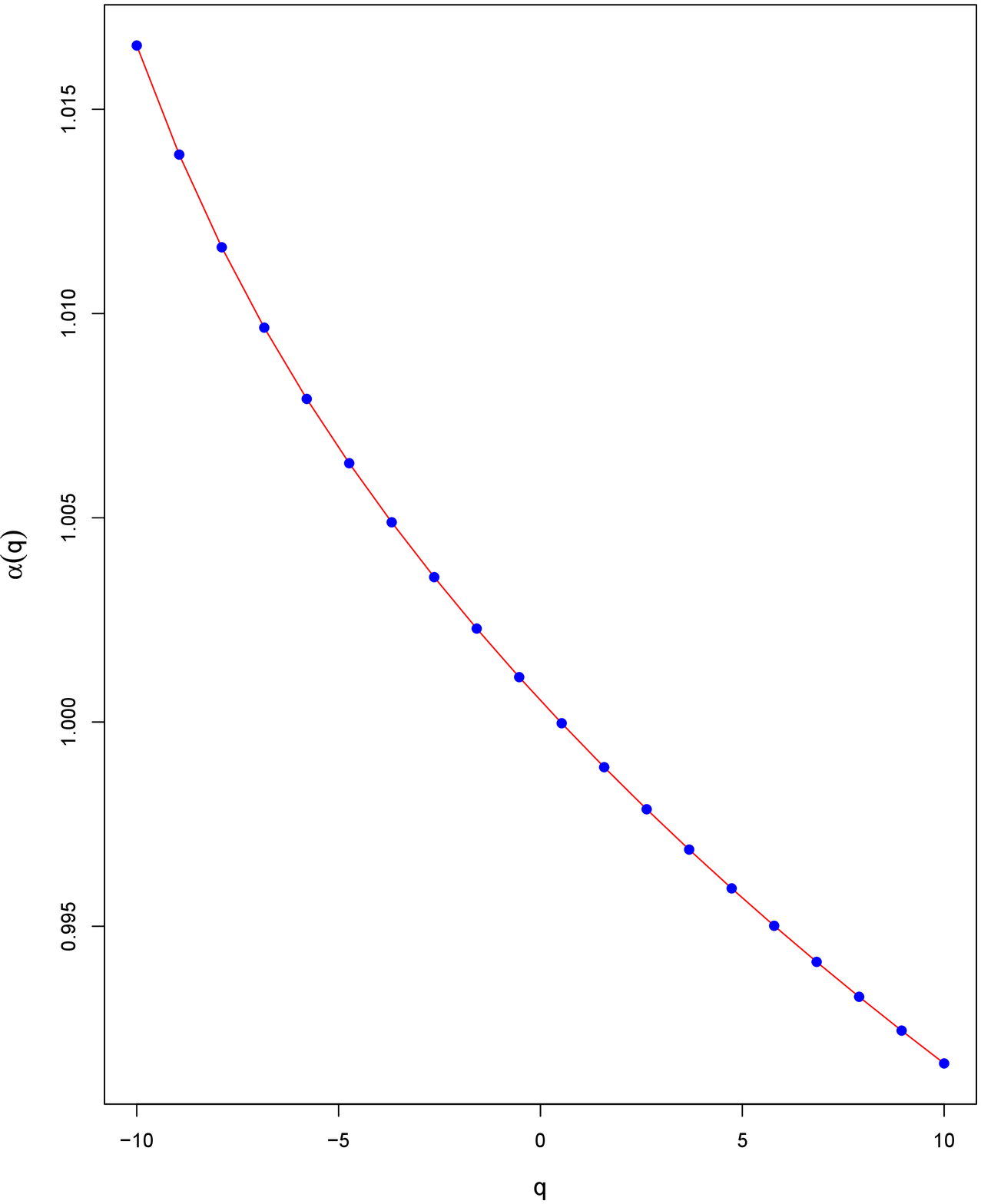}}
    \medskip \\
    \centering
    \subfloat[$f(\alpha)$ versus $\alpha$]{\label{fig16d}
    \includegraphics[width=0.32\textwidth, height=0.20\textheight]{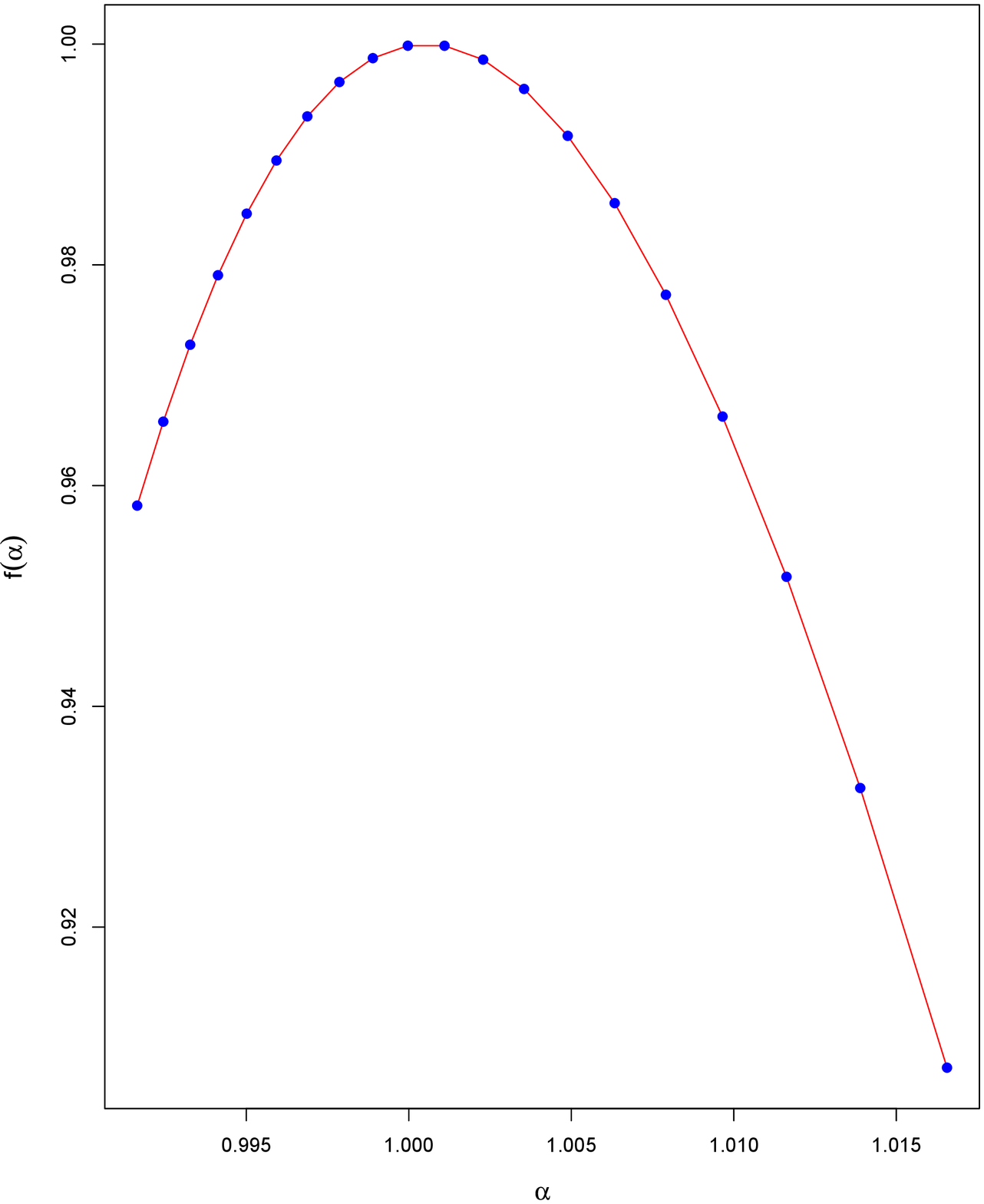}}
    \hfill
    \centering
    \subfloat[Sample R\'enyi function with the fitted Log-Normal model]{\label{fig16e}
    \includegraphics[width=0.32\textwidth, height=0.20\textheight]{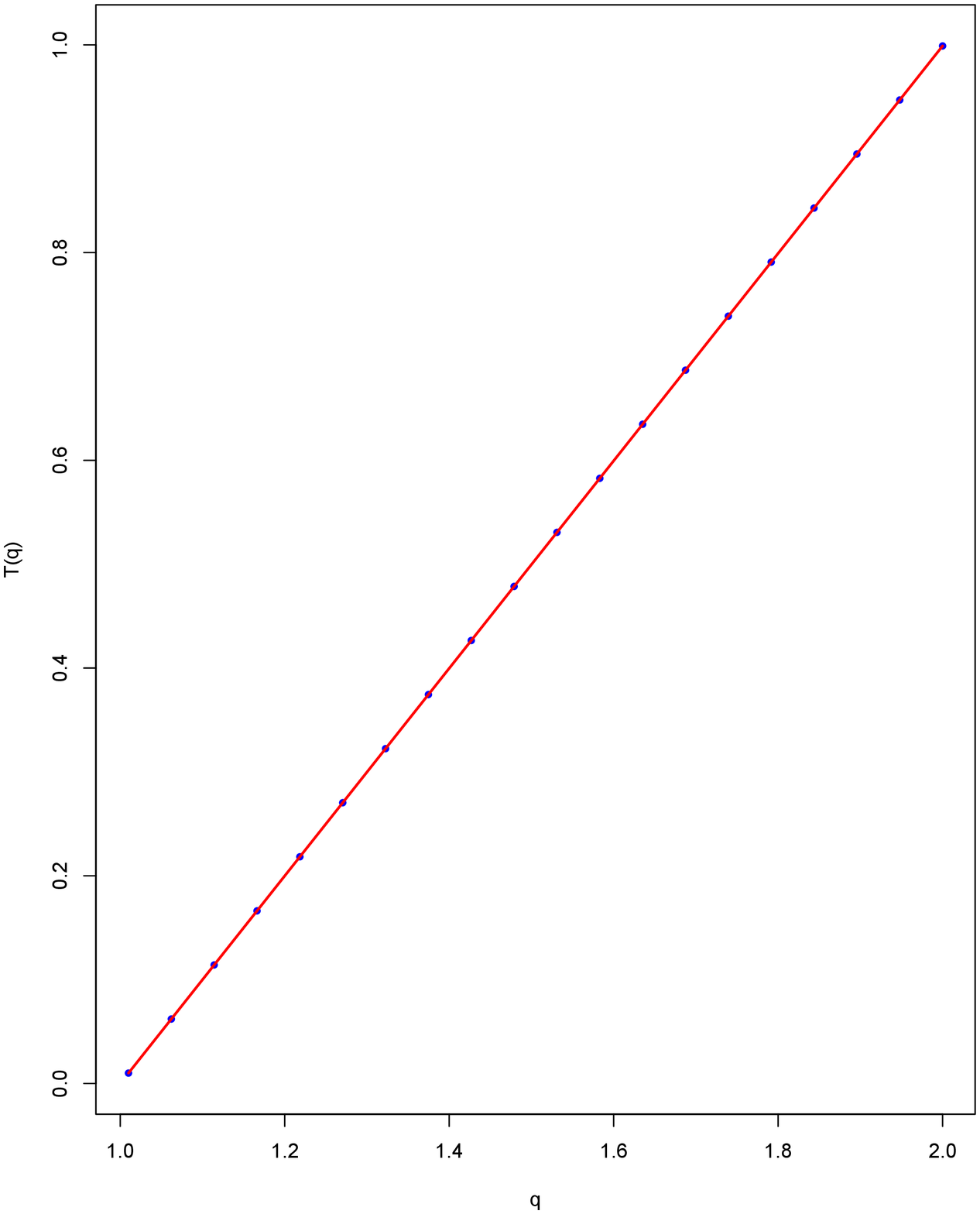}}
    \hfill
    \centering
    \subfloat[Difference between sample R\'enyi function and the fitted Log-Normal model]{\label{fig16f}
    \includegraphics[width=0.32\textwidth, height=0.20\textheight]{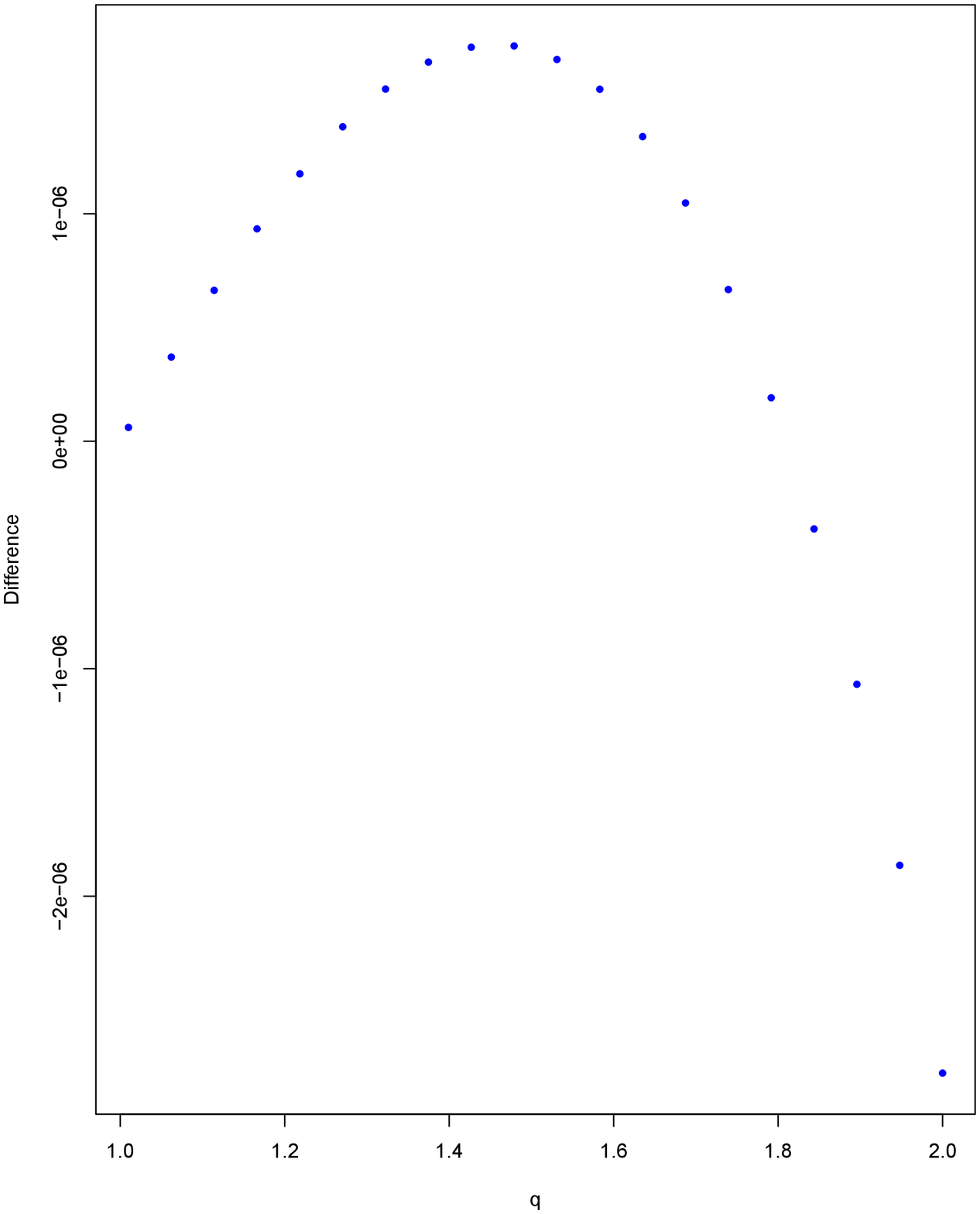}}
    \caption{Whole sky data analysis}
    \label{fig16}
\end{figure}

Extensive numerical studies were conducted for different windows in various sky locations. As in all cases we obtained rather similar results, we restrict our presentation only to few typical examples. The R package rcosmo was used for computations and visualizations, see \cite{fryer2019rcosmo, rcosmo1} for more details. For small windows the function fRen was slightly modified to change the support of the measure $\mu$ from the whole sky to the selected window.

\begin{figure}[!htb]
    \centering\vspace{-1cm}
    \subfloat[$f(\alpha)$ versus $\alpha$ for large window]{\label{fig17a}
    \includegraphics[width=0.32\textwidth, height=0.20\textheight]{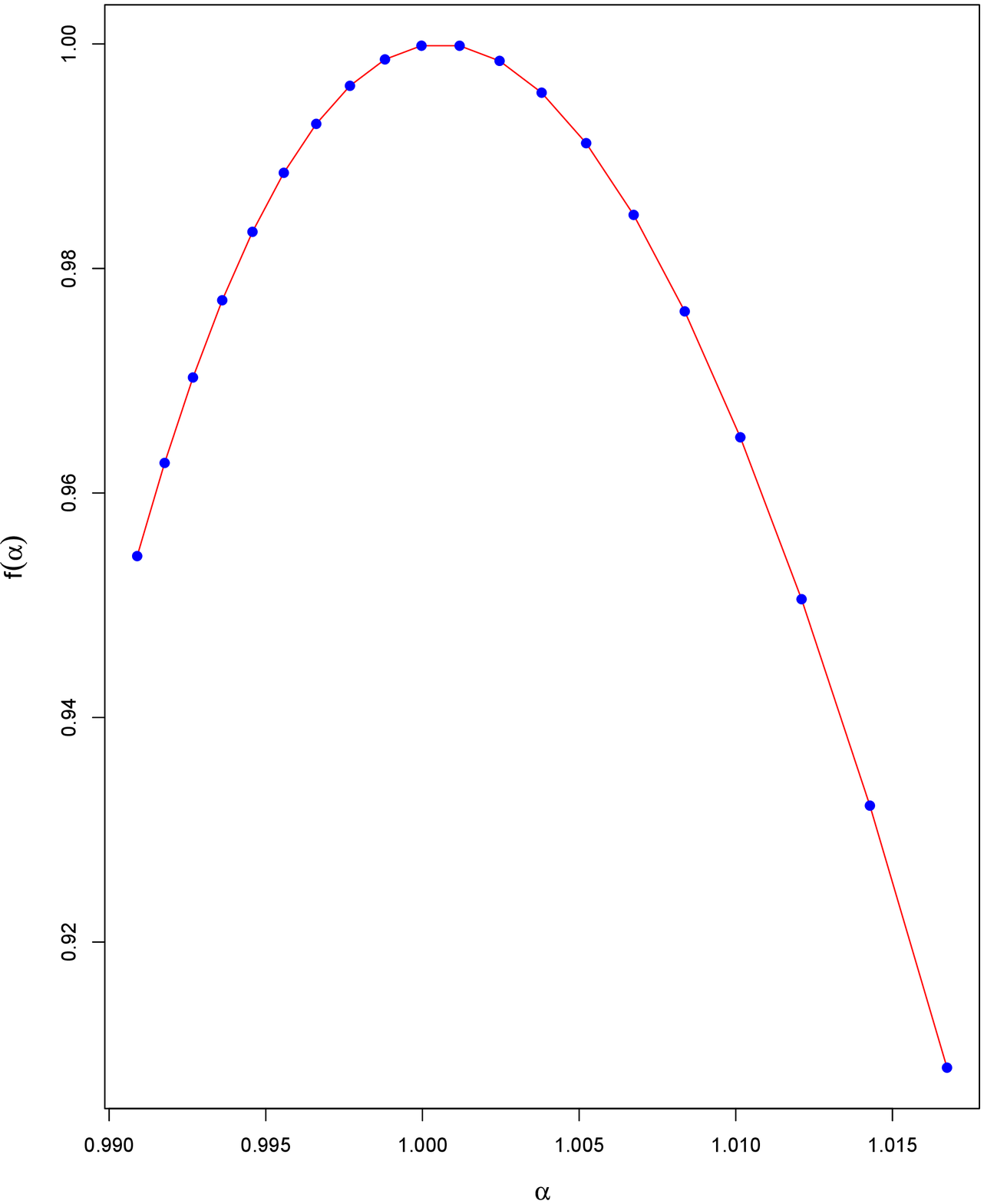}}
    \hfill
    \centering
    \subfloat[Difference with linear function for large window]{\label{fig17b}
    \includegraphics[width=0.32\textwidth, height=0.20\textheight]{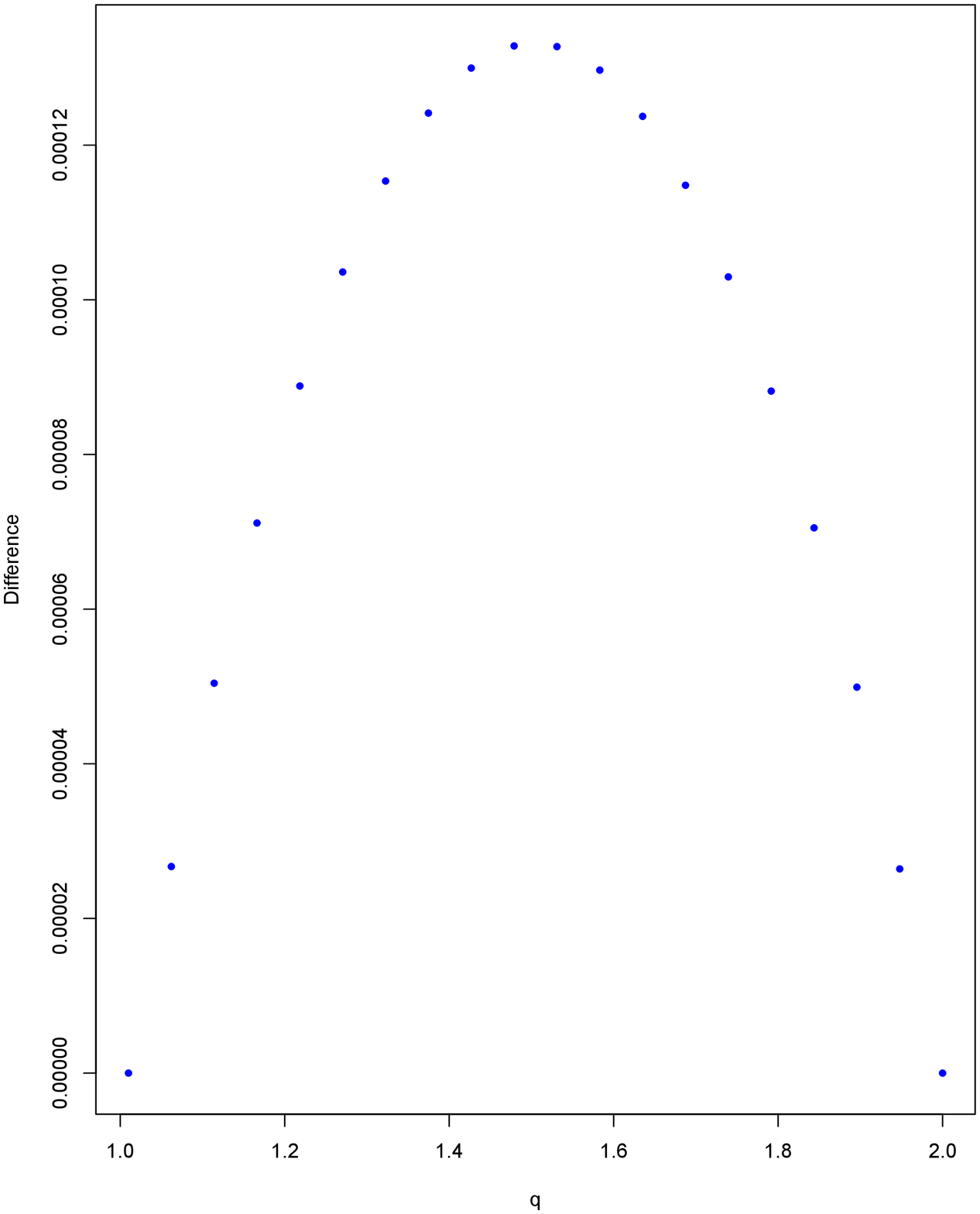}}
    \hfill
    \centering
    \subfloat[Difference with Model 1 for large window]{\label{fig17c}
    \includegraphics[width=0.32\textwidth, height=0.20\textheight]{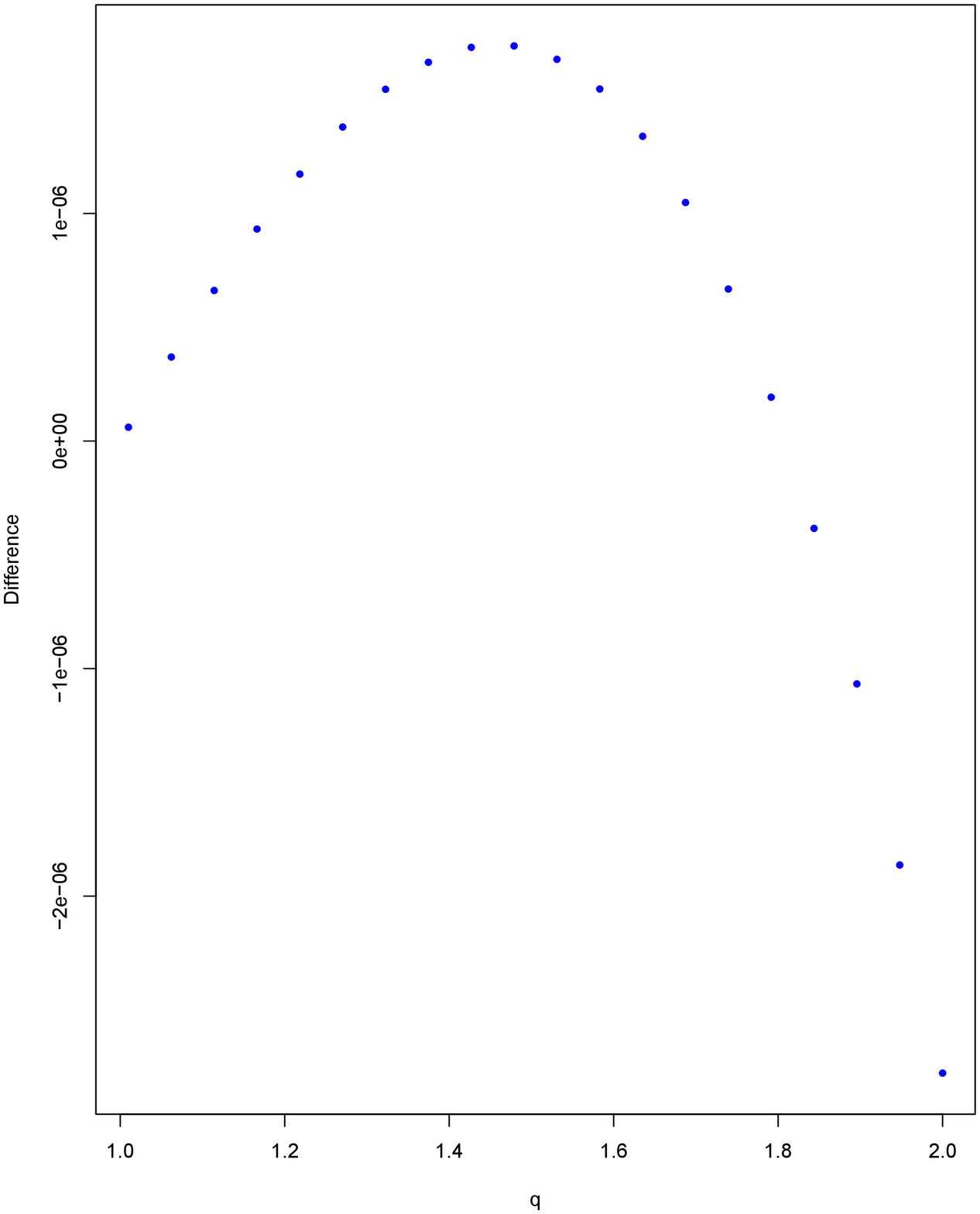}}
    \medskip \\
    \centering
    \subfloat[$f(\alpha)$ versus $\alpha$ for small window]{\label{fig17d}
    \includegraphics[width=0.32\textwidth, height=0.20\textheight]{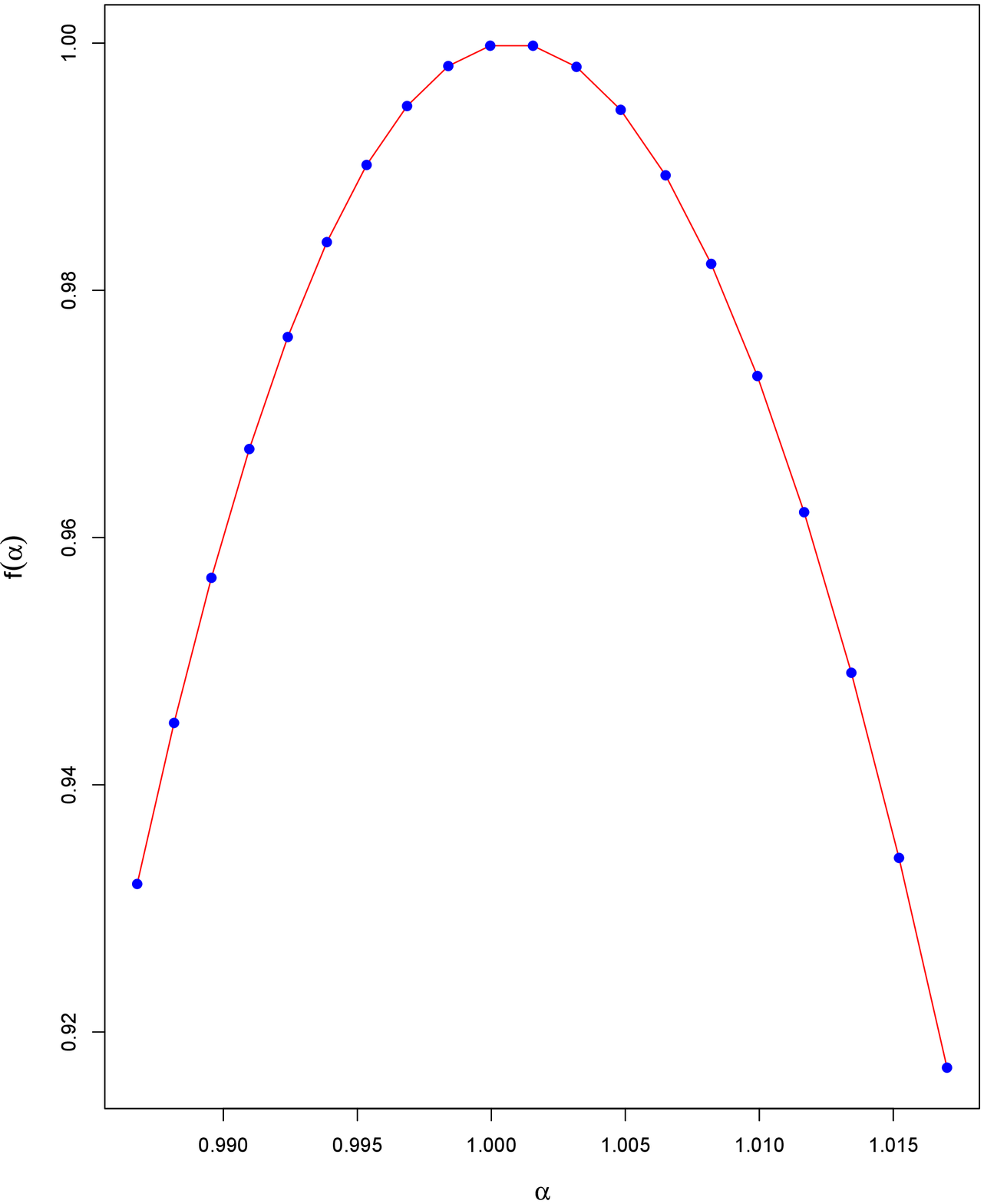}}
    \hfill
    \centering
    \subfloat[Difference with linear function for small window]{\label{fig17e}
    \includegraphics[width=0.32\textwidth, height=0.20\textheight]{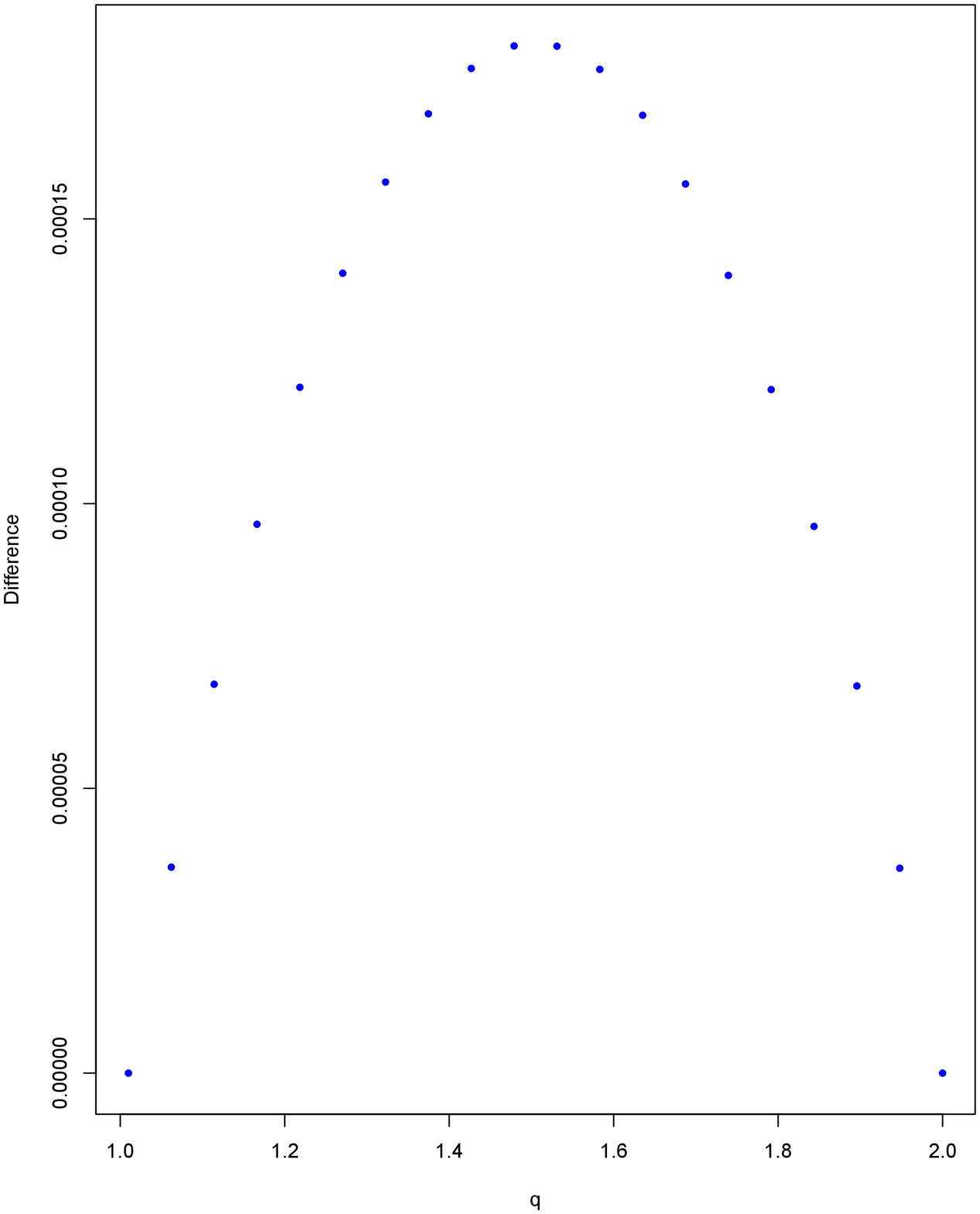}}
    \hfill
    \centering
    \subfloat[Difference with Model 1 for small window]{\label{fig17f}
    \includegraphics[width=0.32\textwidth, height=0.20\textheight]{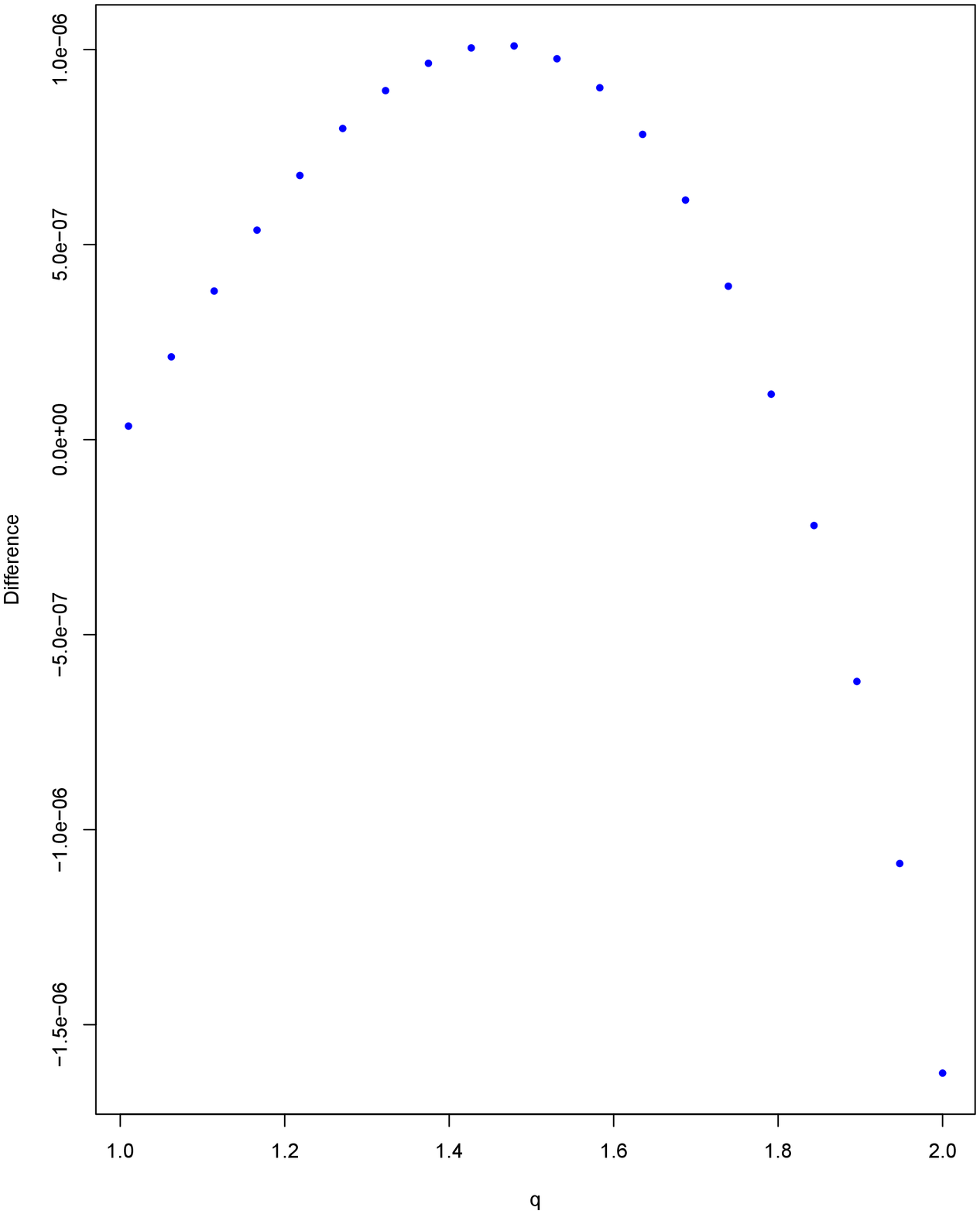}}
    \caption{Analysis of large and small sky windows data}
    \label{fig17}
\end{figure}

First the R\'enyi function was computed for the whole sky. The obtained sample R\'enyi function is shown in Fig.~\ref{fig16a} by dots. The straight line in the Fig.~\ref{fig16a} was drawn to assess departures of the sample R\'enyi function from a linear behaviour. The difference of the sample R\'enyi function and the linear function that connects the points (1,0) and (2,1) is shown in Fig.~\ref{fig16b}. It is clear that the departure from a  linear behaviour is not substantial. Fig.~\ref{fig16c} shows the function $\alpha(q)$ and Fig.~\ref{fig16d} plots the function $f(\alpha(q))$ versus $\alpha(q)$. As it was discussed in Section~\ref{S1:81} to compute $\alpha(q)$ and $f(\alpha(q))$, we used the formula for R\'enyi functions for the range $(-10,10)$ as simulation studies and analysis in~\cite{grahovac2014detecting} suggest the same analytical form of the R\'enyi function as for the range~$[1,2]$. All these plots confirm only very small multifractality of the CMB data. Similar results were also obtained for different sky windows, see, for example, Fig.~\ref{fig17a}, Fig.~\ref{fig17b}, Fig.~\ref{fig17d} and Fig.~\ref{fig17e}. 

The R\'enyi functions, multifractal spectra, similar analysis and plots were produced for different window sizes of the CMB unit sphere. Large, medium, small and very small window sizes with areas 1.231, 0.4056, 0.0596 and 0.0017 were selected, see Fig.~\ref{fig15}. The R\'enyi function was computed for small windows located at different places of the sky sphere such as near the pole, near the equator and other places of the sphere. Although different window sizes of the sphere were investigated, there's not that much of evidence to suggest that we have substantial multifractality. The ranges of $y$ scale in Fig.~\ref{fig16b}, Fig.~\ref{fig17b} and Fig.~\ref{fig17e} suggest that this multifractality is very small. These results and the variations of the values of $\hat{T}(q)$ between windows suggest that collecting data at very fine scales and further tests of hypothesis are required. We plan to develop tests of hypothesis about R\'enyi functions in future publications and use them for new high resolution data that will be available from the next generation CMB experiments CMB-S4~\cite{abazajian2019cmbs4}. As the obtained plots are rather similar we present only two of them for large and small windows in Fig.~\ref{fig17}.

Then, all models from Sections~\ref{S1:6} and \ref{S1:7} were used to fit the empirical R\'enyi function. For the log-normal model we present the results for all windows. For other models, only results for CMB data in a large window are given. Similar results were also obtained for other windows. To fit models to empirical R\'enyi functions several methods were employed. For the log-normal model the simple linear regression approach was used whereas for the other models the non-linear regression approach was applied. 

As the R\'enyi function of the log-normal model is specified by \eqref{eq:8a}, substituting $a=\frac{\sigma_Y^2}{4\ln{b}}$, results in the form $T(q)= a(-q^2+q)+q-1$. Then the R function \enquote{lm} was used for a simple linear regression fit with the intercept~0 to $T(q)-q+1$. The values of the parameter $a$ and the root mean square error for deviations of Model 1 from the empirical R\'enyi function are given in~Table~\ref{table11}.

\begin{table}[!htb]
\caption{Analysis of different sky windows data with Model 1}
\label{table11}
\centering
\begin{tabular}{l l l l l}
\hline\noalign{\smallskip}
\textbf{\thead{Observation \\ window}} & \textbf{[$\alpha_{\min}$, $\alpha_{\max}$]} &  \textbf{$\alpha_{\max}$ - $\alpha_{\min}$} &  \textbf{\thead{a}}  & \textbf{RMSE}\\
\noalign{\smallskip}\hline\noalign{\smallskip}
\textbf{Whole Sky} & {[0.9916, 1.0165]} & 0.024917 & 0.000513 & {\SI{1.3602e-06}{}}\\
\textbf{Large} & {[0.9908, 1.0167]} & 0.025846 & 0.000555 & {\SI{1.3590e-06}{}}\\
\textbf{Medium} & {[0.9893, 1.0159]} & 0.026620 & 0.000629 & {\SI{1.1033e-06}{}}\\
\textbf{Small} & {[0.9867, 1.0170]} & 0.030219 & 0.000745 & {\SI{7.9095e-07}{}}\\
\textbf{Very Small} & {[0.9842, 1.0543]} & 0.070150 & 0.001500 & {\SI{1.3949e-05}{}} \\
\noalign{\smallskip}\hline
\end{tabular}
\end{table}

\begin{figure}[!htb]
    \centering\vspace{-1cm}
    \subfloat[Difference with Model 2]{\label{fig18a}
    \includegraphics[width=0.32\textwidth, height=0.20\textheight]{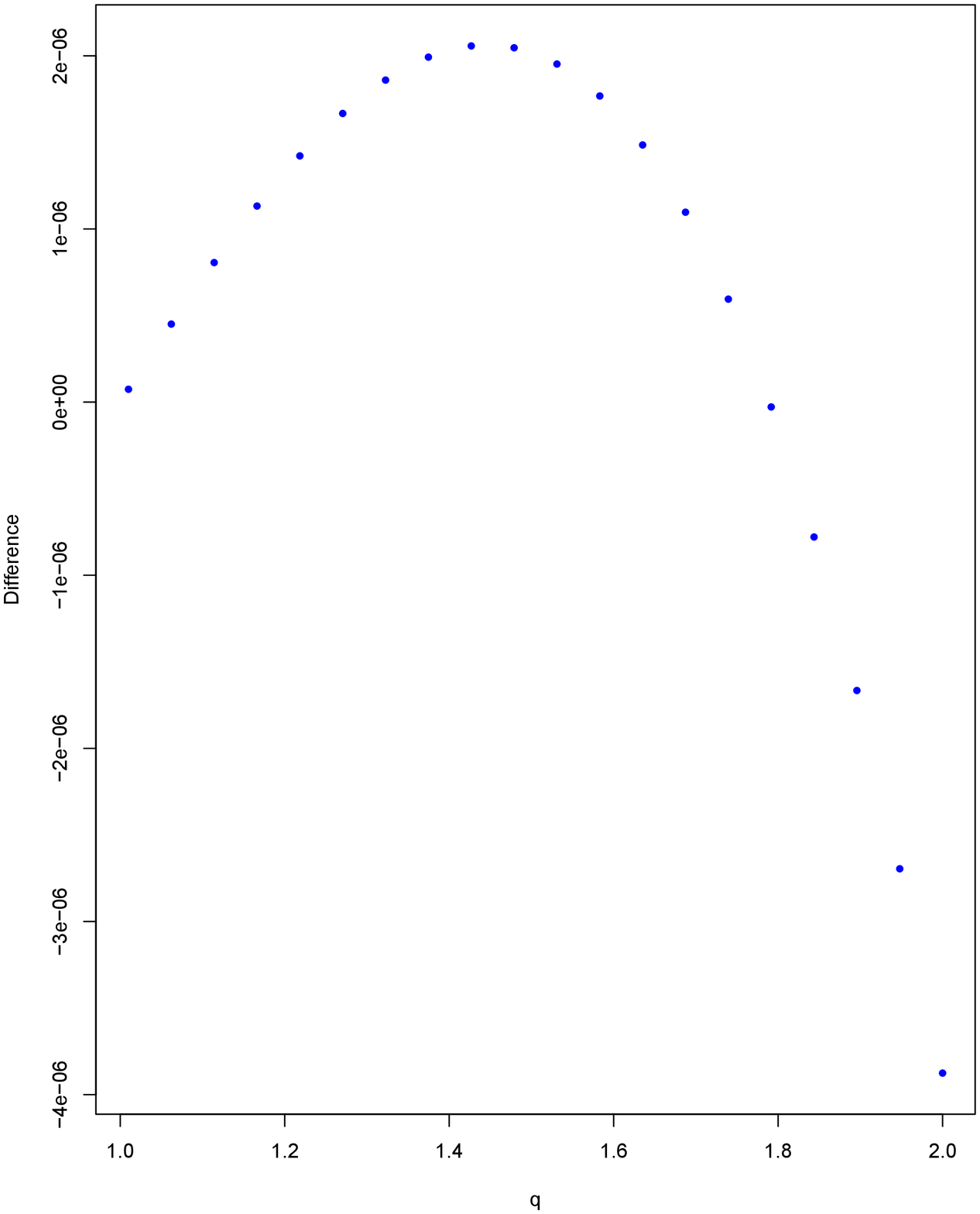}}
    \hfill
    \centering
    \subfloat[Difference with Model 3]{\label{fig18b}
    \includegraphics[width=0.32\textwidth, height=0.20\textheight]{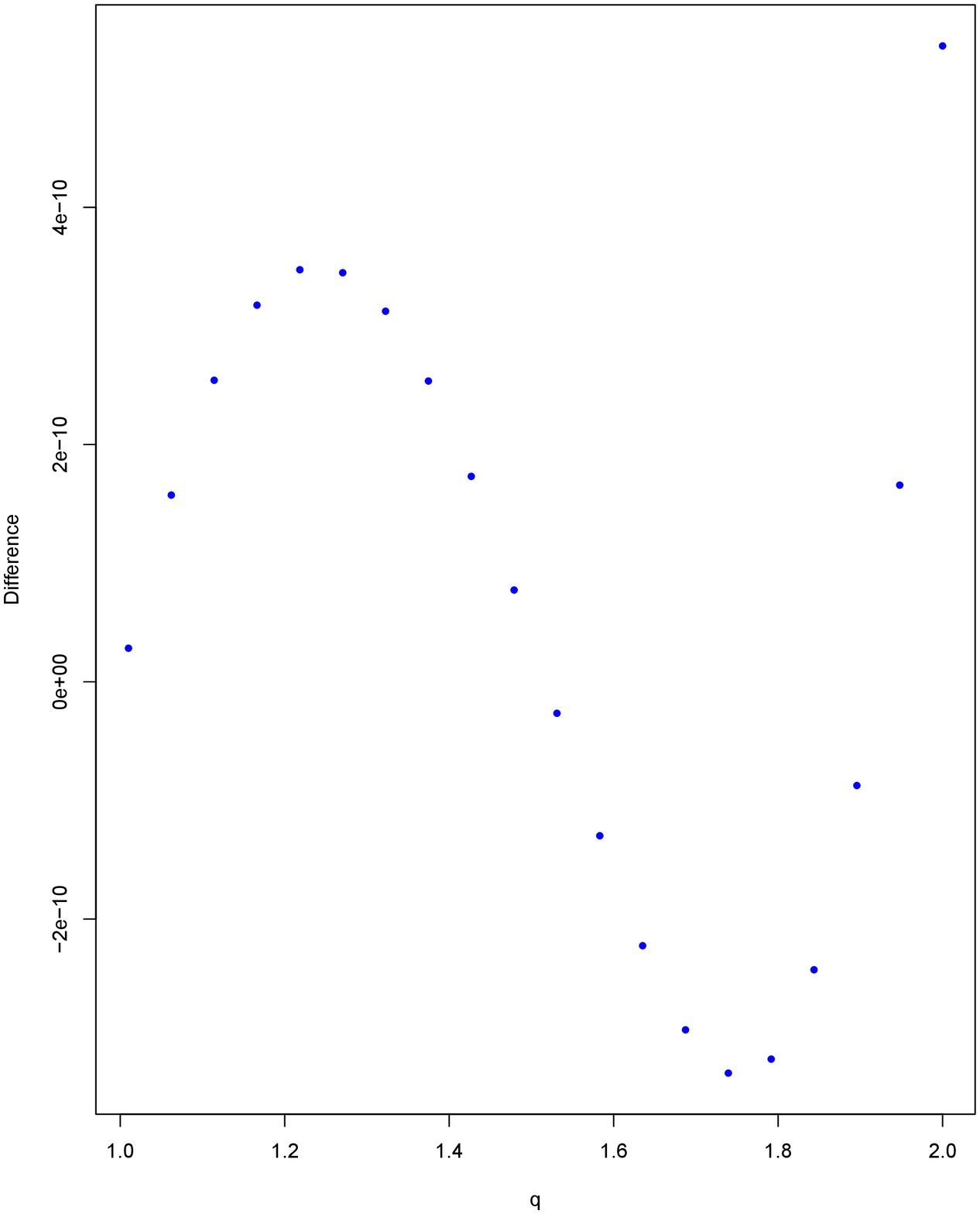}}
    \hfill
    \centering
    \subfloat[Difference with Model 4]{\label{fig18c}
    \includegraphics[width=0.32\textwidth, height=0.20\textheight]{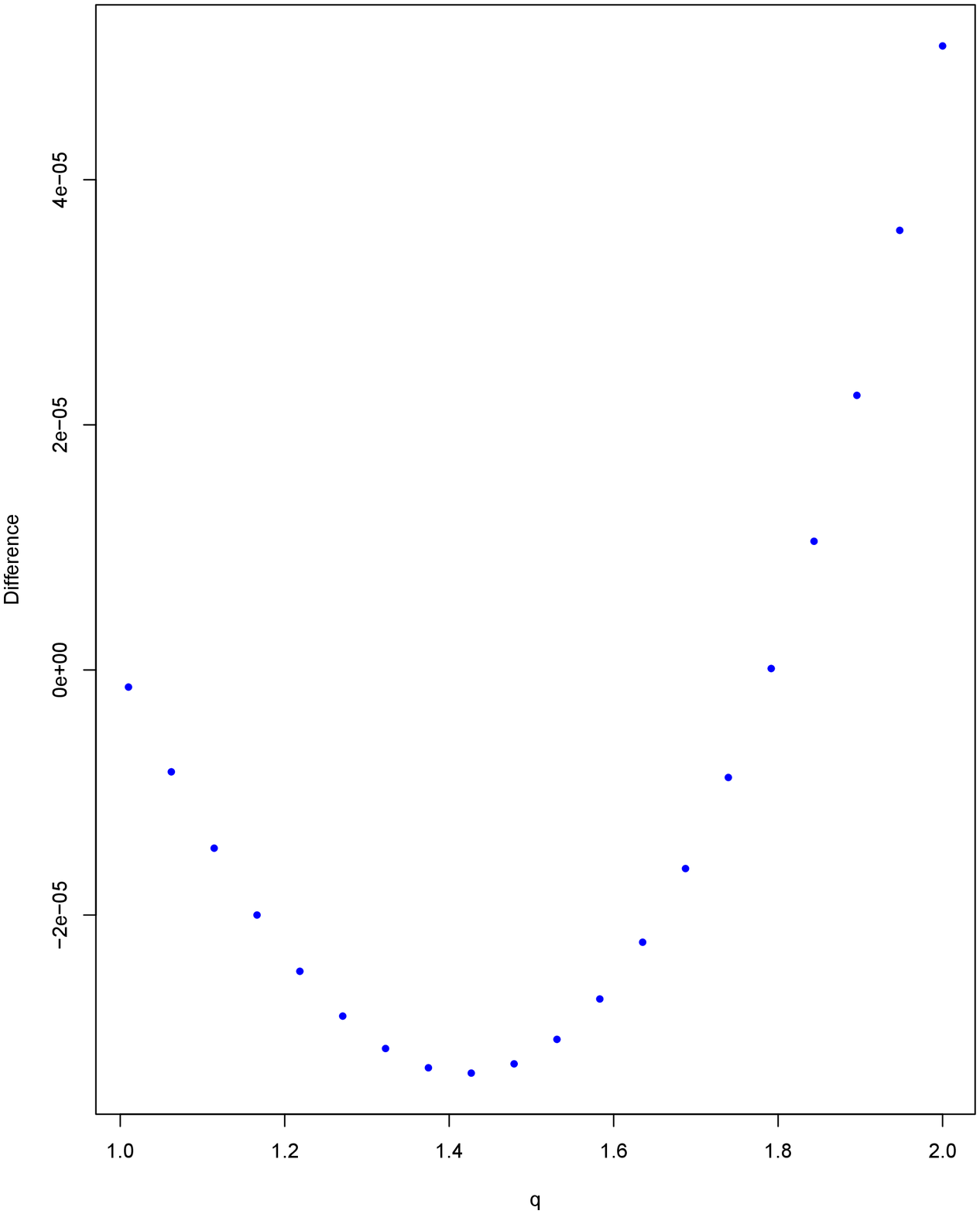}}
    \medskip \\
    \centering
    \subfloat[Difference with Model 5]{\label{fig18d}
    \includegraphics[width=0.32\textwidth, height=0.20\textheight]{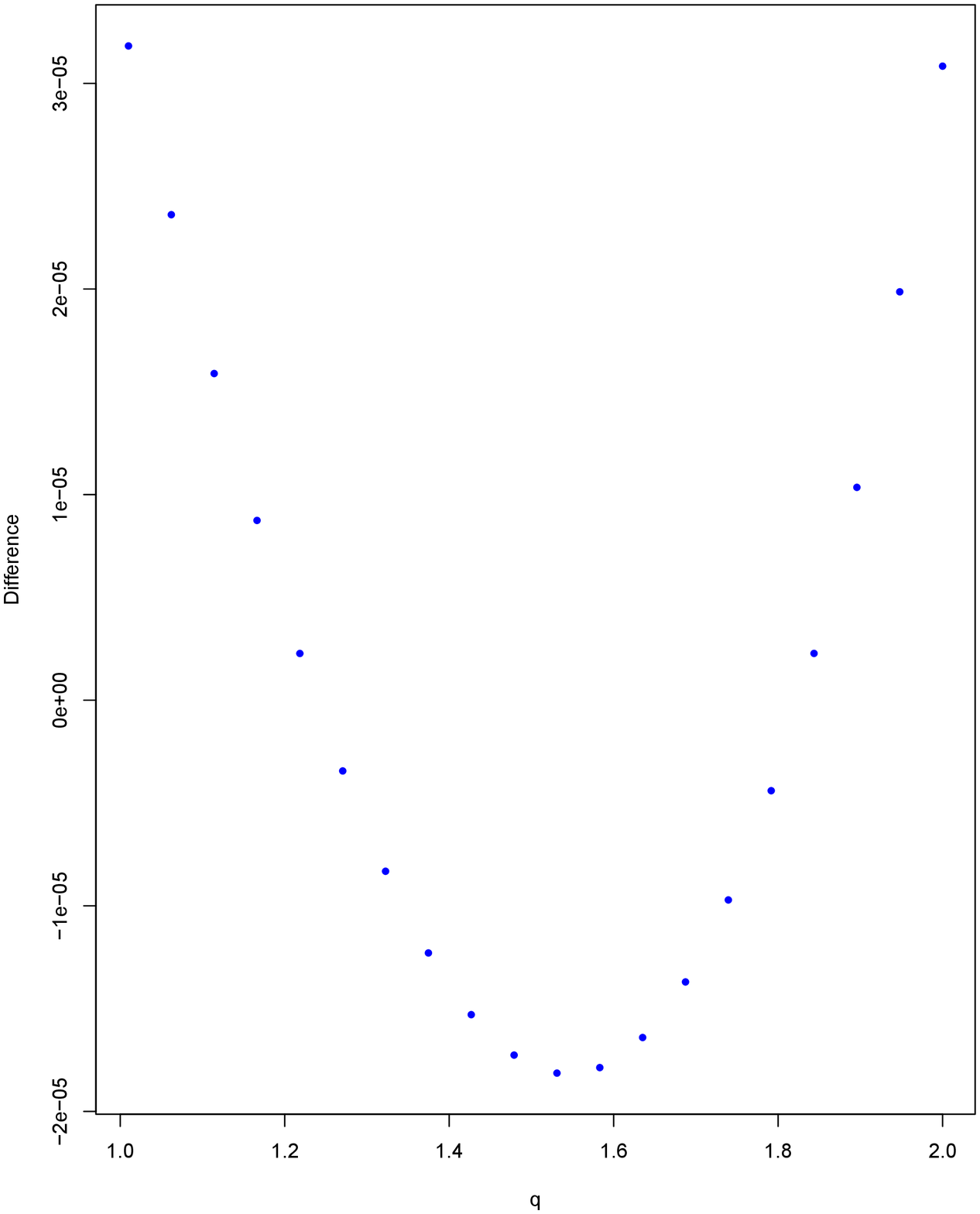}}
    \centering
    \subfloat[Difference with Model 6]{\label{fig18e}
    \includegraphics[width=0.32\textwidth, height=0.20\textheight]{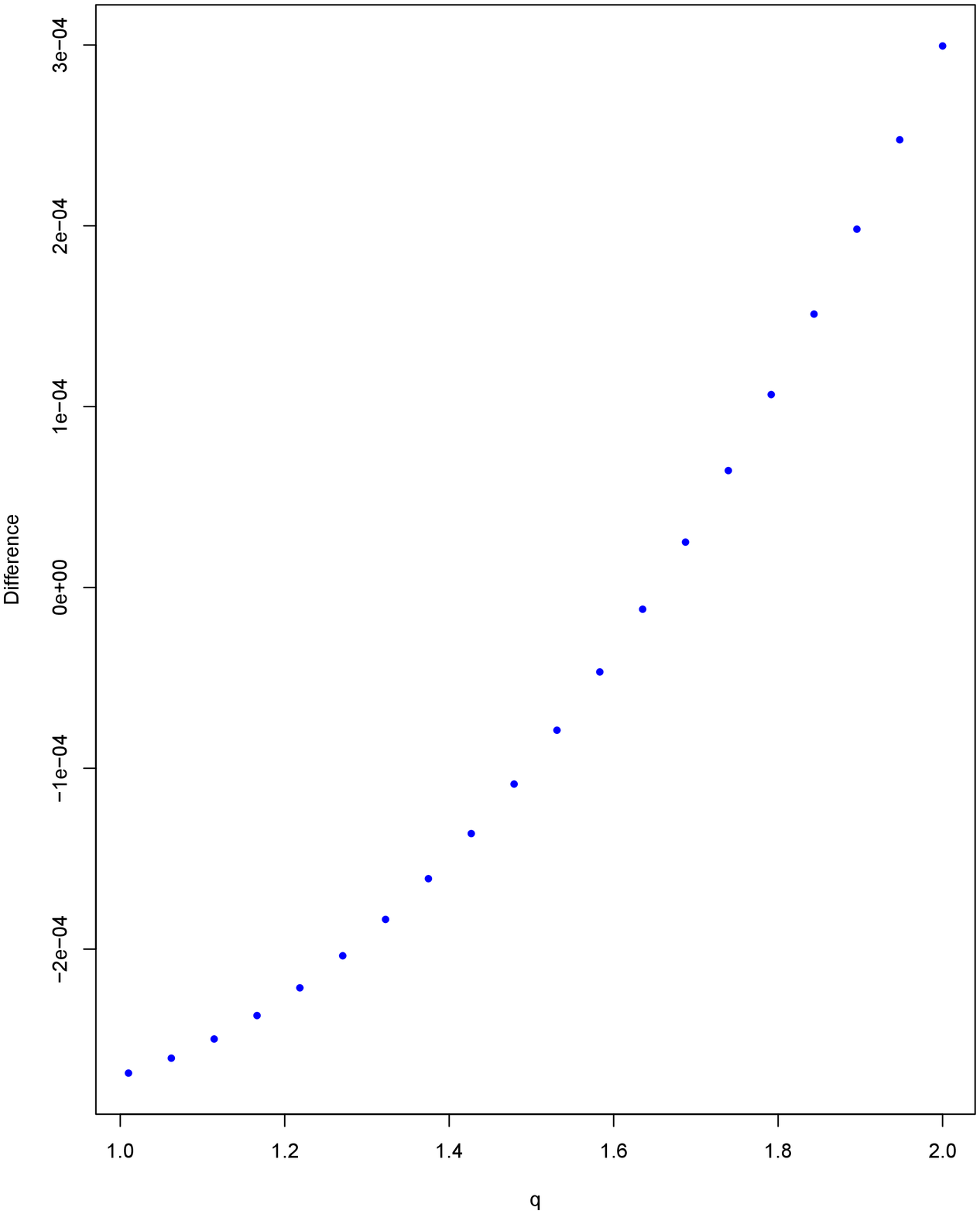}}
    \caption{Differences between the sample R\'enyi function and the fitted model}
    \label{fig18}
\end{figure}

Fig.~\ref{fig16e} demonstrates the fit of the log-normal model(shown in the red colour) to the empirical R\'enyi function. As this plot is rather similar for all other models and windows, we present only the plots of residuals in Fig.~\ref{fig16f}, Fig.~\ref{fig17c}, Fig.~\ref{fig17f} and Fig.~\ref{fig18}.

As the estimated value of $a$ is close to zero, the fit of Model 1 gives an almost degenerated case, when either $\sigma_{Y}^2$ is very small or $b$ is very large, which is consistent with the plot in Fig.~\ref{fig13}. The results in Table~\ref{table11} also confirm that multifractality is very small as for all observation windows $a$ is almost zero and $\alpha_{\max}-\alpha_{\min}$ is very~small.

Next, for the log-gamma model specified by \eqref{eq:11a} we used the reparameterisation $A= \frac{2}{\beta}\ln(b)$, $B={\lambda}^{-1}$ and considered the non-linear model $T(q)-q+1=A^{-1}(\ln(1-Bx)-x\ln(1-B))$. The command \enquote{nlsLM} from the R package minpack.lm with appropriate initial values was used to fit the model to the sample values of $\hat{T}(q)-q+1$. The values of estimated parameters were $\hat{A}=0.029407$ and $\hat{B}=0.005469$ with $RMSE=1.7198 \times 10^{-6}$. The corresponding values of $b$, $\lambda$ and $\beta$ satisfy the assumptions of Theorem~\ref{theo:6.2}.

For the log-negative-inverse-gamma model given by \eqref{eq:12a}, the reparameterisation $A=\frac{1}{2\ln(b)}$, $B=\beta$, $C=\sqrt{\lambda}$ and application of \enquote{nlsLM} resulted in $\hat{A}=0.254719$, $\hat{B}=5.695755$, $\hat{C}=0.386207$ and $RMSE=2.6201\times 10^{-10}$. Note that the obtained $\hat{C}$ corresponds to $\lambda$ that is outside of $L_{\beta,\lambda}$ in Theorem~\ref{theo:6.3}. For parameters satisfying the conditions of Theorem~\ref{theo:6.3}, RMSE is substantially larger. As the results in Theorem~\ref{theo:6.3} might be also true for other parameters(see the discussion in Section~\ref{S1:81}), we used the obtained value of~$\hat{C}$.

 Model~4 was fit by using a linear regression model with the parameter $A=-{\frac{1}{2\log_{2}b}}.$ The estimated parameter $\hat{A}$ was $-0.000667$ and $RMSE = 2.56533 \times 10^{-5}$. For Models~5 and~6 the non-linear regression approach and the R function \enquote{nls} were used. For Model 5 the estimates were found as $\hat{A}=-0.000762$, $\hat{k}\approx 1$ and $RMSE = 1.6347 \times 10^{-5}$. Finally, for Model 6 the estimated parameters were $\hat{A}=0.000269$, $\hat{k}=1$ and $RMSE = 1.8393 \times 10^{-4}$.

Fig.~\ref{fig17c} and Fig.~\ref{fig18} demonstrate that departures of the fitted models from the empirical R\'enyi function are very small, but have different patterns. The numerical studies suggested that Models~2 and~3 are more flexible than the other models. However, to fit these models one has to very carefully choose initial values of the parameters for the nls estimation. Different initial values can lead to different results which can be a potential issue for data which are similarly to CMB show minor multifractality. Also, nls method's rates of convergence for Models~2 and~3 are very slow. Models 1, 4, 5 and 6 have less parameters and are less flexible than Models 2 and 3. However, in many cases they give a reasonable fit very quickly, are robust to the choice of initial values and more computationally efficient.

All models gave a good fit to the empirical R\'enyi functions. The analysis in this section suggests no significant or very small multifractality for the currently available resolution of CMB measurements.\\[-5mm]

\section{Conclusion}
\label{S1:9}

This paper investigates the multifractal behaviour of spherical random fields and some applications to cosmological data from the mission Planck. The aim of this paper is to introduce several multifractal models for random fields on a sphere and to propose simpler models where the R\'enyi function can be computed explicitly. All R\'enyi functions for the specified models exhibit either parabolic or approximately linear behaviours. We present the  R\'enyi function computations for different CMB sky windows located at different places of the sphere. Finally, we fit the specified models to actual CMB data. All models fit to the data. The analysis suggests that there may exist a very minor multifractality of the data.

Some related problems and extensions of the current research that would be interesting for future studies:

\begin{itemize}
    \item Develop statistical tests for different types of R\'enyi functions;
    \item Prove that the theoretical results and the formulae for the R\'enyi functions are also valid for the values of $q$ outside the interval $[1, 2]$, see~\cite{denisov2016limit};
    \item Study other models based on vector random fields (similar to Model 6), where the R\'enyi functions can be computed explicitly;
    \item Develop some approaches to study rates of convergence for the obtained asymptotics, that would serve as analogous of classical convergence rates in central and non-central limit theorems, see~\cite{ANH2015111, anh2019};
    \item Investigate changes of the R\'enyi functions depending on evolutions of random fields driven by SPDEs on the sphere, see \cite{anh2018approximation, broadbridge2019random, broadbridge2020spherically};
    \item Apply the developed models and methodology to other spherical data, in particular, to new high-resolution CMB data from future CMB-S4 surveys~\cite{abazajian2019cmbs4} which will be collecting $3D$ observations.
\end{itemize}
This paper studies data that are modelled as restrictions of $3D$ random fields to the unit sphere. Compared to the available literature, this approach is more consistent with real CMB observations that exist in $\mathbb{R}^3$, but are measured only on $s_2(1),$ see the discussions in \cite{anh2018approximation, broadbridge2019random, broadbridge2020spherically}.  For other applications it would be important to develop similar results for the case of intrinsic spherical random fields, i.e., random fields directly defined on $s_2(1),$ see the discussion about differences of covariance models of random fields with supports in $\mathbb{R}^3$ and $s_2(1)$ in \cite{gneiting2013}.

\section*{Declarations}

{\bf Acknowledgements} This research was partially supported under the Australian Research Council's Discovery Projects funding scheme (project number  DP160101366). The authors are also grateful to the  anonymous referees for their suggestions that helped to improve the style of the paper. We also would like to thank Professor P. Broadbridge for numerous discussions of CMB and underlying physics models.

\noindent{\bf Conflicts of interest/Competing interests}
The authors declare that they have no known competing interests for the results reported in this paper. 

\noindent{\bf Availability of data and material} CMB data from the mission Planck which were used for the analysis purposes are freely available in the NASA/IPAC Infrared Science Archive, see~\cite{[dataset]}.

\noindent{\bf Code availability} All numerical studies were conducted by using Maple 2019.0 and R 3.6.3 software, in particular, the R packages \lq rcosmo\rq~\cite{fryer2019rcosmo, rcosmo1} and \lq RandomFields\rq~\cite{randomfield}. A reproducible version of the code in this paper is available in the folder \enquote{Research materials} from the website~\url{https://sites.google.com/site/olenkoandriy/}. 

\noindent{\bf Authors' contributions} All the authors equally contributed to the paper.


\bibliographystyle{spmpsci}      
\bibliography{Bibliography}   

\begin{appendices}
\section{Proofs}\label{appen}

\begin{proof}[Theorem~{\rm \ref{theo:5.1}}]

By Remark~\ref{muconv}, from the weak convergence of the measures $\mu_k$ to $\mu$ and the assumption of exponential boundedness of the covariance function of the mother field, it follows that
$$E \mu^q_k(B^3) \rightarrow E\mu^q(B^3), \quad k \rightarrow \infty.$$

   By the Lyapunov's inequality, see~\cite[p.162]{Loeve},
$$E \mu^q_k(B^3) \leq (E\mu^2_k(B^3))^{q/2}, \: \text{for} \: q \in [1,2].$$

   Therefore, to guarantee $E\mu^q(B^3)< +\infty, \: q \in [1,2],$ it is sufficient to provide such $b$ and $\sigma^2_{\Lambda}$ that $$\sup_{k\in \mathbb{N}} E\mu^2_k(B^3) < +\infty.$$
By \eqref{eq:4a}, the non-negativity of $\Lambda_{k}(y)$ and independence of $\Lambda^{(i)}$ it holds
\[ E\mu^2_k(B^3) = E \int_{B^3} \int_{B^3} \Lambda_k(y)\Lambda_k(\tilde{y})d\tilde{y}dy
= \int_{B^3} \int_{B^3} E[\Lambda_k(y)\Lambda_k(\tilde{y})]d\tilde{y}dy 
\]\[= \int_{B^3} \int_{B^3} E \prod^k_{i=0} \Lambda^{(i)}(yb^i) \Lambda^{(i)}(\tilde{y}b^i) d\tilde{y}dy
= \int_{B^3} \int_{B^3} \prod^k_{i=0} E\Lambda^{(i)}(yb^i) \Lambda^{(i)}(\tilde{y}b^i) d\tilde{y}dy
\]\[= \int_{B^3} \int_{B^3} \prod^k_{i=0} \biggl(E(\Lambda^{(i)}(yb^{i})-1)(\Lambda^{(i)}(\tilde{y}b^i)-1)+E \Lambda^{(i)}(yb^i) + E \Lambda^{(i)}(\tilde{y}b^i) -1 \biggl) d\tilde{y}dy
\]
{\small  \[= \int_{B^3} \int_{B^3} \prod^k_{i=0}(Cov(\Lambda(yb^i), \Lambda(\tilde{y}b^i))+1)d\tilde{y}dy
= \int_{B^3} \int_{B^3} \prod^k_{i=0}\left(1+ \sigma^2_{\Lambda}\rho_{\Lambda}(\Vert y- \tilde{y}\Vert b^i)\right)d\tilde{y}dy
\]}
\[\leq \int_{B^3} \int_{B^3} \prod^k_{i=0}(1+ \sigma^2_{\Lambda}C e^{-\gamma \Vert y- \tilde{y}\Vert b^i})d\tilde{y}dy
\leq \int_{B^3} \int_{B^3} \prod^{\infty}_{i=0}(1+ \sigma^2_{\Lambda}C  e^{-\gamma \Vert y- \tilde{y}\Vert b^i})d\tilde{y}dy.  
\]
    
   From the inequality $1+a \leq e^a$, it follows that
$$E\mu^2_k(B^3) \leq \int_{B^3} \int_{B^3} \prod^{\infty}_{i=0} e^{\sigma^2_{\Lambda}C e^{-\gamma \Vert y- \tilde{y}\Vert b^i}}d\tilde{y}dy.$$

   Introducing the new variables $z=y$, $\tilde{z} = y-\tilde{y},$ one obtains
$$E\mu^2_k(B^3) \leq \int_{B^3} dz\int_{B^3-B^3} \prod^{\infty}_{i=0} e^{\sigma^2_{\Lambda}C e^{-\gamma \Vert \tilde{z}\Vert b^i}}d\tilde{z},$$
 where $B^3-B^3 = \{\tilde{z}: \tilde{z}=y-\tilde{y}, \: y, \: \tilde{y} \in B^3 \}.$

   Hence, by using the spherical change of variables,
\[E\mu^2_k(B^3) \leq |B^3|\int^{diam(B^3)}_{0} {r}^{2}\prod^{\infty}_{i=0} e^{\sigma^2_{\Lambda}C e^{-\gamma {r} b^i}}dr\]
\[ = \frac{|B^3|}{\gamma^{3}}\int^{\gamma diam(B^3)}_{0} {r}^{2}\prod^{\infty}_{i=0} e^{\sigma^2_{\Lambda}C e^{-r b^i}}d\tau.
\]

   As the exponent $e^{\sigma^2_{\Lambda}C e^{-r b^i}}$ is a decreasing function of $r$, selecting $n(r) = \max(0, -[\log_b(r)]),$ $r>0,$ we obtain
{ \[\prod^{\infty}_{i=0} e^{\sigma^2_{\Lambda}C e^{-r b^i}} \leq \prod^{n(r)-1}_{i=0} e^{\sigma^2_{\Lambda}C e^{-r b^i}}\prod^{\infty}_{i=n(r)} e^{\sigma^2_{\Lambda}C e^{-r b^i}} \]
\[
\leq e^{\sigma^2_{\Lambda}C n(r)}\prod^{\infty}_{i=0} e^{\sigma^2_{\Lambda}C e^{-r b^{i+n(r)}}}
\leq e^{\sigma^2_{\Lambda}C n(r)}\prod^{\infty}_{i=0} e^{\sigma^2_{\Lambda}C e^{-b^i}}.
\]}
   Notice that
{ \[
    \prod^{\infty}_{i=0} e^{\sigma^2_{\Lambda}C e^{-b^i}} = e^{\sigma^2_{\Lambda}C \sum^{\infty}_{i=0} e^{-b^i}} \leq e^{\sigma^2_{\Lambda}C \sum^{\infty}_{i=0} e^{-(1+(b-1)^i)}}\]\[
= e^{\frac{\sigma^2_{\Lambda}C}{e}\sum^{\infty}_{i=0} {e^{-(b-1)}}^i} = e^{{\frac{\sigma^2_{\Lambda}C}{e}} {\frac{1}{1-e^{-(b-1)}}}}  < +\infty. 
\]
}
   Therefore,
\begin{align*}
E\mu^2_k(B^3) &\leq {\frac{|B^3|}{\gamma^{3}}}{e^{\frac{\sigma^2_{\Lambda}C}{e(1-e^{-(b-1)})}}} \int^{\gamma diam(B^3)}_{0} z^{2} e^{\sigma^2_{\Lambda}C n(z)}dz \\ 
&= {\frac{|B^3|}{\gamma^{3}}}{e^{\frac{\sigma^2_{\Lambda}C}{e(1-e^{-(b-1)})}}} \int^{\gamma diam(B^3)}_{0} z^{2}\max\left(1, \: z^{-\frac{\sigma^2_{\Lambda}}{\ln{(b)}}}\right)dz.
\end{align*}
 The integral is finite if $2-\frac{\sigma^2_{\Lambda}C}{\ln{(b)}} > -1,$ \: i.e. $b> e^{\frac{\sigma^2_{\Lambda}C}{3}}.$ \hfill $\Box$
\end{proof}

\begin{proof}[Theorem~{\rm \ref{theo:7.1}}]

 By the definition of  Model~4 it follows that
$$E \Lambda(x) = E(Y^2(x)) = \rho_Y(0) =1, \quad \sigma^2_{\Lambda} = Var \Lambda(x) = E(Y^4(x))-1 =2,$$
$$Cov(\Lambda(x), \Lambda(y)) = E(Y^2(x)-1)(Y^2(y)-1) = 2\rho_{Y}^2(\Vert x-y \Vert).$$

   To compute the covariance, we used the property
\begin{equation}
    E(H_k(Y(x))H_l(Y(y))) = \delta_k^l k! \rho_{Y}^k(\Vert x-y \Vert), \quad x,y \in {\mathbb{R}}^{3}, \label{eq:18a}
\end{equation}
   where $H_k(u), \: k \geq 0, \: u \in \mathbb{R},$ are the Hermite polynomials, see \cite{peccati2011wiener}. For $k=2$, the Hermite polynomial of order $2$ is $H_2(u) = u^2 -1.$

   Thus, Model~4 satisfies Conditions~{\rm\ref{cond1}} and~{\rm\ref{cond2}}.

   Note that the condition $|\rho_{\Lambda}(r)| \leq Ce^{-\gamma r}, \: r>0, \: \gamma>0,$ is equivalent to 
\begin{equation}
    |\rho_{Y}(r)| \leq C^{'}e^{-\gamma^{'}r}, \: r>0, \: \gamma^{'}>0.\label{eq:19a}
\end{equation}

   So, if \eqref{eq:19a} is satisfied, then one can apply Theorems~\ref{theo:4.1} and \ref{theo:4.2} and the R\'enyi function of the limit measure equals to
\[
        T(q) = q-1-\frac{1}{2}\log_b E Y^{2q}(x). 
\]

  Finally, noting that for $p > -1$ and $Z \sim N(\mu, \sigma^2)$
\begin{equation}
     E |Z-\mu|^{p} = \sigma^p \frac{2^{p/2}\Gamma(\frac{p+1}{2})}{\sqrt{\pi}}\label{eq:20a}
\end{equation}
finalises the proof. \hfill $\Box$
\end{proof}

\begin{proof}[Example~{\rm \ref{examp1}}]

 By Remark~\ref{rem:5.1}, it is enough to check that
$$\sup_{k \in N}E \mu_k^4(B^3) = \sup_{k \in N}\int_{B^3}\int_{B^3}\int_{B^3}\int_{B^3} \prod_{i=0}^{k}E\left(\prod_{j=1}^4 Y^2 (y_j b^i)\right)\prod_{j=1}^4 dy_j < +\infty.$$

   Notice, that by Wick's theorem
\begin{equation}
 E\left(\prod_{j=1}^4 Y^2(y_j b^i)\right) = \sum_{p \in P_4^2} \prod_{(j, \tilde{j}) \in p} Cov(Y(y_j b^i), Y(y_{\tilde{j}}b^i)),\label{eq:22a}
\end{equation}
 where the sum is over all parings $p$ of $\{1,1,2,2,3,3,4,4\},$ which are distinct ways of partitioning $\{1,1,2,2,3,3,4,4\}$ into pairs $(i,j)$. The product in \eqref{eq:22a} is over all pairs contained in $p$, see~\cite{janson1997gaussian}.

   Notice that for the pairing $p^{*}=\{(1,1),(2,2),(3,3),(4,4)\}.$
$$\prod_{(j, \tilde{j}) \in p^{*}}Cov(Y(y_j b^i), Y(y_{\tilde{j}}b^i)) = \prod_{j=1}^4 E Y^2(y_j b^i) = 1.$$

   In all other cases of pairing, there is at least one pair $(j, \tilde{j})$ such that $j \neq \tilde{j}$. Therefore, the expectation $E(\prod_{j=1}^4 Y^2(y_j b^i))$ equals $$1+ \sum_{{\substack{p \in P_4^2 \\ p \neq p^{*}}}}\prod_{(j, \tilde{j}) \in p}Cov(Y(y_j b^i), Y(y_{\tilde{j}}b^i)).$$

   As, $1+a < e^a$, it can be estimated by $$\exp\left({\sum_{{\substack{p \in P_4^2 \\ p \neq p^{*}}}}\prod_{(j, \tilde{j})\in p} Cov(Y(y_j b^i), Y(y_{\tilde{j}}b^i))}\right).$$
   As at least for one pairing $(j, \tilde{j}) \in p \neq p^{*}$ it holds that $j \neq \tilde{j},$ then one can use the upper bound$$|Cov(Y(y_j b^i), Y(y_{\tilde{j}}b^i))| \leq \sigma_Y^2 C e^{-\gamma \Vert y_j - y_{\tilde{j}}\Vert b^i},$$
 and the approach from the proof of Theorem~{\rm \ref{theo:5.1}}.

   Namely,

\[
   E\left(\prod_{j=1}^4 Y^2(y_j b^i)\right) \leq e^{\sum_{p \in p_4^2}\prod_{(j, \tilde{j})\in p}\sigma_Y^2 C e^{-\gamma \Vert y_j - y_{\tilde{j}}\Vert b^i}}\]\[ \leq e^{{(\max(\sigma_{\Lambda}^2 C, 1))^4}{\sum_{1 \leq j \leq \tilde{j} \leq 4} e^{-\gamma \Vert y_j - y_{\tilde{j}}\Vert b^i}}}.
\]
 Hence,
\[    \sup_{k \in N}\int_{B^3}\int_{B^3}\int_{B^3}\int_{B^3} \prod_{i=0}^{k}E\left(\prod_{j=1}^4 Y^2 (y_j b^i)\right)\prod_{j=1}^4 dy_j \]\[
     \leq \int_{B^3}\int_{B^3}\int_{B^3}\int_{B^3}  \prod_{i=0}^{\infty}e^{{(\max(\sigma_{\Lambda}^2 C, 1))^4}{\sum_{1 \leq j \leq \tilde{j} \leq 4} e^{-\gamma \Vert y_j - y_{\tilde{j}}\Vert b^i}}} \nonumber \]
     \[ \nonumber
    \leq\left(\int_{B^3}\int_{B^3}\int_{B^3}\int_{B^3} \prod_{i=0}^{\infty}e^{6{(\max(\sigma_{\Lambda}^2 C, 1))^4}{\sum_{1 \leq j \leq \tilde{j} \leq 4} e^{-\gamma \Vert y_j - y_{\tilde{j}}\Vert b^i}}}\prod_{j=1}^4 dy_j\right),
\]
   where the last inequality follows from the generalized  H\"{o}lder's inequality $${\left\Vert\prod_{k=1}^{K}f_k\right\Vert}_{1} \leq \prod_{k=1}^{K} {\Vert f_k \Vert}_{p_k},$$ 
   with $\sum_{k=1}^{K}{p_k}^{-1}=1$. In our case $K=6$ is the number of different $j$ and $\tilde{j}$ satisfying $1 \leq j \leq \tilde{j} \leq 4.$

   Finally, similar to the proof of Theorem~{\rm \ref{theo:5.1}}, from equation  \eqref{eq:21a} we obtain the condition {$b > e^{{\frac{\sigma{(\max(\sigma_{\Lambda}C, 1))}^4}{3}}}.$}
   \hfill $\Box$
\end{proof}

\begin{proof}[Theorem~{\rm \ref{theo:7.2}}]

It follows from \eqref{eq:20a} that
$$E \Lambda(x) = E Y^{2k}(x) = \sigma^{2k}{\frac{2^k \Gamma(k+\frac{1}{2})}{\sqrt{\pi}}} =1,$$
\[
    \sigma^2_{\Lambda} = Var \Lambda(x) = E(Y^{4k}) -1 = {\left(\frac{\sqrt{\pi}}{2^k \Gamma(k+\frac{1}{2})}\right)}^2 {\frac{2^{2k}\Gamma(2k+\frac{1}{2})}{\pi}}-1\]
    \[ = \frac{\sqrt{\pi}\Gamma(2k+ \frac{1}{2})}{{\Gamma}^2(k+\frac{1}{2})} - 1 < +\infty.
\]

   To compute the covariance function we use \eqref{eq:18a} and the following Hermite expansion \cite[page 775]{abramowitz1948handbook}
\[z^{2k} = (2k)! \sum_{i=0}^{k} \frac{H_{2k-2i}(z)}{2^i i!(2k-2i)!}.\]
   Therefore,
\[
    Cov(\Lambda(x), \Lambda(y)) = E(Y^{2k}(x)-1)(Y^{2k}(y)-1) = \frac{\pi E(\tilde{Y}^{2k}(x)\tilde{Y}^{2k}(y))}{2^{2k}{\Gamma}^2(k+\frac{1}{2})}-1 \]\[ = ((2k)!)^2\frac{\pi}{2^{2k}{\Gamma}^2(k+\frac{1}{2})}\sum_{i=0}^k \frac{E[H_{2k-2i}(\tilde{Y}(x))H_{2k-2i}(\tilde{Y}(y))]}{2^{2i}(i!)^2((2k-2i)!)^2}-1 \] 
\begin{equation}\label{eq:23a} = ((2k)!)^2\frac{\pi}{2^{2k}{\Gamma}^2(k+\frac{1}{2})}\sum_{i=0}^k \frac{\tilde{\rho}^{2k-2i}(\Vert x-y \Vert)}{2^{2i}(i!)^2(2k-2i)!}-1, 
\end{equation}
   where $\tilde{Y}(x) = Y(x)/\left(\frac{\sqrt{\pi}}{2^k \Gamma(k+\frac{1}{2})}\right)^{1/{2k}}$ is a zero-mean unit variance Gaussian HIRF with the covariance function
$$\tilde{\rho}(\Vert x-y \Vert) = \left (\frac{2^k \Gamma(k+\frac{1}{2})}{\sqrt{\pi}}\right)^{1/k}\rho_Y(\Vert x-y \Vert).$$

   Notice, that for $i=k$ in \eqref{eq:23a} by the Legendre duplication formula
\begin{align*}
   \frac{((2k)!)^2 \pi}{2^{2k} {\Gamma}^2(k+\frac{1}{2})2^{2k}(k!)^2} 
&=  \frac{{\Gamma}^2(2k+1)\pi}{2^{4k}{\Gamma}^2(k+\frac{1}{2})k^2{\Gamma}^2(k)} 
= \frac{(2k)^2{\Gamma}^2(2k)\pi}{2^{4k}k^2 2^{2-4k}\pi{\Gamma}^2(2k)} = 1.
\end{align*}

   Hence,
$$Cov(\Lambda(x), \Lambda(y))= \frac{((2k)!)^2\pi}{2^{2k}{\Gamma}^2(k+\frac{1}{2})}\sum_{i=0}^{k-1}\frac{(2^k \Gamma(k+\frac{1}{2}))^{2+\frac{2i}{k}}}{2^{2i}(i!)^2(2k-2i)!\pi^{1/2k}}\tilde{\rho}^{2k-2i}(\Vert x-y \Vert).$$

   Therefore, if $|\tilde{\rho}(r)| \leq C^{'}e^{-\gamma^{'}r}, \: r>0,\: \gamma^{'}>0,$ then  the covariance function of Model 5 satisfies the condition $|\rho_{\Lambda}(r)| \leq Ce^{-\gamma r}, r>0, \: \gamma>0,$ and the R\'enyi function equals
\[
   T(q) = q-1-\frac{1}{2}\log_{b} E Y^{2kq}(x)
    = q-1-\frac{1}{2}\log_{b}\left(\frac{2^{kq}\Gamma(kq+\frac{1}{2})}{\sqrt{\pi}}\right). \]\hfill $\Box$
   \end{proof}

\begin{proof}[Theorem~{\rm \ref{theo:7.3}}]  
 
      By properties of the chi-square distribution, it follows that
   $$E \Lambda(x) = \frac{2}{k}E Y(x) =1, \quad Var \Lambda(x) = \frac{4}{k^2}Var Y(x) = \frac{2}{k} < +\infty,$$
   $$Cov(\Lambda(x), \Lambda(y)) = \frac{4}{k^2}\rho_Y(\Vert x-y \Vert).$$
   
      Notice that if $Y(x) = \frac{1}{2}(Z_1^2(x)+...+Z_k^2(x)), \: x \in {\mathbb{R}}^3,$ where $Z_i(x),$ $i=1,...,k,$ are independent zero-mean unit variance components of $k$-dimensional vector Gaussian HIRF with a covariance function $\rho_Z(r), \: r \geq 0$ of each component,~then
   $$Cov(\Lambda(x), \Lambda(y)) = \frac{4}{k^2}\cdot \frac{k}{2}\rho_Z^2(\Vert x-y \Vert) = \frac{2}{k}\rho_Z^2(\Vert x-y \Vert).$$ 
   Therefore, Model 6 satisfies Conditions~{\rm\ref{cond1}} and~{\rm\ref{cond2}} and $|\rho_\Lambda(r)| \leq Ce^{-\gamma r}, \: r>0$, $\gamma >0$,
   if $|\rho_{Y}(r)| \leq C'e^{-\gamma^{'}r}$ or $|\rho_{Z}(r)| \leq C'e^{-\gamma^{'}r}, \: r \geq 0, \gamma^{'}>0.$ 
   
      Then, the corresponding R\'enyi function is given by \eqref{eq:25a}. \hfill $\Box$
  
      \end{proof}
   \end{appendices}
\end{document}